\makeatletter
\newif\ifsiamcls %If SIAM class available
\IfFileExists{siamart250211.cls}{\siamclstrue}{\siamclsfalse}
\ifsiamcls % true case
  \documentclass[final,onefignum,onetabnum]{siamart250211}
\else % false case, not available
  \documentclass[11pt]{article}
  \usepackage[margin=1in]{geometry}
  \usepackage{authblk}
  \usepackage[utf8]{inputenc}
 \usepackage[T1]{fontenc}
\usepackage{lmodern}
  \newcommand{\headers}[2]{}%
  \newenvironment{keywords}{\paragraph{Keywords.}}{}%
  \newenvironment{AMS}{\paragraph{AMS subject classifications.}}{}%
  \newcommand{\email}[1]{\href{mailto:#1}{#1}}%
\fi % END IF SIAM class available
\makeatother

\usepackage{url}
\usepackage[hidelinks]{hyperref}
\usepackage{graphicx}
\usepackage{algorithm}
\usepackage{algpseudocode}

\usepackage{amsmath,amssymb}
\usepackage{mathtools}
\usepackage{bm}
\usepackage{booktabs}
\usepackage{microtype}

\usepackage{xcolor}
\usepackage[normalem]{ulem}
\newcommand{\MSadd}[1]{#1}%      added text kept, rendered in black (arXiv view)
\newcommand{\MSdel}[1]{}%        deleted text suppressed
\newcommand{\MSnote}[1]{}%       reviewer annotations suppressed
\newif\ifarxivmode
\arxivmodetrue
\ifarxivmode
  \newcommand{\Ksadd}[2]{#2}
  \newcommand{\SLadd}[2]{#2}
  \newcommand{\SZadd}[2]{#2}
  \newcommand{\RTadd}[2]{#2}
\else
  \definecolor{slteal}{RGB}{0,128,128}
  \definecolor{szorange}{RGB}{200,90,0}
  \definecolor{rtpurple}{RGB}{128,0,128}
  \newcommand{\Ksadd}[2]{{\color{blue}#2\textsuperscript{[K#1]}}}
  \newcommand{\SLadd}[2]{{\color{slteal}#2\textsuperscript{[SL#1]}}}
  \newcommand{\SZadd}[2]{{\color{szorange}#2\textsuperscript{[SZ#1]}}}
  \newcommand{\RTadd}[2]{{\color{rtpurple}#2\textsuperscript{[RT#1]}}}
\fi
\newcommand{\SLdel}[1]{}%
\newcommand{\SLnote}[1]{}%
\newcommand{\SZdel}[1]{}%
\newcommand{\SZnote}[1]{}%
\newcommand{\RTdel}[1]{}%
\newcommand{\RTnote}[1]{}%

\ifdefined\diffmode
  \definecolor{r98zgreen}{RGB}{0,128,0}
  \newcommand{\Rzzedit}[2]{{\color{r98zgreen}#2}}
\else
  \newcommand{\Rzzedit}[2]{#2}
\fi
\usepackage{tikz}
\usetikzlibrary{arrows.meta,positioning,calc,fit,shapes.geometric,%
  decorations.pathreplacing,patterns,matrix}
\DeclareUnicodeCharacter{2019}{'}
\DeclareUnicodeCharacter{2014}{---}
\DeclareUnicodeCharacter{2013}{--}

\headers{PMP with Smoothed Hamiltonian and Adaptive Refinement}{Ksous et al.}
\title{Pontryagin-Based Solver with Smoothed Hamiltonian,
  Adaptive $\Delta t$, and PA-Bundle Refinement}

\newcommand{\samethanks}[1][\value{footnote}]{\protect\footnotemark[#1]}
\author{%
  \MSdel{Selim}\MSadd{Salim} Ksous\MSadd{\thanks{\MSadd{Mathematics,
    Computer, Electrical and Mathematical Sciences and Engineering
    (CEMSE) Division, King Abdullah University of Science and
    Technology (KAUST), Thuwal, Saudi Arabia
    (\email{salim.ksous@kaust.edu.sa},
    \email{sebastianleonardo.lalvaysegovia@kaust.edu.sa},
    \email{erik.vonschwerin@kaust.edu.sa},
    \email{raul.tempone@kaust.edu.sa}).}}}%
  \and Sebastian Lalvay Segovia\MSadd{\samethanks[1]}%
  \and Mattias Sandberg\MSadd{\thanks{\MSadd{Department of
    Mathematics, KTH Royal Institute of Technology, Stockholm,
    Sweden (\email{msandb@kth.se}, \email{szepessy@kth.se}).}}}%
  \and Erik von Schwerin\MSadd{\samethanks[1]}%
  \and Anders Szepessy\MSadd{\samethanks[2]}%
  \and Ra\'ul Tempone\MSadd{\samethanks[1]}%
}

\makeatletter
\@ifundefined{newsiamremark}{%
  \newtheorem{remark}{Remark}%
}{%
  \newsiamremark{remark}{Remark}%
}
\makeatother
\begin{document}
\ifsiamcls\maketitle\else\date{}\maketitle\fi

\Rzzedit{S1}{\renewcommand{\thefootnote}{\fnsymbol{footnote}}%
\footnotetext[8]{This manuscript was prepared with the assistance
of Anthropic Claude, OpenAI ChatGPT and Grammarly as editing and
structural-review tools; all mathematical content, numerical
experiments, and final wording decisions are the responsibility of
the authors.}%
\renewcommand{\thefootnote}{\arabic{footnote}}}

\begin{abstract}

%% === inlined from abstract_pmp.tex (R98z-cc arxiv-flat) ===
% ------------------------------------------------------------------
% Paper 0 abstract — shared by all top-level entry points
% Synced and refactored in R98p (von Schwerin proofreading round).
% ------------------------------------------------------------------
We present a Pontryagin-based numerical solver for deterministic
optimal control problems in Bolza form. The solver regularizes the
(generally nonsmooth) Hamiltonian by applying a log-sum-exp (soft-min)
smoothing to a concave piecewise-affine (PA) bundle surrogate,
producing a smoothed Hamiltonian~$H_\delta$ that is $C^\infty$ and
concave in the costate and continuously differentiable (indeed
$C^{1,1}$) in the state under the standing assumptions. The resulting
two-point boundary value problem (TPBVP) from the Pontryagin minimum
principle is discretized by a symplectic Euler scheme. \Rzzedit{S2}{The smoothed
discrete TPBVP is solved by damped Newton iteration on the full-space
all-at-once nonlinear system in the discrete state and costate
variables\Rzzedit{B2-S2}{ (the primary formulation used in the numerical experiments)}; a classical multiple-shooting realization is also discussed
as an alternative implementation in~\S5.} A unified adaptive outer
loop jointly controls three error sources: the time-step distribution
$\Delta t$, the number of planes $M$ in the PA surrogate, and the
smoothing parameter $\delta$. The a~posteriori time-step error
indicators follow the error-density framework of Karlsson, Larsson,
Sandberg, and Tempone (2015). The closest related methods are
classical Pontryagin-minimum-principle shooting,
Hamilton--Jacobi--Bellman / viscosity PDE solvers, occupation-measure
relaxations, max-plus approximations, and direct collocation. Four
unconstrained-state numerical examples---the hypersensitive scalar
problem of Karlsson, Larsson, Sandberg, and Tempone (2015), a scalar
nonsmooth Hamiltonian benchmark, a singular tracking problem, and an
attitude control allocation benchmark whose Hamiltonian minimum is
supplied by a quadratic-programming oracle---illustrate the intended
workflow and demonstrate the roles of time refinement, PA-bundle
enrichment, and smoothing continuation. A smooth linear--quadratic
regulator with closed-form Riccati reference solution is retained as
a calibration consistency check in Appendix~A.

We adopt a value-first scope throughout. The goal is to approximate
the value function at the initial state, and all three a~posteriori
error indicators (time, bundle, and smoothing) target computable
diagnostics of the resulting value error; the discrete state, costate,
and control returned by the symplectic-Euler Newton solve are
computational means rather than separately certified outputs. After
the outer loop converges, an admissible trajectory is recovered post
hoc by plugging the converged costate into the Hamiltonian-minimizing
control selector and integrating the induced state dynamics forward.
%% === end of abstract_pmp.tex ===

\end{abstract}



\section{Introduction}\label{sec:introduction}

%% === inlined from introduction_pmp.tex (R98z-cc arxiv-flat) ===
% ------------------------------------------------------------------
% §1 Introduction
% ------------------------------------------------------------------

\MSdel{\paragraph{Program map}
This work develops an adaptive PMP-shooting solver and serves as the
companion to a multi-manuscript deterministic-optimal-control program.
A forthcoming \emph{value-first costate-induced} foundation companion
proves a value-equivalence theorem and supplies the selector-well-posed
domain on which a second \emph{adaptive a posteriori} companion adds a
balanced refinement controller; a third \emph{state-constrained}
companion extends the costate framework to explicit state constraints
and compares its direct $(P2_K)$ solver to the present PMP-shooting
solver; a fourth
\emph{two-level} companion restores routing delays via an outer
Lagrangian-relaxation layer. This work is independent of those
companions in its theory and numerics; the cross-work relation is
made explicit in Section~\ref{sec:discussion}.}%
\MSdel{Program map paragraph removed at R10 per user directive: this work is now presented standalone, without forward references to companion works in the introduction.}

\paragraph{Contribution}
This work presents a Pontryagin-minimum-principle (PMP) shooting
solver for the (generally nonsmooth) Pontryagin two-point
boundary-value problem (TPBVP) \Ksadd{1}{arising from continuous-time
deterministic optimal control in Bolza form}, based on a log-sum-exp
smoothing of a concave piecewise-affine Hamiltonian surrogate, with
adaptive control of the time step $\Delta t$, the PA bundle size $M$,
and the smoothing parameter $\delta$, validated on four
unconstrained-state benchmarks.

\MSadd{\paragraph{Value-first scope}
The goal of this work is to approximate the value function $U(0,x_0)$.\ All
a~posteriori error indicators and all adaptive-refinement decisions target the value-error
$|J_h-U(0,x_0)|$ (or computable proxies of it); the controller refines only those layers whose
contribution to the value-error exceeds tolerance.\ We do \emph{not} claim separate pathwise
error bounds on the discrete state $x_h$, costate \RTadd{R62-M1}{$p_h$}, or control $\alpha_h$ produced by
the symplectic-Euler Newton solve --- these are computational means, not certified outputs.\
\RTadd{R59-M7}{After the outer loop converges with costate \RTadd{R62-M1}{$p_h$}, an admissible
\emph{admissible-trajectory} pair $(x_h^{\mathrm{adm}},\alpha_h)$ is recovered post hoc by a
single forward sweep at the converged \RTadd{R62-M1}{$p_h$} using the Hamiltonian-minimizing control
selector;\ the three-step procedure (selector definition, closed-loop integration, recorded
control) is stated in detail in~Section~\ref{ssec:outer-loop} (post-convergence recovery paragraph) to keep
the introduction light.}\ The reconstructed pair is admissible by construction and serves as
the admissible trajectory for the \RTadd{R72-B1}{primal upper bound $J^{\mathrm{pr}}(\alpha_h)\ge U(0,x_0)$};\
we do not claim it as a separately-certified approximation of any optimal trajectory.\ The
present work does \emph{not} certify a primal--dual value gap:\ lower verification functions
and certified intervals are developed in the costate-optimization companion works and are
reserved for that work.}

Many optimal control problems in engineering, robotics, and
computational physics feature nonlinear dynamics, constraints, and
nonsmooth control profiles, making their numerical solution
challenging. The Pontryagin \MSdel{Maximum}\MSadd{minimum} \RTadd{R62-m2}{principle}
(PMP)~\cite{Pontryagin1962}\MSadd{\ (also known as the maximum
principle; we use the minimum convention throughout this work
to match the $\min$~sign in~\eqref{eq:hamiltonian})} provides a
powerful necessary-condition framework, but numerical solvers based
on PMP typically require differentiability of the Hamiltonian for
Newton-type convergence---which fails when the optimal control
switches discontinuously (bang--bang regimes). Structure-preserving time
integration (e.g., symplectic Runge--Kutta
schemes~\cite{Chyba2009,SanzSerna2016}) improves stability over long
horizons but does not, by itself, resolve the nonsmoothness
difficulty. Direct transcription
methods~\cite{Betts1998,Elnagar1995,Rao2014} avoid this issue by
discretizing the problem into a nonlinear program, at the expense of
structure preservation and---for fine meshes---large sparse systems.

This work proposes a solver for the PMP two-point boundary value
problem that combines three ingredients:
\begin{itemize}
\item \textbf{Soft-min Hamiltonian smoothing.}\;
 A piecewise-affine (PA) bundle surrogate of the concave Hamiltonian
 is smoothed via a log-sum-exp (soft-min)
 formula~\cite{Nesterov2005}, producing a smoothed Hamiltonian
 $H_\delta$ that is $C^\infty$ and concave in the costate (and
 sufficiently differentiable in the state under the standing
 assumptions) to enable Newton-based solvers.
\item \textbf{PA bundle approximation.}\;
 The Hamiltonian is accessed through a \emph{value-and-subgradient
 oracle}; the PA surrogate is built incrementally from oracle
 responses and refined adaptively when the modeling error exceeds a
 tolerance.
\item \textbf{Three-source adaptive refinement.}\;
 A unified outer loop controls the time mesh $\Delta t$, the bundle
 size~$M$, and the smoothing level~$\delta$ via computable
 a~posteriori error indicators. \SZadd{AN3}{The solver has a two-level
 structure:\ an \emph{outer} adaptive loop that revises $\Delta t$, $M$
 and~$\delta$ between solves, and an \emph{inner} damped-Newton iteration
 that solves the discrete TPBVP at frozen bundle and smoothing.}\;
 The time-step adaptivity follows the
 error-density framework
 \SZadd{AR4}{of}~\cite{KarlssonLarssonSandbergTempone2015}.
 \MSadd{The three indicators play asymmetric roles in the outer
 loop:\ the time indicator is a \emph{max}-type stopping quantity
 driving mesh refinement, while the bundle and smoothing indicators
 are \emph{sum}-type quantities driving bundle enrichment and
 $\delta$-halving;\ the outer loop includes a time-balance guard so
 that further mesh refinement is suppressed once the time indicator
 is an order of magnitude below the dominant non-time indicator
 (Section~\ref{ssec:outer-loop}).}
\end{itemize}
\Rzzedit{S3}{The resulting solver uses a symplectic Euler discretization
of the smoothed canonical system, solved by damped Newton iteration on
the full-space all-at-once nonlinear system; a classical multiple-shooting
realization is also discussed in Section~\ref{sec:discussion} as an
alternative implementation. The same outer loop adapts all three error
sources toward user-prescribed tolerances.}

\paragraph{Relation to prior work}
The method occupies a middle ground between classical indirect
shooting~\cite{Bryson1975} (efficient but sensitive to
nonsmoothness), Hamilton--Jacobi--Bellman PDE
solvers~\cite{FlemingSoner2006} (globally optimal but limited by the
curse of dimensionality), occupation-measure
relaxations~\cite{Lasserre2008} and max-plus
methods~\cite{Akian2008} (global but computationally heavy), and
direct collocation~\cite{Betts1998,Elnagar1995,Rao2014} (robust
constraint handling via large sparse KKT systems, but with loss of
Hamiltonian structure). By combining smoothing, adaptivity, and
structure preservation within an indirect PMP framework, the solver
aims to retain the efficiency of shooting methods while mitigating
their sensitivity to nonsmooth switches. For a modern survey of
trajectory optimization approaches, see Rao~\cite{Rao2014}.

\paragraph{Scope and contribution}
The contribution of this work is algorithmic rather than
theorem-driven. We combine (i)~a concave piecewise-affine surrogate
of the minimization Hamiltonian, (ii)~log-sum-exp smoothing of that
surrogate, and (iii)~adaptive control of the time step~$\Delta t$,
the bundle size~$M$, and the smoothing level~$\delta$ within an
indirect, symplectic-Euler solver. The solver architecture is close
in spirit to several existing lines of work: classical multiple
shooting~\cite{BockPlitt1984} and the symplectic-PMP / a~posteriori
adaptivity lineage
\SZadd{AR4}{of~\cite{SandbergSzepessy2006, KarlssonLarssonSandbergTempone2015}}. We do not
claim to have first implemented any single one of these ingredients;
rather, we claim that their specific combination into a three-source
adaptive indirect solver, in the form presented here, is new to the
best of our knowledge. The mathematical content is an
\emph{algorithmic note}: we present the regularized discrete solver
and its three error controls.
\MSdel{State-constrained extensions---the viability-based tangential Hamiltonian, the relaxed discrete-step variant, the log-barrier interior-point smoothing, and the associated constrained value-equivalence theorem---are the subject of a forthcoming costate-optimization companion work and are not treated here.}
We do not claim new convergence rates in $\Delta t$, $M$, or
$\delta$ for the coupled adaptive scheme; the a~posteriori indicators
of Section~\ref{ssec:indicators} are used as computable surrogates
that align with the error decomposition
\SZadd{AR4}{of~\cite{KarlssonLarssonSandbergTempone2015, SandbergSzepessy2006}}, and the adaptive
refinement strategy is validated empirically on the four
unconstrained-state benchmarks of Section~\ref{sec:numerics}.

\paragraph{Structure}
Section~\ref{sec:formulation} introduces the problem class, standing
assumptions, and the hierarchy of Hamiltonian objects ($H$, $\bar H$,
$H_\delta$). Section~\ref{sec:algorithm} presents the discrete
solver: the symplectic Euler TPBVP, the Newton iteration, the three
a~posteriori indicators, and the adaptive outer loop.
\RTadd{R62-B1}{Section~\ref{sec:numerics} provides numerical
experiments on four unconstrained-state benchmark problems:\ the
hypersensitive scalar problem, a nonsmooth scalar Hamiltonian
benchmark, a singular tracking problem, and a QP-oracle attitude
control-allocation benchmark.\ Appendix~\ref{apx:lqr} records a
smooth LQR/Riccati calibration check.}
Section~\ref{sec:discussion} discusses connections to prior work and
limitations, and Section~\ref{sec:conclusion} gives concluding remarks.
%% === end of introduction_pmp.tex ===

\section{Problem Formulation, Surrogate, and Smoothing}\label{sec:formulation}

%% === inlined from Pformulation_v2.tex (R98z-cc arxiv-flat) ===
% ------------------------------------------------------------------
% §2 Problem Formulation, Surrogate, and Smoothing
% ------------------------------------------------------------------

We consider deterministic optimal control problems in Bolza form over
a finite time horizon~$[0,T]$:
\begin{equation}\label{eq:bolza}
\begin{aligned}
\min_{a(\cdot)} &\quad g(x(T)) + \int_0^T \ell(x(t),a(t),t)\,dt \\
\text{s.t.} &\quad \dot x(t) = f(x(t),a(t),t),\qquad x(0)=x_0, \\
 &\quad a(t)\in\mathcal{A}\subset\mathbb{R}^m \quad\text{a.e.\ }t\in[0,T].
\end{aligned}
\end{equation}
Here $x(t)\in\mathbb{R}^d$ is the state, $a(t)$ the control, $f$ the
controlled drift, $\ell$ the running cost, and $g$ the terminal cost.
\MSadd{The state dimension is $d\in\mathbb{N}$ and the control
dimension is $m\in\mathbb{N}$, so $\mathcal{A}\subset\mathbb{R}^m$.}
The \MSadd{admissible} control set is denoted~$\mathcal{A}$; see
Assumption~(A2) below.
\Rzzedit{B2-S33}{For scalar box constraints we write
\begin{equation*}
\Pi_{[a,b]}(z) \;:=\; \min\bigl(b,\max(a,z)\bigr)
 \;=\; \operatorname{median}(a,z,b)
\end{equation*}
for the scalar box projection of $z\in\mathbb{R}$ onto $[a,b]\subset\mathbb{R}$;
this notation appears in the numerical examples of Section~\ref{ssec:ex-hypersensitive}
and the LQR calibration of Appendix~\ref{apx:lqr}.}

\MSdel{Scope (state constraints) paragraph removed at R10 per user directive: this work is presented as a standalone unconstrained-state PMP-shooting solver, without forward references to a state-constrained companion in the body.}

\subsection{Standing assumptions}\label{ssec:assumptions}

\begin{enumerate}
\item[\MSadd{(A0)}] \MSadd{\emph{Existence and normality.}\;
 An optimal pair $(x^\star,a^\star)$ exists for the Bolza
 problem~\eqref{eq:bolza}, and the normal Pontryagin minimum
 principle applies to it (i.e.\ the abnormal multiplier is~$1$,
 in the sense of~\cite[Ch.~22]{Vinter2000}, and the costate
 equation~\eqref{eq:pmp} is satisfied).
 Standard sufficient conditions for~(A0)
 \MSdel{include:\
 $\mathcal{A}$ is compact and the integrand is continuous in~$a$;
 $f$ is control-affine with convex velocity set;\ or
 $\mathcal{A}$ is closed and unbounded with $\ell$ coercive in~$a$
 uniformly on compact $(x,t)$-sets, so a Filippov-type argument
 applies. We do not redo the corresponding existence theory and
 refer to standard texts such as~\cite{Vinter2000,Clarke2013}.}%
 \MSadd{are
 (\cite[Theorem~7.4.5]{CannarsaSinestrari2004}):
 \begin{enumerate}
 \item the control set $\mathcal{A}$ is compact;
 \item there exists $K>0$ such that
 $|f(x_2,a,t)-f(x_1,a,t)|\le K\,|x_2-x_1|$ for all
 $x_1,x_2\in\mathbb{R}^d$, $a\in\mathcal{A}$, $t\in[0,T]$;
 \item for every $R>0$ there exists $\beta_R>0$ such that
 $|\ell(x_2,a,t)-\ell(x_1,a,t)|\le\beta_R\,|x_2-x_1|$ for all
 $x_1,x_2\in\overline{B_R}$, $a\in\mathcal{A}$, $t\in[0,T]$;
 \item for every $(x,t)\in\mathbb{R}^d\times[0,T]$ the set
 \[
 \bigl\{(c,v)\in\mathbb{R}^{d+1}:
 \exists\,a\in\mathcal{A},\ v=f(x,a,t),\
 c\ge\ell(x,a,t)\bigr\}
 \]
 is convex (Cesari-type augmented-set convexity; the slack variable $c$ here is a
 running-cost upper bound and is not the costate);
 \item for every $x\in\mathbb{R}^d$ and $a\in\mathcal{A}$, the
 maps $t\mapsto f(x,a,t)$ and $t\mapsto\ell(x,a,t)$ are
 continuous on~$[0,T]$;
 \item $g\in C(\mathbb{R}^d)$.
 \end{enumerate}
 We do not redo the corresponding existence theory and refer
 to~\cite{CannarsaSinestrari2004,Vinter2000,Clarke2013} for
 the standard development.}}
\item[(A1)] \emph{Regularity (local).}\;
 \MSadd{Let $\mathcal{D}\subset\mathbb{R}^d\times[0,T]$ be a
 compact region containing the iterates of the discrete solver
 and the optimal trajectory of~(A0).}\;
 The drift $f:\mathbb{R}^d\times\mathcal{A}\times[0,T]\to\mathbb{R}^d$ and
 the running cost $\ell:\mathbb{R}^d\times\mathcal{A}\times[0,T]\to\mathbb{R}$
 are continuous in $(x,a)$ for a.e.\ $t$ and measurable in~$t$ for
 every $(x,a)$. Moreover, $f$ and $\ell$ are continuously
 differentiable in~$x$ for each $(a,t)$, with derivatives that are
 uniformly bounded and Lipschitz in~$x$ \MSadd{on
 $\mathcal{D}\times\mathcal{A}$ (when $\mathcal{A}$ is unbounded
 these uniformity conditions are required only on compact
 $a$-sub\-sets that contain~$a^\star$ and the iterates'
 controls)}.\;The terminal cost
 $g:\mathbb{R}^d\to\mathbb{R}$ satisfies $g\in C^1(\mathbb{R}^d)$.\MSadd{\
 \SZadd{AM7}{Admissible controls are taken in $L^\infty([0,T];\mathcal{A})$
 -- in particular distributional / Dirac-measure controls
 (such as $a(t)=\delta(t-t_0)$) are excluded -- so the Carath\'eodory
 framework below applies.}\;
 In particular, (A1) implies via the Carath\'eodory existence
 theorem that, for every measurable admissible control
 $a:[0,T]\to\mathcal{A}$, the ODE
 $\dot x(t)=f(x(t),a(t),t)$, $x(0)=x_0$, has a unique maximal
 solution; for trajectories that remain in the compactness
 region~$\mathcal{D}$, the maximal interval contains~$[0,T]$, which
 is the situation for the optimal trajectory of~(A0) and for the
 solver iterates by construction.}
\item[(A2)] \MSadd{\emph{Hamiltonian minimizer attainment.}\;
 For every $(p,x,t)$ in the region of interest, the map
 \[
 a \;\mapsto\; p\cdot f(x,a,t)+\ell(x,a,t)
 \]
 attains its minimum over~$\mathcal{A}$. A sufficient condition
 is that $\mathcal{A}$ be compact and the integrand continuous
 in~$a$, or that $\mathcal{A}$ be closed and the integrand coercive
 in~$a$ uniformly on compact $(p,x,t)$-sets. This is a pointwise
 condition needed to define the Hamiltonian and the oracle;\
 problem-level existence of an optimal control belongs to~(A0).}\;
 \MSadd{(The LQR calibration run of Appendix~\ref{apx:lqr} is formally
 stated with $\mathcal{A}=\mathbb{R}$ and uses algorithmic box bounds
 in the solver run; see~\S\ref{apx:lqr} for details.)}
\item[\MSadd{(A3)}] \RTadd{R66-B1}{\emph{Attainment and active-branch
 regularity.}\;
 For every queried $(p,x,t)$, the Hamiltonian minimum
 $H(p,x,t)=\min_{a\in\mathcal{A}}\{p\cdot f(x,a,t)+\ell(x,a,t)\}$
 is attained. On open regions of $(p,x,t)$ where a minimizing branch
 $a^\star(p,x,t)$ is unique and locally stable (i.e.\ depends
 continuously on $(p,x,t)$), the branch-reduced Hamiltonian
 \[
 H^{\mathrm{br}}(p,x,t)
 \;=\;
 p\cdot f(x,a^\star(p,x,t),t)+\ell(x,a^\star(p,x,t),t)
 \]
 is differentiable in~$x$ and satisfies the envelope identity \Rzzedit{B1-S6}{(see \cite[Thm.~4.5.4]{Vinter2000})}
 \[
 \nabla_x H^{\mathrm{br}}(p,x,t)
 \;=\;
 p\cdot\nabla_x f(x,a^\star,t)+\nabla_x\ell(x,a^\star,t).
 \]
 At switching points where the active branch changes, derivatives are
 read in the Clarke/generalized-gradient or quasi-Newton sense
 \cite{Clarke2013,Vinter2000};\ \Rzzedit{A1e}{the Newton system of \S\ref{ssec:newton} uses}
 the smoothed Hamiltonian $H_\delta$ of~\S\ref{ssec:smoothing}, whose
 derivatives are classical for fixed $\delta>0$.\ This
 active-branch reading subsumes both the classically-smooth and the
 bang--bang/singular regimes:\ in particular, the bundle inequality
 $\bar H\ge H$ of~\S\ref{ssec:PA-surrogate} is a structural
 consequence of the PA construction and does \emph{not} require
 uniqueness or strict convexity of the minimizer.

 \emph{Remark (sufficient condition).}\;
 Strict convexity of $a\mapsto p\cdot f(x,a,t)+\ell(x,a,t)$ on
 $\mathcal{A}$ for every $(p,x,t)$ is a sufficient condition for
 active-branch uniqueness and is the standard textbook hypothesis.\
 It is \emph{not} assumed in the bang--bang and singular examples
 of~\S\ref{ssec:ex-nonsmooth}--\ref{ssec:ex-singular}, where strict
 convexity fails---either identically on $\mathcal{A}$ for
 affine-in-$a$ Hamiltonians (classical bang--bang) or off the
 switching set $\{(p,x,t):p\cdot\partial_a f+\partial_a\ell=0\}$
 (singular regimes);\ active-branch regularity in the sense above
 remains in force.}
\item[\MSadd{(A4)}] \RTadd{R70-B2}{\emph{Differentiability of the
 smoothed Hamiltonian.}\;
 For each fixed $\delta>0$, the soft-min Hamiltonian $H_\delta$
 of~\eqref{eq:H_delta} is $C^\infty$ in the costate $p$.\ Its
 regularity in $x$ is inherited from the coefficient functions
 $g_i(x,t)=f(x,a_i^\star,t)$ and $d_i(x,t)=\ell(x,a_i^\star,t)$
 of~\S\ref{ssec:PA-surrogate}:\ under (A1) it is $C^{1,1}$ in $x$,
 and under (A5) it is $C^2$ in $(p,x)$ on the region used by the
 Newton solve.\ The gradients
 in~\eqref{eq:grad_H_delta} are classical under (A1), while the
 full Newton Jacobian of Section~\ref{ssec:newton} uses the
 strengthened regularity (A5).\ On the switching set of the
 unsmoothed problem, the Newton step is interpreted as a smoothed
 quasi-Newton step in the sense of~(A5).}
\item[(A5)] \emph{Higher regularity for the Newton analysis.}\;
 Whenever the exact Newton Jacobian of the discrete system of
 Section~\ref{ssec:newton} is invoked, we also assume that
 $f$ and~$\ell$ are twice continuously differentiable in~$x$ on the
 relevant region, with uniformly bounded second derivatives, and
 that $g\in C^2(\mathbb{R}^d)$. Under (A5) the smoothed
 Hamiltonian $H_\delta$ is $C^2$ in~$(p,x)$ and the Jacobian blocks
 $\partial_x\nabla_pH_\delta$, $\partial_x\nabla_xH_\delta$ exist in
 the classical sense. Under (A1) alone, the Jacobian is
 interpreted as a generalized or quasi-Newton Jacobian.
\end{enumerate}

\subsection{Hamiltonian and the
 \MSdel{Pontryagin maximum principle}%
 \MSadd{Pontryagin minimum principle}}
\label{ssec:pmp}

\MSadd{%
\emph{Naming convention.}\; Because the Hamiltonian
$H(p,x,t)=\min_{a\in\mathcal{A}}\{p\cdot f+\ell\}$
is defined as a minimum over controls, the resulting necessary
conditions take the form of a \emph{minimum} principle; we refer to
it throughout this work as the \emph{Pontryagin minimum principle}
(PMP). This matches the sign convention used in
every formula that follows; readers more familiar with the ``maximum
principle'' form should recover it by flipping the sign of~$H$.}%
\MSnote{MS, p.~3: rename ``maximum'' $\to$ ``minimum'' principle to
match the $\min$ in the definition of~$H$.}

Define the minimization Hamiltonian
\begin{equation}\label{eq:hamiltonian}
H(p,x,t) \;:=\; \min_{a\in\mathcal{A}}\bigl\{p\cdot f(x,a,t)+\ell(x,a,t)\bigr\},
\end{equation}
where $p\in\mathbb{R}^d$ is the costate (adjoint) variable.
\MSdel{Under the minimization convention adopted here, the map}%
\MSadd{The map}
$p\mapsto H(p,x,t)$ is \emph{concave} (as a pointwise minimum of
affine functions of~$p$).

Under (A0) and assumptions \MSdel{(A1)--(A3)}\MSdel{(A1) and (A3)}\MSadd{(A1) and (A2)}, the
Pontryagin \MSdel{maximum} \MSadd{minimum} principle (PMP) provides
the following necessary conditions for an optimal
trajectory~$x^*$: there exists a costate
arc~$p(\cdot)$ satisfying, at points of differentiability of~$H$,
the \SZadd{AM5}{\emph{canonical system}\;\Rzzedit{S5}{(Hamilton's
canonical equations)}}
\begin{equation}\label{eq:pmp}
\begin{aligned}
\dot x(t) &= \phantom{-}\nabla_p H(p(t),x(t),t), \\
\dot p(t) &= -\nabla_x H(p(t),x(t),t),
\end{aligned}
\qquad
\begin{aligned}
x(0) &= x_0, \\
p(T) &= \nabla g(x(T)),
\end{aligned}
\end{equation}
together with the minimality condition
$H(p(t),x(t),t) = p(t)\cdot f(x(t),a^*(t),t)+\ell(x(t),a^*(t),t)$
for a.e.\ $t$.
\Rzzedit{B2-S19}{We write
\begin{equation}\label{eq:astar-selector}
a^\star(p,x,t) \;\in\;
\mathop{\mathrm{arg\,min}}_{a\in\mathcal{A}}\bigl\{p\cdot f(x,a,t)+\ell(x,a,t)\bigr\}
\end{equation}
for any measurable selection of the pointwise Hamiltonian minimizer; under~(A2) the
argmin set is nonempty for each $(p,x,t)$, and under~(A3) the
selection is unique and continuous on the active-branch regions
of~\eqref{eq:astar-selector}. At switching points the selection is
read in the Clarke generalized sense~\cite{Clarke2013,Vinter2000};\
the discrete solver uses the smoothed
Hamiltonian~$H_\delta$ of~\S\ref{ssec:smoothing}, whose gradients are
classical for fixed $\delta>0$.}
In general the Hamiltonian~$H$ defined
by~\eqref{eq:hamiltonian} need not be classically differentiable
\SLdel{R98m1}{in~$x$}\SLadd{R98m1}{in $(p,x)$} (e.g., at switching times where the minimizer $a^*(p,x,t)$ is
non-unique). At such points, the \SLdel{R98m1}{costate equation}\SLadd{R98m1}{canonical system}
in~\eqref{eq:pmp} is to be interpreted as the \SLdel{R98m1}{differential inclusion}\SLadd{R98m1}{differential inclusions} \SLadd{R98m1}{$\dot x(t)\in\partial_p^C H(p(t),x(t),t)$,
}$\dot p(t)\in-\partial_x^C H(p(t),x(t),t)$, where $\partial^C$
denotes Clarke's generalized gradient~\cite{Clarke2013}. Under~(A1)
\Rzzedit{S6}{the envelope identity of (A3) yields the classical
$x$-derivative} whenever the minimizer $a^*(p,x,t)$ is unique. For numerical
purposes we replace~\eqref{eq:pmp} by the same canonical system
driven by the smoothed Hamiltonian~$H_\delta$ of
Section~\ref{ssec:smoothing}, for which classical derivatives exist
everywhere.

\begin{figure}[t]
 \centering
 \includegraphics[width=\linewidth]{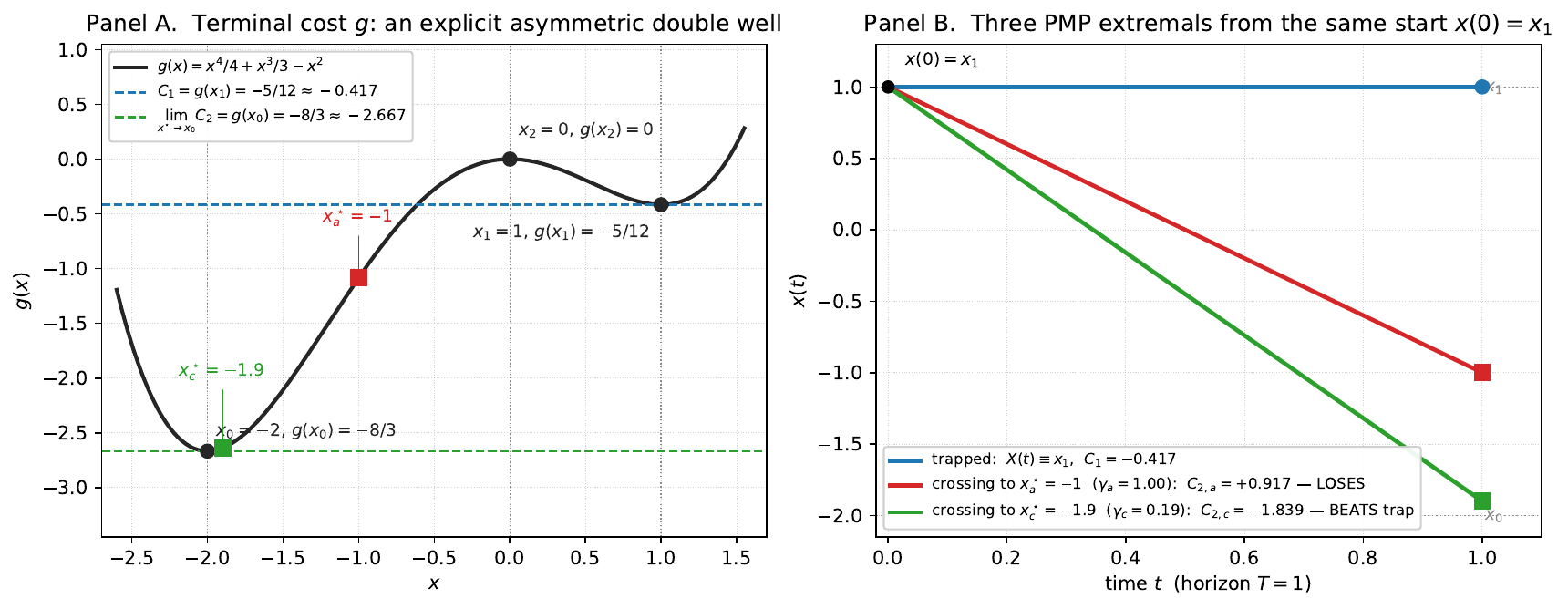}
 \caption{\RTadd{R50}{Pontryagin local-minimum trap for the explicit
 asymmetric double well $g(x)=x^4/4+x^3/3-x^2$ of
 Remark~\ref{rmk:pmp-trap}, with horizon $T=1$.
 \textbf{Panel~A:} the cost~$g$ as a function of the terminal
 state~$x$, with critical points $x_0=-2$ (global minimum,
 $g(x_0)=-8/3$), $x_2=0$ (local maximum, $g(x_2)=0$), and $x_1=1$
 (local-but-not-global minimum, $g(x_1)=-5/12$). The blue dashed
 line marks the trapped cost $C_1=g(x_1)$; the green dashed line
 marks the limiting crossing cost $\lim_{x^\star\to x_0}C_2=g(x_0)$.
 The two coloured squares show the candidate crossing-extremal
 targets used in Panel~B.
 \textbf{Panel~B:} the three Pontryagin extremals starting from
 \RTadd{R69-C6}{$x(0)=x_1=1$}. The trapped extremal (blue) sits at \RTadd{R69-C6}{$x(t)\equiv 1$} for
 all~$\gamma$, with cost $C_1=-5/12$. The shallow crossing extremal
 (red) reaches $x^\star_a=-1$ at $t=T$ with the matched
 $\gamma_a=1$ but pays a strictly larger cost
 $C_{2,a}=+11/12>C_1$, so the trap wins. The deep crossing extremal
 (green) reaches $x^\star_c=-1.9$ at $t=T$ with the matched
 $\gamma_c\approx 0.190$ and pays $C_{2,c}\approx-1.84\ll C_1$, so
 the global well wins. All numbers are evaluated in closed form
 from~\eqref{eq:pmp-trap-g}--\eqref{eq:pmp-trap-c2}; no numerical
 solver is involved.}}
 \label{fig:pmp-trap}
\end{figure}

% R50: warning remark on PMP local-minimum trap (Anders' suggestion;
% adapted from a teaching example in the RWTH NM-SDE lecture notes).
\begin{remark}[\RTadd{R50}{PMP: necessary, not sufficient---a double-well warning}]\label{rmk:pmp-trap}
\Rzzedit{B2-S7}{The purpose of this remark is to exhibit a concrete
one-dimensional example in which the shooting PMP system has multiple
roots, demonstrating that solving the PMP system to high accuracy
does not by itself guarantee global optimality and motivating the
warm-start procedures of~\S\ref{ssec:init}.}\
The PMP system~\eqref{eq:pmp} provides \emph{necessary} conditions for
optimality; solutions of it need not be optimal. When the terminal cost~$g$
is non-convex, the shooting \Rzzedit{S8}{two-point boundary value
problem} can admit multiple extremals, and a
Newton-type iteration warm-started near a suboptimal one will converge to
it. We illustrate this with an explicit one-dimensional construction in
which all of the standing assumptions~(A0)--(A5) of~\S\ref{ssec:assumptions}
hold.

Fix $T>0$ and a control-penalty weight $\gamma>0$ (chosen distinct from
the Hamiltonian-smoothing parameter $\delta$ of~\S\ref{ssec:smoothing}),
and consider the deterministic \Rzzedit{S8}{optimal control problem} with \RTadd{R69-C6}{scalar state
$x(t)\in\mathbb{R}$, scalar control $\alpha(t)\in\mathbb{R}$, dynamics
$\dot x(t)=\alpha(t)$ with initial condition $x(0)=x_1$ (specified
below), and Bolza cost}
\begin{equation}\label{eq:pmp-trap-cost}
C(\alpha) \;=\; \tfrac{\gamma}{2}\!\int_0^T\!\alpha(t)^2\,dt + g(\RTadd{R69-C6}{x(T)}).
\end{equation}
Take the explicit asymmetric double well
\begin{equation}\label{eq:pmp-trap-g}
g(x) \;=\; \frac{x^4}{4}+\frac{x^3}{3}-x^2,
\qquad
g'(x) \;=\; x(x+2)(x-1).
\end{equation}
The critical points of~$g$ are
\[
x_0=-2,\qquad x_2=0,\qquad x_1=1,
\]
with values
$g(x_0)=-\tfrac{8}{3}$, $g(x_2)=0$, $g(x_1)=-\tfrac{5}{12}$, so that
$x_0$ is the global minimum, $x_1$ is a local-but-not-global minimum,
$x_2$ is the local maximum separating the two wells, and
$g(x_0)<g(x_1)$. The Hessian check
$g''(x)=3x^2+2x-2$ gives $g''(x_0)=6>0$, $g''(x_2)=-2<0$,
$g''(x_1)=3>0$, confirming the well structure.

\SLdel{R98m2}{The Hamiltonian
$H(p,x,\alpha)=p\alpha+\tfrac{\gamma}{2}\alpha^2$}\SLadd{R98m2}{The pre-minimization integrand
$\widetilde H(p,x,\alpha):=p\alpha+\tfrac{\gamma}{2}\alpha^2$} is strictly convex
in~$\alpha$, so the pointwise minimizer is unique:
$\alpha^*(p)=-p/\gamma$, giving the \SLdel{R98m2}{reduced }Hamiltonian
$H(p,x)\SLadd{R98m2}{:=\min_\alpha \widetilde H(p,x,\alpha)}=-p^2/(2\gamma)$. The costate equation $\dot p=-\partial_x H=0$
forces $p\equiv\mathrm{const}$, hence \RTadd{R69-C6}{$x(t)=x_1-pt/\gamma=1-pt/\gamma$} is
affine. The terminal condition $p(T)=g'(\RTadd{R69-C6}{x(T)})$ collapses the BVP to the
scalar fixed-point equation
\begin{equation}\label{eq:pmp-trap-fixed-point}
p \;=\; g'\!\bigl(1-pT/\gamma\bigr),
\end{equation}
which admits at least two solutions:
\begin{enumerate}
\item \emph{Trapped extremal.}\;
 $p=0$, $\alpha\equiv 0$, \RTadd{R69-C6}{$x(t)\equiv 1$}, with cost
 $C_1=g(1)=-\tfrac{5}{12}$. Admissible because $g'(1)=0$.
\item \emph{Crossing extremal.}\;
 Pick any $x^\star\in(x_0,x_2)=(-2,0)$, where $g'(x^\star)>0$. Set
 $p=g'(x^\star)$ and $\gamma=g'(x^\star)\,T/(1-x^\star)$. The affine
 trajectory \RTadd{R69-C6}{$x(t)=1-(1-x^\star)\,t/T$} then lands at~$x^\star$ at time~$T$,
 with cost
 \begin{equation}\label{eq:pmp-trap-c2}
 C_2 \;=\; g(x^\star)+g'(x^\star)\,\frac{1-x^\star}{2}.
 \end{equation}
\end{enumerate}

\emph{Cost comparison.}\;
As $x^\star\to x_0=-2$ with the matching~$\gamma\to 0$, smoothness of~$g$
gives $g(x^\star)\to g(x_0)=-\tfrac{8}{3}$ and $g'(x^\star)\to 0$, hence
\[
C_2 \;\longrightarrow\; g(x_0) \;=\; -\tfrac{8}{3}
\;<\; -\tfrac{5}{12} \;=\; g(x_1) \;=\; C_1,
\]
a cost gap of $C_1-C_2\to \tfrac{9}{4}$. By continuity, the crossing
extremal beats the trapped one for every~$x^\star$ in a sufficiently small
neighborhood of~$x_0$, but not for every $x^\star\in(x_0,x_2)$:
Figure~\ref{fig:pmp-trap} illustrates the dichotomy with two
explicit choices, a shallow target $x^\star_a=-1$ for which the
crossing extremal loses ($C_{2,a}=+11/12>-5/12=C_1$, with the
matched~$\gamma_a=1$) and a deep target $x^\star_c=-1.9$ for which
it wins ($C_{2,c}\approx-1.84<-5/12$, with the matched
$\gamma_c\approx 0.19$).

\emph{Reading.}\; Two morals follow. First, solving the PMP shooting
system to high accuracy does not by itself guarantee global optimality;
comparing several extremals (multi-start, homotopy, or a coarse global
benchmark) is the companion ingredient required to discriminate among
them. Second, the trap here lives in the BVP shooting
equation~\eqref{eq:pmp-trap-fixed-point} having multiple roots, \emph{not}
in the inner Hamiltonian minimization, which is uniquely solved by strict
convexity of~$H$ in~$\alpha$. Consequently, neither the PA-softmin
smoothing of~\S\ref{ssec:smoothing} nor the time-adaptive refinement of
the outer loop addresses this failure mode: they target Hamiltonian
non-convexity in~$\alpha$ and time-discretization error, respectively.
Escaping a Pontryagin local minimum is an orthogonal capability that lies
outside the scope of this work.
\end{remark}

\subsection{Piecewise-affine bundle surrogate $\bar H$}
\label{ssec:PA-surrogate}

\Ksadd{2}{The construction in this subsection is independent of
pointwise uniqueness of the Hamiltonian minimizer $a^*(p,x,t)$: each
oracle response contributes one affine majorant of $p\mapsto H(p,x,t)$
at the sampled $(x,t)$, regardless of whether the minimizer attaining
that majorant is unique. \RTadd{R70-M1}{The active-branch regularity
assumption~(A3) is invoked only where the envelope derivative is
needed;\ the structural inequality $\bar H\ge H$ below does not
require uniqueness, branch stability, or strict convexity of the
minimizer.}}

Because $p\mapsto H(p,x,t)$ is concave, it is the infimum of all
affine functions lying above it. We maintain a finite collection of
affine \emph{majorants}
$\{L_i(p,x,t)=g_i(x,t)\cdot p+d_i(x,t)\}_{i=1}^M$, sampled at
previously probed costate values, and define the piecewise-affine (PA)
surrogate
\begin{equation}\label{eq:pa_model}
\bar H(p,x,t) \;:=\;
 \min_{1\le i\le M}\bigl\{g_i(x,t)\cdot p + d_i(x,t)\bigr\}.
\end{equation}
By construction $\bar H(p,x,t)\ge H(p,x,t)$ pointwise, and adding
more majorants can only tighten the bound: $H\le\bar H_{M+1}\le\bar H_M$.

\emph{Oracle assumptions.}\label{ssec:oracle}\;
We assume access to a \emph{Hamiltonian oracle} that, given a query
point $(p,x,t)$, returns:
(i)~the value $H(p,x,t)$,
(ii)~a supergradient $g\in\partial_p H(\cdot,x,t)|_{p}$
(equivalently, the \Rzzedit{S9}{controlled drift} $f(x,a^*,t)$ at an optimal control $a^*$),
and (iii)~the optimal control $a^*$ itself.
Each oracle call produces one new affine majorant
$L_{M+1}(p')=g\cdot p'+[H(p,x,t)-g\cdot p]$ at the sample
point~$(x,t)$; the extension of this majorant to a neighborhood of
$(x,t)$ is discussed in Section~\ref{ssec:extension}. The PA
surrogate is refined adaptively whenever the modeling error
$\bar H-H$ sampled along the current trajectory exceeds a
tolerance~$\varepsilon_{\mathrm{PA}}$.

For a one-dimensional concave function on a compact interval under
Lipschitz regularity, the uniform approximation error of a PA bundle
with $M$~planes is $\mathcal{O}(1/M)$ (see~\cite{HiriartUrrutyLemarechal1993},
Theorem~IX.2.2.1, applied to the convex function~$-H$).
\MSdel{No multidimensional rate is claimed here without further analysis.}%
\MSnote{MS, p.~4: strike this defensive sentence.}

\subsection{Extension of the surrogate coefficients to all $(x,t)$}
\label{ssec:extension}

The oracle assumptions above produces, at each query point
$(p_i^*,x_i^*,t_i^*)$, a slope vector $g_i^*\in\mathbb{R}^d$ and an
intercept $d_i^*\in\mathbb{R}$ that define a single affine majorant
of $p\mapsto H(p,x_i^*,t_i^*)$ \emph{at that particular $(x,t)$
value}. The Newton solver, however, needs to evaluate $\bar H$,
$H_\delta$, and---critically---$\nabla_xH_\delta$ at \emph{every}
mesh point $(x_n,t_n)$, which changes from one Newton iteration to
the next. Treating $g_i$ and $d_i$ as functions of $(x,t)$ defined
on a neighborhood of the current trajectory is therefore an
essential and nontrivial step. In particular, the extension must
satisfy two requirements:
\begin{enumerate}
\item[(E1)] \emph{Majorant preservation.}\;
 $L_i(p,x,t)=g_i(x,t)\cdot p+d_i(x,t)\ge H(p,x,t)$ for all
 $(p,x,t)$ in the relevant region (at least along the current and
 nearby iterates), so that $\bar H\ge H$ remains valid.
\item[(E2)] \emph{Smoothness in $(x,t)$.}\;
 The maps $x\mapsto g_i(x,t)$ and $x\mapsto d_i(x,t)$ must be
 differentiable, since $\nabla_xH_\delta$ (needed for the costate
 equation and the Newton Jacobian) involves
 $\nabla_xg_i$ and $\nabla_xd_i$.
\end{enumerate}
Note that (E1) and (E2) are automatically compatible with concavity
of $p\mapsto\bar H(p,x,t)$: a minimum of affine functions of~$p$ is
concave in~$p$ regardless of how the coefficients depend on~$(x,t)$.

\RTadd{R51}{Throughout this work we assume the problem data $f(x,a,t)$ and
$\ell(x,a,t)$ are available analytically (not only through an oracle),
which is the setting of all four benchmarks in~\S\ref{sec:numerics}. In
this setting the natural and cleanest construction is the
\emph{control-based} (frozen-control) extension defined below; it
satisfies~(E1) and~(E2) by inspection and gives a globally valid
majorant $\bar H\ge H$. The QP-oracle attitude control allocation of
\S\ref{ssec:ex-B} is a special case within this setting:\ the
oracle's evaluation involves an internal optimisation (a parametric
QP), but $f$ and $\ell$ are still available analytically, so the
control-based extension applies directly. Other settings---black-box
stochastic oracles, learned Hamiltonians, mean-field reductions, and
other cases where a clean $(f,\ell)$ decomposition is not
available---admit local-model surrogates with diagnostic-only error
indicators, and are deferred to future work
(\S\ref{sec:discussion}).}

\smallskip
Record the \emph{optimal control} $a_i^*$ returned by the oracle at
$(p_i^*,x_i^*,t_i^*)$ and define the $i$-th majorant as
\begin{equation}\label{eq:control_extension}
L_i(p,x,t) \;:=\; p\cdot f(x,a_i^*,t)+\ell(x,a_i^*,t)
\qquad\text{for all $(p,x,t)$}.
\end{equation}
That is, $g_i(x,t)=f(x,a_i^*,t)$ and $d_i(x,t)=\ell(x,a_i^*,t)$ are
obtained by \emph{freezing the control} at~$a_i^*$ and letting
$(x,t)$ vary.

\emph{Majorant preservation (E1):}\; By definition of the Hamiltonian,
$H(p,x,t)=\min_{a\in\mathcal{A}}\{p\cdot f(x,a,t)+\ell(x,a,t)\}$, so
evaluating the integrand at the fixed control $a_i^*\in\mathcal{A}$
always gives an upper bound:
$L_i(p,x,t)\ge H(p,x,t)$ for \emph{all} $(p,x,t)$. The majorant
property holds globally, not merely at the sample point.

\emph{Smoothness (E2):}\; Since $f$ and $\ell$ are $C^1$ in $x$
by~(A1), the coefficients inherit the same regularity, and the
required derivatives are
\begin{equation}\label{eq:coeff_derivs}
\nabla_x g_i(x,t)=\nabla_x f(x,a_i^*,t),\qquad
\nabla_x d_i(x,t)=\nabla_x \ell(x,a_i^*,t).
\end{equation}
These are computable from the problem data at negligible extra cost.

\emph{Concavity in $p$:}\; Each $L_i$ is affine in~$p$ for every
$(x,t)$, so $\bar H=\min_i L_i$ is concave in~$p$ automatically.

\subsection{Smoothed Hamiltonian $H_\delta$}\label{ssec:smoothing}

The PA surrogate~$\bar H$ is piecewise affine and hence nonsmooth.
To obtain a differentiable model suitable for Newton-type methods, we
apply a \emph{soft-min} (log-sum-exp) smoothing with parameter
$\delta>0$:
\begin{equation}\label{eq:H_delta}
H_\delta(p,x,t) \;:=\;
-\delta\log\!\left(\sum_{i=1}^M
 \exp\!\left(-\frac{g_i(x,t)\cdot p+d_i(x,t)}{\delta}\right)\right).
\end{equation}
For each fixed~$(x,t)$, the map $p\mapsto H_\delta(p,x,t)$ is concave
and $C^\infty$ in~$p$. Under assumption~(A1), the coefficient
functions $g_i(x,t)=f(x,a_i^*,t)$ and $d_i(x,t)=\ell(x,a_i^*,t)$ are
$C^1$ in~$x$ with Lipschitz derivatives, so the smoothed Hamiltonian
is $C^{1,1}$ in~$x$. Higher regularity in~$x$ is inherited only if
$f$ and~$\ell$ are correspondingly smoother in~$x$ (we make this
explicit in assumption~(A5) of Section~\ref{ssec:assumptions} when
the Newton Jacobian is analyzed). Via the coefficient
extension~\eqref{eq:control_extension} of
Section~\ref{ssec:extension}, the required derivatives are
\begin{equation}\label{eq:grad_H_delta}
\nabla_p H_\delta = \sum_{i=1}^M w_i\,g_i(x,t),\qquad
\nabla_x H_\delta = \sum_{i=1}^M w_i\bigl[\nabla_x g_i(x,t)\cdot p
 + \nabla_x d_i(x,t)\bigr],
\end{equation}
where the softmax weights are
\begin{equation}\label{eq:softmax_weights}
w_i \;=\;
\frac{\exp\!\bigl(-L_i(p,x,t)/\delta\bigr)}
 {\sum_{j=1}^M \exp\!\bigl(-L_j(p,x,t)/\delta\bigr)},
\qquad L_i(p,x,t)=g_i(x,t)\cdot p+d_i(x,t).
\end{equation}
The uniform smoothing bias satisfies~\cite{Nesterov2005}
\begin{equation}\label{eq:smoothing_bias}
0 \;\le\; \bar H(p,x,t) - H_\delta(p,x,t) \;\le\; \delta\log M.
\end{equation}
\RTadd{R62-M2}{For fixed bundle size~$M$, the soft-min bias
in~\eqref{eq:smoothing_bias} is $\mathcal{O}(\delta)$ as
$\delta\to 0$;\ when $M$ grows during bundle refinement, the smoothing
budget should be read through the product $\delta\log M$. The
bound~\eqref{eq:smoothing_bias} controls only the smoothing bias of
$H_\delta$ \emph{relative to the current PA bundle $\bar H$};\ the
residual PA modeling error $\bar H - H$ is a separate quantity,
monitored by the indicator $\eta_{\mathrm{PA}}$
of~\S\ref{ssec:indicators}, and is not absorbed
by~\eqref{eq:smoothing_bias}.}

\SLadd{N1}{%
\begin{remark}[Structural $1/\delta$ factor in the Hessian of $H_\delta$
 and Newton conditioning]\label{rem:H_delta_hessian}
Differentiating $\nabla_p H_\delta=\sum_i w_i\,g_i$
in~\eqref{eq:grad_H_delta} once more in~$p$ and using
$\partial_p w_i = -w_i(g_i-\sum_j w_j g_j)/\delta$ gives the closed-form
identity
\Rzzedit{B3-S11}{\begin{equation}\label{eq:H_delta_hessian}
\begin{split}
\nabla^2_{pp} H_\delta(p,x,t)
&\;=\; -\frac{1}{\delta}\,\mathrm{Cov}_w\!\bigl[g_i(x,t)\bigr] \\
&\;=\; -\frac{1}{\delta}\left[
\sum_{i=1}^M w_i\,g_i g_i^{\!\top}
\;-\;\Bigl(\sum_{i=1}^M w_i g_i\Bigr)
 \Bigl(\sum_{i=1}^M w_i g_i\Bigr)^{\!\top}
\right].
\end{split}
\end{equation}}
The Hessian therefore carries a structural~$1/\delta$ factor weighted
by the (positive semidefinite) covariance of the bundle directions
under the softmax weights. When one affine branch dominates ($w_i\to 1$
for a single $i$) the covariance is small and the curvature stays
bounded; near switching regions, however, two or more branches are
simultaneously active, the covariance is order one, and the curvature
of~$H_\delta$ grows like $\mathcal{O}(\delta^{-1})$.
Consequently, the statement that $H_\delta$ is $C^2$ in $(p,x)$
under~(A5) is understood for each fixed $\delta>0$; the conditioning
of the Newton Jacobian of Section~\ref{ssec:newton} deteriorates
proportionally to $\delta^{-1}$ as $\delta\downarrow 0$, which
motivates the $\delta$-annealing schedule of
Section~\ref{ssec:outer-loop} and the $\eta_\delta$ indicator that
tracks the bound~\eqref{eq:smoothing_bias} along iterations
(cf.~\cite[Example~3.2, p.~A959]{KarlssonLarssonSandbergTempone2015}).
\end{remark}%
}

\MSadd{%
Figure~\ref{fig:H_surrogate} illustrates the three Hamiltonian
objects $H$, $\bar H$, and $H_\delta$ on a one-dimensional slice at
fixed $(x, t)$; the technical reading is given in the caption.%
}

\MSadd{%
\begin{figure}[htbp]
\centering
\begin{tikzpicture}[xscale=1.8, yscale=1.4, >=Stealth,
 font=\footnotesize]
 % Axes
 \draw[->, thick] (-0.1, 0) -- (5.7, 0) node[below right] {$p$};
 \draw[->, thick] (0, -0.1) -- (0, 3.3) node[above left] {};
 \node[rotate=90, anchor=south] at (-0.35, 1.6) {value};
 % Three affine majorants (dashed grey)
 \draw[gray!70, thick, dashed] (0, 0.5) -- (5.3, 3.15);
 \draw[gray!70, thick, dashed] (0, 1.5) -- (5.3, 1.5);
 \draw[gray!70, thick, dashed] (0, 3.0) -- (5.3, 0.35);
 \node[gray!80] at (5.6, 3.05) {$\ell_1$};
 \node[gray!80] at (5.6, 1.5) {$\ell_2$};
 \node[gray!80] at (5.6, 0.45) {$\ell_3$};
 % PA envelope bar H (blue, piecewise linear): min of the three
 \draw[blue!70!black, very thick]
 (0, 0.5) -- (2, 1.5) -- (3, 1.5) -- (5, 0.5);
 \node[blue!70!black, anchor=south west] at (0.3, 0.7) {$\bar H(p)$};
 % Smoothed H_delta (orange): the exact log-sum-exp soft-min of the
 % three affine pieces, sampled at 21 equispaced points and drawn
 % smoothly. By construction this curve is C^infty and concave.
 % Parameters used: delta = 0.15, M = 3,
 % L_1 = 0.5 + 0.5 p, L_2 = 1.5, L_3 = 3 - 0.5 p.
 \draw[orange!95!black, very thick]
 plot[smooth, tension=0.55] coordinates {
 (0.00, 0.500) (0.25, 0.625) (0.50, 0.749) (0.75, 0.873)
 (1.00, 0.995) (1.25, 1.113) (1.50, 1.224) (1.75, 1.321)
 (2.00, 1.393) (2.25, 1.438) (2.50, 1.452) (2.75, 1.438)
 (3.00, 1.393) (3.25, 1.321) (3.50, 1.224) (3.75, 1.113)
 (4.00, 0.995) (4.25, 0.873) (4.50, 0.749) (4.75, 0.625)
 (5.00, 0.500)
 };
 \node[orange!95!black, anchor=north] at (2.5, 1.42) {$H_\delta(p)$};
 % True concave H (red smooth parabola-like)
 \draw[red!70!black, very thick, domain=0:5, smooth, samples=80]
 plot (\x, {-0.15*(\x-2.5)*(\x-2.5) + 1.28});
 \node[red!70!black, anchor=north east] at (4.6, 0.76) {$H(p)$};
 % Smoothing-bias indicator (double arrow between bar H and H_delta)
 % at p=2 where the gap is largest (1.5 - 1.393 = 0.107, within
 % delta log M = 0.15 * log 3 ~ 0.165)
 \draw[<->, thin] (2.0, 1.393) -- (2.0, 1.5)
 node[midway, right, font=\tiny] {$\le\delta\log M$};
 % Tick marks on x-axis at the kinks
 \draw[thin] (2, 0) -- (2, -0.06);
 \draw[thin] (3, 0) -- (3, -0.06);
\end{tikzpicture}
\caption{One-dimensional schematic of the Hamiltonian objects at
fixed $(x,t)$: three affine majorants $\ell_i(p)=g_i\!\cdot\!p+d_i$
(dashed grey) whose pointwise minimum is the concave piecewise-affine
bundle surrogate~$\bar H$ (blue) and satisfies
$\bar H\ge H$ for the exact concave Hamiltonian~$H$ (red).
\RTadd{R72-M1}{The soft-min smoothing $H_\delta$ (orange) is
$C^\infty$ and concave in~$p$, with $x$-regularity inherited from the
coefficient functions, and satisfies $0\le\bar H-H_\delta\le\delta\log M$
uniformly~\eqref{eq:smoothing_bias}.}}
\label{fig:H_surrogate}
\end{figure}
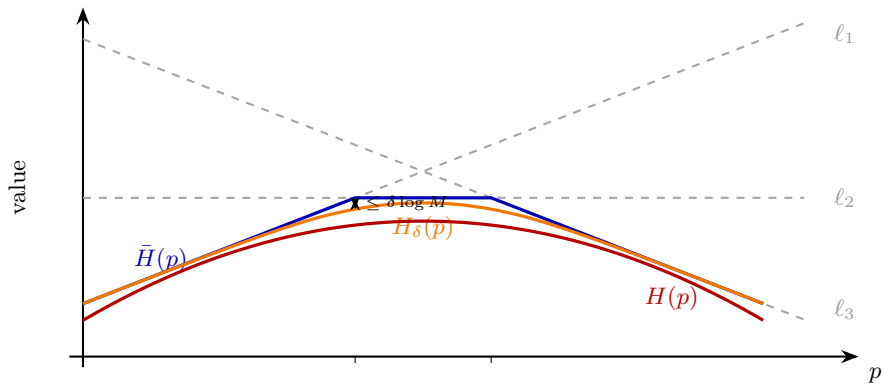%
}

\emph{Alternative smoothing.}\;
A Moreau--Yosida (infimal convolution) regularization could be used
instead, at the cost of solving an auxiliary minimization per
evaluation. We adopt the soft-min formula~\eqref{eq:H_delta} because
its derivatives are available in closed form via~\eqref{eq:grad_H_delta}.
\MSadd{\ In some settings the exact Hamiltonian admits a closed-form
expression whose single nonsmooth term can be regularized
analytically (e.g.\ $x^{10}-|p|\to x^{10}-\sqrt{p^2+\delta^2}$);
Section~\ref{sec:numerics} shows that in such cases the PA bundle
collapses to an exact finite set and the solver uses the analytic
smoothing directly.}%

\subsection{Notation summary}\label{ssec:notation}

The solver works with \MSdel{four}\Ksadd{3}{three} Hamiltonian objects:

\smallskip
\begin{center}
\begin{tabular}{lp{0.65\textwidth}}
\toprule
Symbol & Role \\
\midrule
$H$ & Exact (oracle) Hamiltonian~\eqref{eq:hamiltonian};
 concave in~$p$, possibly nonsmooth in~$p$ and~$x$. \\
$\bar H$ & PA bundle surrogate~\eqref{eq:pa_model}; concave,
 piecewise affine in~$p$, satisfies $\bar H\ge H$. \\
$H_\delta$ & Soft-min smoothing of~$\bar H$~\eqref{eq:H_delta};
 $C^\infty$ and concave in~$p$, and $C^{1,1}$
 (resp.~$C^2$ under~(A5)) in~$x$; used for all
 discrete dynamics. \\
\bottomrule
\end{tabular}
\end{center}
\smallskip

\noindent
All discrete dynamics in the solver are driven by~$H_\delta$; the
exact Hamiltonian~$H$ and the PA surrogate~$\bar H$ are used only for
error monitoring and bundle refinement. The coefficient functions
in~$H_\delta$ are $g_i(x,t)=f(x,a_i^*,t)$ and
$d_i(x,t)=\ell(x,a_i^*,t)$ (Section~\ref{ssec:extension}), where
$\{a_i^*\}_{i=1}^M$ are the controls returned by previous oracle
calls; the $x$-derivatives $\nabla_xg_i=\nabla_xf(x,a_i^*,t)$ and
$\nabla_xd_i=\nabla_x\ell(x,a_i^*,t)$ are computed from the problem
data and feed directly into~\eqref{eq:grad_H_delta}.
%% === end of Pformulation_v2.tex ===

\section{Discrete Solver and Adaptive Refinement}\label{sec:algorithm}

%% === inlined from Algorithm_v2.tex (R98z-cc arxiv-flat) ===
% ------------------------------------------------------------------
% §3 Discrete Solver and Adaptive Refinement
% ------------------------------------------------------------------

This section describes the complete numerical method. We first write
the discrete two-point boundary value problem (TPBVP) obtained by
applying the symplectic Euler scheme to the smoothed canonical
system~\eqref{eq:pmp}, then describe the damped Newton solver, define
the three a~posteriori error indicators, and present the adaptive
outer loop.

% ------------------------------------------------------------------
\subsection{Symplectic Euler discretization}\label{ssec:symp-euler}

Let $0=t_0<t_1<\cdots<t_N=T$ be a partition of $[0,T]$ with steps
$\Delta t_n=t_{n+1}-t_n$. The $p$-implicit symplectic Euler scheme
\SZadd{AM16}{(``$p$-implicit'' here means the costate $p_{n+1}$ is the
unknown solved-for at each step, while the state $x_{n+1}$ is
computed explicitly from $p_{n+1}$)}
for the smoothed canonical system driven by $H_\delta$
(Section~\ref{ssec:smoothing}) reads
\SLadd{P4}{}%
\begin{equation}\label{eq:symp_euler}
\begin{aligned}
x_{n+1} &= x_n + \Delta t_n\,\nabla_p H_\delta(p_{n+1},x_n,t_n), \\
p_n &= p_{n+1} + \Delta t_n\,\nabla_x H_\delta(p_{n+1},x_n,t_n),
\end{aligned}
\end{equation}
for $n=0,\dots,N-1$, with the two-point boundary conditions
\begin{equation}\label{eq:bc}
x_0 = x^0,\qquad p_N = \nabla g(x_N).
\end{equation}
The state equation in~\eqref{eq:symp_euler} is \emph{explicit} in
$x_{n+1}$ once $p_{n+1}$ is known, while the costate update is
\emph{implicit}---a structure that mirrors the forward-state /
backward-costate nature of the continuous system and preserves
its symplectic geometry~\cite{SanzSerna2016,Chyba2009}.

% ------------------------------------------------------------------
\subsection{Newton residual and Jacobian}\label{ssec:newton}

\paragraph{Full-space (all-at-once) nonlinear system.}
We adopt a full-space formulation in which both the discrete state
and costate values are treated as unknowns. Stacking
\begin{equation}\label{eq:unknowns}
w \;=\; (x_1,\dots,x_N,\;\MSdel{p_0,}p_1,\dots,p_N)^\top
 \;\in\;\mathbb{R}^{\MSdel{(2N+1)d}\MSadd{2Nd}},
\end{equation}
\MSnote{MS, p.~7: $p_0$ is not needed as an unknown; the discrete
costate equation at $n=0$ defines $p_0$ externally once the remaining
$(x,p)$ are known. Dropping it gives $2Nd$ unknowns / $2Nd$ equations.}%
we define the residual $R(w)$ by the state, costate, and boundary
blocks
\begin{equation}\label{eq:residual}
\begin{aligned}
F_n^{x}(w) &:= x_{n+1}-x_n-\Delta t_n\,
 \nabla_p H_\delta(p_{n+1},x_n,t_n)=0,
 && n=0,\dots,N-1, \\
F_n^{p}(w) &:= p_n-p_{n+1}-\Delta t_n\,
 \nabla_x H_\delta(p_{n+1},x_n,t_n)=0,
 && n=\MSdel{0,\dots,N-1}\MSadd{1,\dots,N-1}, \\
F^{T}(w) &:= p_N-\nabla g(x_N)=0,
\end{aligned}
\end{equation}
with the prescribed initial state $x_0=x^0$ held fixed (so $x_0$ is
not among the unknowns).
\MSadd{The boundary-layer costate $p_0$ is recovered after the solve
by the discrete adjoint equation at $n=0$:
$p_0 = p_1 + \Delta t_0\,\nabla_xH_\delta(p_1,x_0,t_0)$.}
System~\eqref{eq:residual} provides
\MSdel{$2Nd+d=(2N+1)d$} \MSadd{$Nd+(N-1)d+d = 2Nd$}
scalar equations, matching the dimension of~$w$.

\Rzzedit{A1a}{\emph{An algebraically equivalent variant.}\;
An algebraically equivalent full-space variant keeps $p_0$ as an unknown alongside
$(x_1,\dots,x_N,p_1,\dots,p_N)$, i.e.
$z=(x_1,\dots,x_N,p_0,\dots,p_N)^\top\in\mathbb{R}^{(2N+1)d}$, with the
matching boundary equation $F^p_0(z):=p_0 - p_1 - \Delta t_0
\nabla_xH_\delta(p_1,x_0,t_0)=0$ appended to the residual. This system has
$Nd+Nd+d=(2N+1)d$ scalar equations matching the
dimension of~$z$. The two formulations differ only in whether $p_0$ is
recovered from the discrete adjoint at $n=0$ post-solve (the convention~$w$ used here)
or solved as an unknown (the variant convention~$z$); the local time-coupling structure of the residual,
the Jacobian sparsity pattern, and the asymptotic complexity bound
$\mathcal{O}(Nd^3)$ stated below are identical between the two.}

\paragraph{Jacobian sparsity.}
With the natural ordering of~\eqref{eq:unknowns}, each residual block
couples only to variables at adjacent time indices $n$ and $n+1$.
Consequently, the Jacobian $J=\partial R/\partial w$ \Rzzedit{S12}{of
the residual $R$ defined in~\eqref{eq:residual}} is \emph{block-banded}: the only nonzero blocks are between $(x_n,p_{n+1})$
and $(x_{n+1},p_n)$ through the $(x,p)$-dependence of
$\nabla_pH_\delta$ and $\nabla_xH_\delta$, plus the terminal coupling
$F^T$ between $p_N$ and~$x_N$.
\MSadd{Equivalently, if the unknowns are permuted into the
time-interleaved order
$w=(x_1,p_1,x_2,p_2,\dots,x_N,p_N)^\top\in\mathbb{R}^{2Nd}$, then
the Jacobian is block-tridiagonal of size $2Nd\times 2Nd$ with
$2d\times 2d$ blocks and bandwidth $O(d)$, so Gaussian elimination
costs $O(Nd^3)$ by the standard banded-LU formula.
This interleaved view makes the complexity statement below
immediate.}%
\MSnote{MS, p.~8: suggested making the $O(Nd^3)$ count
self-evident by writing $w$ in interleaved order.}
The linear system $J\,\Delta w=-R$ can
therefore be solved in $\mathcal{O}(Nd^3)$ operations by a single
forward/backward block sweep, i.e.\ \emph{linear} in~$N$ for fixed
state dimension~$d$. Second derivatives of $H_\delta$ that appear in
$J$ exist under the strengthened regularity~(A5)
(Section~\ref{ssec:assumptions}); if only (A1) is assumed, $J$ should
be interpreted as a generalized or quasi-Newton Jacobian.

\MSadd{%
Figure~\ref{fig:stencil-sparsity} summarises both the local
symplectic-Euler coupling on the time grid and the resulting
block-tridiagonal sparsity of the Newton Jacobian.%
}

\MSadd{%
\begin{figure}[htbp]
\centering
\begin{tikzpicture}[>=Stealth, font=\footnotesize]
 \begin{scope}
 \node[anchor=south, font=\small\bfseries] at (2.1, 3.1)
 {(a) Symplectic-Euler stencil};
 \draw[->, thick] (-0.3, 0) -- (4.7, 0) node[below] {$t$};
 \foreach \i/\lbl in {0/{n{-}1}, 1/{n}, 2/{n{+}1}, 3/{n{+}2}} {%
 \pgfmathsetmacro{\xpos}{0.5 + \i*1.2}
 \draw[fill=black] (\xpos, 0) circle (1.5pt);
 \node[font=\tiny, below=1pt] at (\xpos, 0) {$t_{\lbl}$};
 \node[blue!70!black, font=\small] at (\xpos, 2.15) {$x_{\lbl}$};
 \node[red!70!black, font=\small] at (\xpos, 0.65) {$p_{\lbl}$};
 }
 \foreach \i in {0,1,2} {%
 \pgfmathsetmacro{\xs}{0.5 + \i*1.2}
 \pgfmathsetmacro{\xt}{0.5 + (\i+1)*1.2}
 \draw[->, blue!70!black, thick]
 (\xs+0.14, 2.15) to[bend left=22] (\xt-0.14, 2.15);
 }
 \foreach \i in {3,2,1} {%
 \pgfmathsetmacro{\xs}{0.5 + \i*1.2}
 \pgfmathsetmacro{\xt}{0.5 + (\i-1)*1.2}
 \draw[->, red!70!black, thick]
 (\xs-0.14, 0.65) to[bend left=22] (\xt+0.14, 0.65);
 }
 \node[blue!70!black, anchor=west, font=\scriptsize]
 at (-0.3, 2.80)
 {$x_{n+1}=x_n+\Delta t_n\,\nabla_p H_\delta(p_{n+1},x_n,t_n)$};
 \node[red!70!black, anchor=west, font=\scriptsize]
 at (-0.3, -0.75)
 {$p_n=p_{n+1}+\Delta t_n\,\nabla_x H_\delta(p_{n+1},x_n,t_n)$};
 \end{scope}
 \begin{scope}[xshift=7.8cm]
 \node[anchor=south, font=\small\bfseries] at (1.5, 3.1)
 {(b) Block-tridiagonal $J$};
 \def\bs{0.52}
 \foreach \i in {0,...,4} {%
 \pgfmathsetmacro{\yb}{(4-\i)*\bs}
 \fill[blue!25] (\i*\bs, \yb) rectangle (\i*\bs+\bs, \yb+\bs);
 }
 \foreach \i in {0,...,3} {%
 \pgfmathsetmacro{\ya}{(4-\i)*\bs}
 \pgfmathsetmacro{\yc}{(3-\i)*\bs}
 \fill[orange!30] (\i*\bs+\bs, \ya)
 rectangle (\i*\bs+2*\bs, \ya+\bs);
 \fill[orange!30] (\i*\bs, \yc)
 rectangle (\i*\bs+\bs, \yc+\bs);
 }
 \foreach \k in {0,...,5} {%
 \draw[gray!60, thin] (\k*\bs, 0) -- (\k*\bs, 5*\bs);
 \draw[gray!60, thin] (0, \k*\bs) -- (5*\bs, \k*\bs);
 }
 \draw[black, thick] (0, 0) rectangle (5*\bs, 5*\bs);
 \node[anchor=north, font=\tiny] at (2.5*\bs, -0.05)
 {$w=(x_1,p_1,\dots,x_N,p_N)^\top$};
 \node[anchor=west, font=\tiny] at (5*\bs + 0.25, 4*\bs + 0.25)
 {\tikz\draw[fill=blue!25] (0,0) rectangle (0.28,0.28);\
 diagonal};
 \node[anchor=west, font=\tiny] at (5*\bs + 0.25, 3*\bs + 0.25)
 {\tikz\draw[fill=orange!30] (0,0) rectangle (0.28,0.28);\
 off-diag};
 \node[anchor=west, font=\tiny] at (5*\bs + 0.25, 2*\bs + 0.25)
 {$2d\times 2d$ blocks};
 \node[anchor=west, font=\tiny] at (5*\bs + 0.25, 1*\bs + 0.25)
 {LU cost $O(Nd^3)$};
 \end{scope}
\end{tikzpicture}
\caption{Coupling and sparsity structure of the discrete TPBVP.
(a)~Symplectic-Euler stencil on the time grid: the forward state
update (blue) and the backward costate update (red) each require one
evaluation of $\nabla_{p,x}H_\delta$ and couple only adjacent
time indices. (b)~Resulting Newton-Jacobian sparsity in the
interleaved ordering $w=(x_1,p_1,\dots,x_N,p_N)^\top$:
block-tridiagonal with $2d\times 2d$ blocks, giving an
$\mathcal{O}(Nd^3)$ per-iteration cost.}
\label{fig:stencil-sparsity}
\end{figure}%
}

\paragraph{Implementation options.}
The same mathematical system~\eqref{eq:residual} admits several
linear-algebra realizations. A \emph{condensed single-shooting}
variant eliminates the states via the explicit forward sweep, leaving
a dense $Nd\times Nd$ system in the costates only (the sensitivity
chain $\partial x_n/\partial p_j$ for $j\le n$ is the source of the
density). A \emph{multiple-shooting} variant in the sense of Bock and
Plitt~\cite{BockPlitt1984} introduces $S$ segment boundary variables
$x^{(s)}$ with explicit continuity defects; it yields a
block-tridiagonal Jacobian of block size~$d$ and is solved in
$\mathcal{O}(Sd^3)$ operations. Both realizations are
algebraically equivalent to~\eqref{eq:residual}. In the experiments
of Section~\ref{sec:numerics} we use the full-space formulation
directly.

\paragraph{Damped Newton iteration.}
Given an initial guess $w^{(0)}$ (e.g., from a coarser mesh, or an
explicit forward sweep from a zero costate), iterate for
$k=0,1,2,\dots$:
\begin{enumerate}
\item Residual and Jacobian: evaluate $R(w^{(k)})$ and the
 \Rzzedit{B2-S14}{block-banded Jacobian $J(w^{(k)})=\partial R/\partial w$} from~\eqref{eq:residual}.
\item Newton step: solve $J^{(k)}\Delta w^{(k)}=-R^{(k)}$ by a
 block-banded elimination.
\item Damped update: set
 $w^{(k+1)}=w^{(k)}+\alpha_k\,\Delta w^{(k)}$, where the step
 length $\alpha_k\in(0,1]$ is chosen by Armijo backtracking on the
 \SLadd{N5a}{infinity-norm residual:\ starting from $\alpha=1$, halve
 until $\|R(w+\alpha\Delta w)\|_\infty\le
 (1-c_1\alpha)\|R(w)\|_\infty$ with $c_1=10^{-4}$ \Rzzedit{A1b}{(the $\ell^\infty$
 acceptance test, rather than the squared $\ell^2$ merit, is a natural choice for
 the full-space discretization because it controls the largest
 residual component directly)}}.
\item \SLadd{N5b}{Stop when $\|R(w^{(k)})\|_\infty<\varepsilon_{\mathrm{Newton}}$,
 or after $k_{\max}$ iterations}.\Rzzedit{A1c}{}
\end{enumerate}
\RTadd{R70-M2}{Newton failure (residual not below
$\varepsilon_{\mathrm{Newton}}$ after $k_{\max}$ steps) is handled
by the nonlinear least-squares fallback documented below; the
next outer iteration's $\delta$-action is then driven by
$\eta_\delta$ as usual (no explicit Newton-robustness
$\delta$-raise branch is implemented in the present solver).}

\medskip
The full Newton solver is stated as Algorithm~\ref{alg:newton}.

\begin{algorithm}[H]
\caption{Damped Newton solver for the smoothed TPBVP (full-space system)}\label{alg:newton}
\begin{algorithmic}
\Require Mesh $\{t_n\}_{n=0}^N$, smoothed Hamiltonian $H_\delta$,
 \MSdel{initial guess $w^{(0)}=(x_1^{(0)},\dots,x_N^{(0)},p_0^{(0)},\dots,p_N^{(0)})$,
 fixed $x_0=x^0$,}%
 \MSadd{prescribed initial state $x_0=x^0$ (the remaining
 components of $w^{(0)}$ are generated internally by the horizon
 continuation of Section~\ref{ssec:init}),}
 tolerance
 \RTadd{R64}{$\varepsilon_{\mathrm{Newton}}$,}
 maximum iterations $k_{\max}$%
\MSnote{MS, p.~8: ``fixed $x_0$ sufficient'' --- the rest of $w^{(0)}$
is supplied by horizon continuation, not by the caller.}
\Ensure Approximate solution $(x_0,\dots,x_N,p_0,\dots,p_N)$
\For{$k=0,1,\dots,k_{\max}-1$}
 \State Evaluate residual $R(w^{(k)})=(F^x,F^p,F^T)$
 via~\eqref{eq:residual}
 \State Assemble block-banded Jacobian $J^{(k)}=\partial R/\partial w$
 \State Solve \MSdel{$J^{(k)}\Delta w=-R^{(k)}$}%
 \MSadd{$J^{(k)}\,\Delta w^{(k)}=-R(w^{(k)})$}
 by block-banded elimination%
 \MSnote{MS, p.~8: use the same $w^{(k)}$ / $\Delta w^{(k)}$
 notation as introduced in line~5 above.}
 \State \SLadd{N5a}{Line search: find $\alpha_k$ by Armijo backtracking on $\|R\|_\infty$}
 \State $w^{(k+1)}\leftarrow w^{(k)}+\alpha_k\,\MSdel{\Delta w}\MSadd{\Delta w^{(k)}}$
 \If{\SLadd{N5b}{$\|R\|_\infty<\varepsilon_{\mathrm{Newton}}$}}
 \State \Return converged solution
 \EndIf
\EndFor
\State \Return failure to the outer solve wrapper \RTadd{R66-M1}{(the
outer loop invokes the least-squares fallback of the
following paragraph and, on the next outer iteration, refines whichever
layer's indicator-to-tolerance ratio exceeds tolerance)}
\end{algorithmic}
\end{algorithm}

\Ksadd{4,6}{A third common cause of Newton stall, in addition to large
$\delta$ (poor smoothing-driven conditioning) and coarse mesh (residual
discretization bias), is information-starvation in the PA bundle: when
$M$ is small or the sampled oracle responses are clustered, the
surrogate $\bar H$ misses the local curvature of $H$ and the Newton
iterates leave its basin. In this regime the outer loop's
bundle-enrichment branch (Algorithm~\ref{alg:outer} below) and the
horizon-continuation warm-start of~\S\ref{ssec:init} are the two
designed responses; a poorly initialized iterate is most efficiently
cured by the latter.}

\MSadd{%
\paragraph{Sparse linear algebra.}
\MSdel{The Jacobian of the full-space system in
$w=(x_1,\dots,x_N,p_1,\dots,p_N)$ is block-tridiagonal under the
interleaved ordering of~\S\ref{ssec:newton} and is assembled in
compressed-sparse form at every Newton step.}%
\Ksadd{5}{The Jacobian of the full-space system, after the interleaved
re-indexing $w=(x_1,p_1,\dots,x_N,p_N)^\top$ of~\S\ref{ssec:newton}
(alternating one state node with one costate node so that consecutive
nodes share a stencil), is block-tridiagonal and is assembled in
compressed-sparse form at every Newton step.}\RTadd{R62-B2code}{\ The
interleaved ordering is a mathematical re-indexing used for stating the
block-tridiagonal structure;\ in the implementation, the natural
grouped layout $(x_1,\dots,x_N,p_0,\dots,p_N)^\top$ is equivalent under
a permutation and is used directly for the sparse-LU factorization
without materializing the interleaved permutation.} The resulting linear
system is solved with a sparse LU factorization
(\texttt{scipy.sparse.linalg.splu}), which preserves the
$\mathcal{O}(Nd^3)$ cost without incurring the constant factor of
dense linear algebra. Sparse assembly is essential when adaptive
refinement drives~$N$ into the hundreds or low thousands, where a
dense solve would dominate the wall-clock budget.

\paragraph{Nonlinear least-squares fallback.}\Rzzedit{B2-S18}{\label{para:LS-fallback}}
When the damped-Newton iteration stalls---typically because the
current iterate is outside Newton's basin of attraction on a
nonsmooth or singular example---the solver falls back to a nonlinear
least-squares \Rzzedit{B2-S18}{(LS in what follows)} minimization of the residual norm
$\tfrac12\|R(w)\|_2^2$ using a
\SLadd{N6a}{trust-region reflective nonlinear least-squares fallback
 (\texttt{scipy.optimize.least\_squares} with
 \texttt{method="trf"})}. The fallback is
\emph{not} a replacement for Newton:\ it is used as a robustness
mechanism on a per-outer-iteration basis;\
\SLadd{N6b}{the least-squares output is returned as the solution of the
current outer iteration and is then used as a warm start for the
\emph{next} outer iteration's Newton solve (the Newton loop does not
resume within the same outer iteration after the fallback finishes)}.
\SLadd{47-1}{In practice the fallback acts as a safety valve at the
coarsest meshes of the most nonsmooth examples and is rarely active
once the mesh has been refined and the smoothing~$\delta$ has been
reduced;\ the runs of \S\ref{sec:numerics} arm the fallback with
$\varepsilon_{\mathrm{Newton}}=10^{-10}$ and $k_{\max}=50$
(common-setup paragraph of~\S\ref{sec:numerics}).}%
}%
\MSnote{2026-04-24 update, \S3.10--3.11: document the sparse LU
assembly and the nonlinear least-squares fallback; both are now part
of the implementation but were absent from the draft.}

\RTadd{R62-B7}{\paragraph{Newton-failure handling and accuracy-driven $\delta$-halving.}
\Rzzedit{S16}{The smoothing parameter~$\delta$ is acted on by a single
mechanism in the present solver, which we call \emph{accuracy-driven
$\delta$-halving} and define here. The controller logic that selects
this action when needed is the overload-ratio outer loop of
\S\ref{ssec:outer-loop}.}\ Once the
inner damped-Newton loop of Algorithm~\ref{alg:newton} (or the
least-squares fallback of the preceding paragraph, when Newton stalls)
returns a solution, the outer-loop controller of~\S\ref{ssec:outer-loop}
consults the smoothing-bias indicator $\eta_\delta$ defined
in~\S\ref{ssec:indicators}.\ If $\eta_\delta>\varepsilon_\delta$, the
smoothing bias $\bar H-H_\delta$ is the dominant contributor to the
value-error budget, and the controller \emph{halves} $\delta$ (or
applies the overload-ratio scaling of Algorithm~\ref{alg:outer}).\
Newton stalls themselves are absorbed by the LS fallback on the
current outer iteration;\ the LS output is then used as a warm start
for the next outer iteration's Newton solve, with the next iteration's
$\delta$-action driven by $\eta_\delta$ as usual.\ The
robustness/accuracy distinction is therefore handled \emph{within} the
inner solve (Newton$\to$LS fallback) rather than by a separate outer
$\delta$-raise branch;\ extending the controller with an explicit
Newton-robustness $\delta$-raise (raising $\delta$ on the next outer
iteration whenever Newton/LS stalls) is a natural enhancement and is
left to future work.}

% ------------------------------------------------------------------
\subsection{Initialization by horizon continuation}\label{ssec:init}

The Newton iteration on the full-space system
$F(w)=0$ converges quadratically inside a neighbourhood of the solution
whose radius shrinks with the horizon $T$, the nonlinearity of $f$,
and the stiffness of the adjoint equation. A naive cold start of the
form $p_n^{(0)}\equiv\nabla g(x_0)$ is adequate for short horizons
but frequently fails for the moderate-to-long horizons encountered in
practice. We therefore recommend the following cheap, deterministic
initialization scheme, which we have found to be very robust in our
experiments and adds negligible cost relative to a single Newton
solve at the target horizon.

\paragraph{Horizon continuation.}
Fix an integer $K\ge 0$ and define a geometric schedule of horizons
\begin{equation}\label{eq:horizon_schedule}
T_k \;:=\; 2^{\,k-K}\,T,\qquad k=0,1,\dots,K,
\end{equation}
so that $T_0=2^{-K}T$ is small enough that the sub-problem on
$[0,T_0]$ is nearly trivial and $T_K=T$ is the target horizon.
We solve the sequence of sub-problems~$(P_k)$ obtained by replacing
$T$ with $T_k$ in~\eqref{eq:bolza}, using the
\SLadd{47-2}{converged state--costate trajectory $(x^{(k-1)},p^{(k-1)})$}
of~$(P_{k-1})$ as the initial iterate for the Newton solve of~$(P_k)$.

\begin{itemize}
\item \emph{Base case $k=0$.} For $T_0\ll 1$ the running cost
 contributes $O(T_0)$ and the state barely moves, so the terminal
 adjoint condition $p(T_0)=\nabla g(x(T_0))$ is well approximated by
 the constant guess $p_n^{(0)}\equiv\nabla g(x_0)$ on the
 coarse mesh. A single Newton iterate usually yields a solution
 accurate to machine precision.
\item \emph{Continuation step $k-1\to k$.} Let $(x^{(k-1)},p^{(k-1)})$
 be the converged discrete trajectory on $[0,T_{k-1}]$, sampled at
 nodes $0=t_0^{(k-1)}<\cdots<t_{N_{k-1}}^{(k-1)}=T_{k-1}$. We build
 an initial guess on $[0,T_k]$ by \emph{time stretching}:
 \begin{equation}\label{eq:time_stretch}
 \tilde p_n^{(k)} \;:=\; p^{(k-1)}\!\left(
 \tfrac{T_{k-1}}{T_k}\,t_n^{(k)}\right),\qquad
 \tilde x_n^{(k)} \;:=\; x^{(k-1)}\!\left(
 \tfrac{T_{k-1}}{T_k}\,t_n^{(k)}\right),
 \end{equation}
 evaluated by piecewise-linear interpolation between the nodes of
 the previous mesh.
\end{itemize}

\Rzzedit{S17}{\paragraph{Optional forward sweep.}
To sharpen the time-stretched warm start of \eqref{eq:time_stretch}
before Newton is invoked, we may run one forward sweep of the
smoothed canonical state equation along the time-stretched
costate~$\tilde p$. By the coefficient extension of
\S\ref{ssec:extension}, this equation reads
\begin{equation}\label{eq:warm-start-vfield}
 \dot x \;=\; \nabla_p H_\delta(p, x, t) \;=\;
 \sum_{i=1}^{M} w_i(x, p, t)\, f(x, a_i^\star, t),
\end{equation}
where the softmax weights $w_i$ are those of
\eqref{eq:softmax_weights} from \S\ref{ssec:smoothing} and the
$a_i^\star$ are the oracle-returned controls at the PA bundle
samples. Equation~\eqref{eq:warm-start-vfield} is the equation that
the Newton residual $F(w) = 0$ of \S\ref{ssec:newton} enforces, so
integrating it along $\tilde p$ produces a state guess that is
already consistent with the Newton system before Newton is invoked.
Each $a_i^\star$ is individually admissible by \S\ref{ssec:oracle},
so the weighted sum on the right-hand side of
\eqref{eq:warm-start-vfield} is well-defined for arbitrary admissible
sets $\mathcal{A}$, convex or not, without any post-processing step.}

\MSnote{R98z Phase B v2 (F1): retracts the earlier "softmin control
selector $a_\delta^\star$" definition and the convex-admissibility
caveat. The selector $a_\delta^\star$ would have been needed only if
the warm-start sweep integrated $\dot x = f(x, a_\delta^\star, t)$
(the control-averaged form); since we chose C1 (the drift-averaged
form, eq:warm-start-vfield), the selector and the caveat are not
load-bearing. Equations eq:softmin\_selector and eq:hard\_mode\_idx
are also retracted; eq:warm-start-vfield is preserved. The previous
"Optional, ... $u^\star_\delta$ ... sharpens the guess" inline
definition + the ``Convex-admissibility caveat'' MSadd block +
Remark 3 (\texttt{rmk:warm-start-state-eq}) are all absorbed into
the single paragraph above.}
\SLnote{R98z Phase B v2: this absorbs the content of the SL-added
Remark 3 from R47-3.}%

\paragraph{Compatibility with the outer loops.}
Horizon continuation is orthogonal to, and composes transparently
with, the other two continuation mechanisms already present in our
solver:
(i) the smoothing homotopy $\delta_0>\delta_1>\cdots>\delta$ of the
outer loop, which uses a strongly smooth Hamiltonian at large
$\delta$ as a trivial warm start for successively sharper ones; and
(ii) the bundle-enrichment homotopy $M=1\to M=2\to\cdots$, which
grows the PA surrogate one plane at a time. A robust default is to
run horizon continuation as the outermost loop, with $\delta$- and
$M$-continuation nested inside at each horizon $T_k$. In all
numerical examples of Section~\ref{sec:numerics}, starting from a constant
guess $p\equiv\nabla g(x_0)$ on~$T_0=T/8$ and doubling the horizon
three times sufficed for full Newton convergence without a single
line-search failure.

\paragraph{Alternatives and fall-backs.}
\SLadd{47-4}{Horizon continuation with $K=3$ doublings and a constant
initial costate $p\equiv\nabla g(x_0)$ on $T_0=T/8$ sufficed in every
benchmark of~\S\ref{sec:numerics} (see~\S\ref{ssec:numerics-summary}).
If horizon continuation stalls (e.g.\ when the qualitative structure
of the optimal control changes across horizons), recommended remedies
are a Riccati warm start (linearise $f$ and $L$ along an uncontrolled
forward simulation and solve the resulting LQR by a single backward
Riccati sweep) and classical numerical continuation in a scalar
homotopy parameter deforming the dynamics from a linear surrogate to
the target system~\cite{BockPlitt1984}.}

% ------------------------------------------------------------------
\subsection{A posteriori error indicators}\label{ssec:indicators}

After the Newton solve converges, we evaluate three indicators that
quantify the distinct error sources. Each indicator is defined at the
\emph{discrete} level, using only quantities already available from
the solve.

\MSadd{\paragraph{Value-first alignment.}
Consistent with the value-first scope of \S\ref{sec:introduction}, each of the three
indicators $\eta_{\mathrm{time}}, \eta_{\mathrm{PA}}, \eta_\delta$ is a \emph{value-error}
\RTadd{R66-M2}{diagnostic} for $|J_h - U(0,x_0)|$ in the construction-order telescope
\[
 |J_h - U(0,x_0)| \;\lesssim\; E_{\mathrm{PA}} + E_\delta + E_h + E_{\mathrm{opt}},
\]
where $E_{\mathrm{PA}}$ collects the PA-surrogate bias controlled by $\eta_{\mathrm{PA}}$,
$E_\delta$ the soft-min smoothing bias controlled by $\eta_\delta$, $E_h$ the time-discretization
bias controlled by $\eta_{\mathrm{time}}$, and $E_{\mathrm{opt}}$ the inner Newton residual at
termination. No indicator is a state-error, costate-error, or control-error bound; the
controller refines only the layers whose value-error \RTadd{R66-M2}{diagnostic} exceeds the corresponding
tolerance.}

\paragraph{(I) Time-step error indicator $\eta_{\mathrm{time}}$.}
\SZadd{AR1}{\emph{Role:}\ $\eta_{\mathrm{time}}$ is the max-type
indicator that drives mesh refinement, in contrast to the sum-type
$\eta_{\mathrm{PA}}$ and $\eta_\delta$ that bound contributions to the
value (see item~(IV) below).\ }%
Following \Rzzedit{S19}{the error-density framework
of~\cite{KarlssonLarssonSandbergTempone2015}}, the
time-discretization error of the symplectic Euler
scheme~\eqref{eq:symp_euler} admits the a~posteriori representation
\begin{equation}\label{eq:time_error_rep}
\RTadd{R69-C7}{\bar U_h(x_0,0) - U(x_0,0)} \;=\;
\sum_{n=0}^{N-1}\Delta t_n^2\,\rho_n \;+\; R,
\qquad |R|\le C'\Delta t_{\max}^2,
\end{equation}
where the error density is
$\rho_n=-\nabla_pH_\delta(p_{n+1},x_n,t_n)\cdot
\nabla_xH_\delta(p_{n+1},x_n,t_n)/2$
and \RTadd{R69-C7}{$\bar U_h$, $U$} are the discrete and continuous
optimal values, respectively (aligned with the value-first scope
notation $U(0,x_0)$ of~\S\ref{sec:introduction}). \RTadd{R62-H1}{The implementation evaluates $\rho_n$
directly from this leading-order Hamiltonian-gradient formula at every
node;\ the gradients $\nabla_pH_\delta$ and $\nabla_xH_\delta$ are
already required by the Newton residual and Jacobian, so $\rho_n$
adds no new derivative cost.}\ An algebraically equivalent computable
variant $\tilde\rho_n$ of~\cite[Remark~2.5]{KarlssonLarssonSandbergTempone2015}
expresses the same leading-order quantity using only $H_\delta$
evaluations:
\begin{equation}\label{eq:computable_density}
\tilde\rho_n \;:=\;
\frac{H_\delta(p_{n+1},x_n,t_n)}{\Delta t_n}
-\frac{H_\delta(p_n,x_n,t_n)+H_\delta(p_{n+1},x_{n+1},t_{n+1})}{2\Delta t_n}
+\frac{p_n-p_{n+1}}{2}\cdot
\frac{\nabla_pH_\delta(p_{n+1},x_n,t_n)}{\Delta t_n},
\end{equation}
which requires only Hamiltonian evaluations and one partial derivative
at already computed grid values.
\Ksadd{11}{The two densities are equivalent at leading order and both
are approximate: $\rho_n$ in~\eqref{eq:time_error_rep} is the
\emph{leading-order asymptotic density} of the symplectic-Euler
truncation error, while $\tilde\rho_n$
in~\eqref{eq:computable_density} is the algebraically equivalent
\emph{symplectic-Euler-consistent re-expression}
of~\cite[Remark~2.5]{KarlssonLarssonSandbergTempone2015}, obtained by
substituting the discrete velocities $(x_{n+1}-x_n)/\Delta t_n$ and
$(p_n-p_{n+1})/\Delta t_n$ for the gradients
$\nabla_p H_\delta$ and~$\nabla_x H_\delta$ that
appear in~$\rho_n$. The two formulas require the same underlying
regularity of~$H_\delta$ (Taylor-expanding $\tilde\rho_n$ recovers
$\rho_n$ to leading order in $\Delta t_n$ under the symplectic-Euler
relation); $\tilde\rho_n$ is not a more robust object in that sense.
\RTadd{R62-H1}{The implementation evaluates $\rho_n$ rather than
$\tilde\rho_n$: both forms compute the same leading-order density,
but $\nabla_p H_\delta$ and $\nabla_x H_\delta$ are \Rzzedit{S20}{already
computed during the Newton-residual evaluation and reused here}, so
$\rho_n$ is the cheapest
available expression in the present code. The
$\tilde\rho_n$~form replaces $\nabla_x H_\delta$ by two extra
scalar $H_\delta$ evaluations and is the natural choice if the
$\nabla_x H_\delta$ derivative is unavailable;\ the two forms are
algebraically equivalent under the symplectic-Euler identities, so
substituting one for the other does not change the asymptotic
indicator. Neither density is connected to the time-stretching
warm-start of~\S\ref{ssec:init}: the warm-start provides an initial
iterate, while $\rho_n$ is evaluated only \emph{after} Newton
convergence at each outer iterate.}}

The stabilized error indicators are
$\bar r_n:=\bar\rho_n\,\Delta t_n^2$
with
\RTadd{R62-H2}{$\bar\rho_n:=\max(|\rho_n|,\RTadd{R69-C4}{\kappa}\sqrt{\Delta t_{\max}})$}
and floor constant \RTadd{R69-C4}{$\kappa=10^{-6}$ (using $\kappa$ to
avoid clash with the horizon-continuation integer $K$
of~\S\ref{ssec:init})}.\RTadd{R62-H2}{\ Since
$\eta_{\mathrm{time}}$ and the equidistribution coefficients
of~\S\ref{ssec:quasi-norm-advantage} use only $|\bar\rho_n|$, dropping
the sign of $\rho_n$ in $\bar\rho_n$ has no effect on either the
stopping quantity or the step-count coefficients;\ the signed
$\rho_n$ is retained for diagnostic plots when a signed density is
required for visual inspection.}
The aggregate time-step indicator and its companion sum are
\begin{equation}\label{eq:eta_time}
\eta_{\mathrm{time}}^{\max} \;:=\; \max_{0\le n\le N-1}\bar r_n,
\qquad
\eta_{\mathrm{time}}^{\mathrm{sum}} \;:=\; \sum_{n=0}^{N-1}\bar r_n.
\end{equation}
We use $\eta_{\mathrm{time}}^{\max}$ as the \emph{stopping
quantity}: the time discretization is considered acceptable as soon
as the worst interval falls below the per-interval tolerance
$\varepsilon_{\mathrm{time}}/N$, in line with the equidistribution
view of the local error density $\rho_n$ in~\cite{KarlssonLarssonSandbergTempone2015}.
\SZadd{AR2}{Using the \emph{max}-type quantity (rather than the
sum-type companion) as the universal stop is deliberate:\ it
guarantees that \emph{every} discretization interval satisfies the
per-cell threshold, whereas a sum-type stop would admit a
configuration in which a small number of high-density cells are
masked by many low-density cells.\ }%
The companion quantity $\eta_{\mathrm{time}}^{\mathrm{sum}}$ is
reported alongside the max for diagnostic purposes but does not
enter the stopping rule. Where the symbol
$\eta_{\mathrm{time}}$ appears without a superscript in the
remainder of this work, it always refers to the
max-type stopping quantity~$\eta_{\mathrm{time}}^{\max}$.

\paragraph{(II) PA modeling error indicator $\eta_{\mathrm{PA}}$.}
We sample the gap between the PA surrogate and the oracle
Hamiltonian along the current trajectory:
\begin{equation}\label{eq:eta_PA}
\eta_{\mathrm{PA}} \;:=\;
\sum_{n=0}^{N-1}\Delta t_n\,
\bigl[\bar H(p_{n+1},x_n,t_n)-H(p_{n+1},x_n,t_n)\bigr].
\end{equation}
Each evaluation $H(p_{n+1},x_n,t_n)$ calls the oracle and
simultaneously produces a new subgradient that can be used for
bundle enrichment.
\RTadd{R51}{\ Throughout this work the PA coefficient extension is
the control-based one of~\S\ref{ssec:extension}, which gives a global
majorant $\bar H\ge H$; the gap~$\eta_{\mathrm{PA}}\ge 0$ then admits
the upper-bound interpretation used below.}

\paragraph{(III) Smoothing bias indicator $\eta_\delta$.}
The smoothing bias is bounded analytically
by~\eqref{eq:smoothing_bias}; we also monitor it empirically:
\begin{equation}\label{eq:eta_delta}
\eta_\delta \;:=\;
\sum_{n=0}^{N-1}\Delta t_n\,
\bigl[\bar H(p_{n+1},x_n,t_n)-H_\delta(p_{n+1},x_n,t_n)\bigr]
\;\le\; T\,\delta\log M.
\end{equation}

\MSadd{%
\paragraph{(IV) Asymmetric roles of the three indicators.}
The three indicators~\eqref{eq:eta_time},~\eqref{eq:eta_PA},
and~\eqref{eq:eta_delta} are not used interchangeably, and their
mathematical roles are deliberately different:

\begin{itemize}
\item $\eta_{\mathrm{time}}=\max_n\bar r_n$ is a \emph{max-type}
 indicator:\ the time discretization is considered acceptable as soon
 as the worst interval falls below tolerance, because the local
 error density $\rho_n$ is the object that drives mesh refinement.
\item $\eta_{\mathrm{PA}}$ and $\eta_\delta$ are \emph{sum-type}
 (interval-weighted) indicators:\ they accumulate local modeling and
 smoothing gaps and therefore behave as empirical upper bounds on the
 total contribution of those two error sources to the objective.
\RTadd{R70-M3}{\ The one-sided interpretation of
 $\eta_{\mathrm{PA}}$ in this paragraph relies on the control-based
 majorant extension of~\S\ref{ssec:extension};\ if that extension is
 replaced by a local black-box model, $\eta_{\mathrm{PA}}$ becomes
 diagnostic only.}
\end{itemize}

}%
\MSnote{2026-04-24 update: document the asymmetric (max vs sum) roles
of the three indicators and the reported-but-non-stopping
$\eta_{\mathrm{time,sum}}$.}

% ------------------------------------------------------------------
\subsection{Adaptive mesh refinement}\label{ssec:adaptive-mesh}

The time mesh is refined following
Algorithm~2.6 of~\cite{KarlssonLarssonSandbergTempone2015}.
Given a tolerance $\varepsilon_{\mathrm{time}}>0$, a subdivision
factor $M_{\mathrm{sub}}\ge 2$, and a safety parameter $s\in(0,1)$:

\begin{enumerate}
\item If $\eta_{\mathrm{time}}^{\max}<\varepsilon_{\mathrm{time}}/N$,
 the mesh is accepted.
\item Otherwise, every interval with
 $\bar r_n>s\,\varepsilon_{\mathrm{time}}/N$ is subdivided into
 $M_{\mathrm{sub}}$ equal parts, the mesh is updated, and the
 Newton solve is restarted (warm-started from interpolated values).
\end{enumerate}
The same per-interval threshold $\varepsilon_{\mathrm{time}}/N$ is
used inside the outer loop (Algorithm~\ref{alg:outer}) to decide
whether the time discretization is the binding indicator.

\paragraph{Scope of the imported theory.}
For a \emph{fixed} smoothed Hamiltonian~$H_\delta$ and a \emph{fixed}
bundle, the refinement strategy above is exactly the one
\SZadd{AR4}{analyzed in}~\cite{KarlssonLarssonSandbergTempone2015}, and the
following statements apply at that fixed inner stage. Under the
regularity conditions
of~\cite[Theorem~2.7]{KarlssonLarssonSandbergTempone2015}, each
refinement pass strictly decreases the maximal indicator provided
$M_{\mathrm{sub}}^2>c^{-1}$ and $s\le c/M_{\mathrm{sub}}$
\MSdel{(where $c$ is the error-density ratio bound)}\Ksadd{13}{(here
$c\ge 1$ is the error-density ratio bound from~\cite[Theorem~2.7 and
the surrounding discussion]{KarlssonLarssonSandbergTempone2015}: on
every refined subinterval, the ratio of the local maximum of
$|\bar\rho|$ to the local mean of $|\bar\rho|$ is at most $c$. In
plain terms, $c$ measures how non-uniformly the error density varies
on the subinterval; the conditions $M_{\mathrm{sub}}^2>c^{-1}$ and
$s\le c/M_{\mathrm{sub}}$ require the bisection factor to be large
enough and the safety margin to be small enough that this
non-uniformity cannot defeat the contraction.)}%
\RTadd{AR3-mstz}{\ \emph{Sign-convention note.}\;Two conventions for
the comparability constant coexist in the equidistribution
literature:\ the max-over-mean form used in the parenthetical above
gives a $c\ge 1$, while the lower-comparability convention of
\cite[eq.~(2.22)]{KarlssonLarssonSandbergTempone2015} uses $0<c\le 1$,
with the large distortion number being $c^{-1}$. Both refer to the
same structural quantity. In the rest of this section we standardize
on the latter convention:\ $c\in(0,1]$, and the conditions
$M_{\mathrm{sub}}^{2}>c^{-1}$ and $s_1\le c/M_{\mathrm{sub}}$ are
read with that sign, where $s_1\in(0,1)$ is the marking parameter of
Algorithm~2.6, distinct from~$c$}.
The accuracy and
efficiency guarantees of~\cite[Theorems~2.8--2.9]{KarlssonLarssonSandbergTempone2015}
then give
$\limsup_{\mathrm{TOL}\to 0}\mathrm{TOL}^{-1}|\text{error}|\le 1$
and a quasi-optimal step count
$\mathrm{TOL}\cdot N\le C\|\bar\rho/c\|_{L^{1/2}}$.
\SZadd{AR3}{Here the $L^{1/2}$ functional is the quasi-norm
$\|\bar\rho/c\|_{L^{1/2}} :=\bigl(\int_0^T
|\bar\rho(t)/c|^{1/2}\,dt\bigr)^2$ of the normalized local error
density;\ smaller values indicate a sharper time-localization of the
discretization error and a correspondingly larger efficiency margin
for an equidistributing adaptive mesh over a uniform one.\ }%
\RTadd{AR3-elab}{The constant $C$ is moreover not a problem-dependent
fudge factor:\ the pair $(c, C)$ is the \emph{same structural pair of
constants} that threads through all three guarantees
of~\cite{KarlssonLarssonSandbergTempone2015}:\
(i)~the \emph{stopping rule} -- the per-interval threshold
$\varepsilon_{\mathrm{time}}/N$ used in Algorithm~2.6 (above) and the
safety margin $s\le c/M_{\mathrm{sub}}$ are calibrated against the
same equidistribution principle that yields the inequalities below;\
(ii)~the \emph{accuracy} guarantee
$\limsup_{\mathrm{TOL}\to 0}\mathrm{TOL}^{-1}|\text{error}|\le 1$ of
\cite[Theorem~2.8]{KarlssonLarssonSandbergTempone2015}, which uses the
same density-ratio bound~$c$ to control how the actual symplectic-Euler
error relates to the user-prescribed tolerance;\ and
(iii)~the \emph{efficiency} guarantee
$\mathrm{TOL}\cdot N\le C\|\bar\rho/c\|_{L^{1/2}}$
of \cite[Theorem~2.9]{KarlssonLarssonSandbergTempone2015}, which is
sharp in the asymptotic limit:\
$\lim_{\mathrm{TOL}\to 0^+}\mathrm{TOL}\cdot N=\|\bar\rho/c\|_{L^{1/2}}$
for the optimal equidistributing adaptive mesh, i.e.\ $C\to 1$.
The closer the initial mesh is to the asymptotic regime, the closer
both $c$ and $C$ can be chosen to $1$; picking $(c, C)$ a priori on a
coarse initial mesh is, however, not always obvious. \emph{Good
practice} is therefore to track the empirical parent-to-children
quotients of the local error density across each subdivision pass and
confirm that they stay within the $(c, C)$ used in the stopping and
refinement rules above;\ if the empirical quotients exceed the assumed
bound the refinement loop is operating outside the regime in which
Theorems 2.7--2.9 are tight, and either the initial mesh should be
refined or the assumed $(c, C)$ relaxed.
In short, $c$ and $C$ are not independent decorations of one inequality
but the very constants that the refinement decision, the stopping rule,
and the resulting accuracy / efficiency bounds all share.\ }%
In the present
solver, the outer loop (Section~\ref{ssec:outer-loop}) modifies not
only the time mesh but also the Hamiltonian model via bundle
enrichment and the smoothing level~$\delta$. Consequently, the
combined three-source adaptive scheme is an \emph{algorithmic
extension} of \cite{KarlssonLarssonSandbergTempone2015}, and should be read as a heuristic
guided by, but not directly covered by, their theorems. Its
effectiveness is demonstrated empirically in
Section~\ref{sec:numerics}.

% ------------------------------------------------------------------
\subsection{Adaptive outer loop}\label{ssec:outer-loop}

The three indicators are controlled by a single outer loop that
refines the time mesh, the PA bundle, and the smoothing level in a
prioritized sequence.

\RTadd{R70-B3}{\emph{Inactive-layer convention.}\;
If a layer is inactive by construction, its indicator-to-tolerance ratio is set to
zero and its branch is skipped. In particular, in analytic-smoothing
mode $\eta_{\mathrm{PA}}\equiv 0$ structurally and we set
$r_{\mathrm{PA}}=0$ rather than dividing by
$\varepsilon_{\mathrm{PA}}=0$ in Algorithm~\ref{alg:outer};\ in
exact-smooth mode both $r_{\mathrm{PA}}$ and $r_\delta$ are set to
zero. The ratio definitions in Algorithm~\ref{alg:outer} are read as
piecewise, with the inactive branches returning zero.}

\begin{algorithm}[H]
\caption{Adaptive outer loop}\label{alg:outer}
\begin{algorithmic}
\Require Initial mesh $\{t_n\}_{n=0}^{N_0}$, initial PA bundle
 $\bar H_0$ with $M_0$ planes, initial smoothing
 $\delta_0>0$, tolerances
 $\varepsilon_{\mathrm{time}},\varepsilon_{\mathrm{PA}},
 \varepsilon_\delta$, overload parameter $\beta_{\mathrm{bal}}\in(0,1)$\RTadd{R73-M5}{, time-balance demotion threshold $\beta_{\mathrm{tb}}\in(0,1)$}
\Ensure Converged solution $(x,p)$ with all indicators below
 tolerances
\State $\delta\leftarrow\delta_0$;\;
 $\bar H\leftarrow\bar H_0$;\; mesh $\leftarrow\{t_n\}$
\While{not converged}
 \State Build $H_\delta$ from current $\bar H$ and $\delta$
 via~\eqref{eq:H_delta}
 \State Solve TPBVP by Algorithm~\ref{alg:newton}
 $\to(x^{(k)},p^{(k)})$
 \State Compute $\eta_{\mathrm{time}}^{\max}$, $\eta_{\mathrm{PA}}$,
 $\eta_\delta$ via~\eqref{eq:eta_time}--\eqref{eq:eta_delta}
 \State $r_{\mathrm{time}}\leftarrow
 \eta_{\mathrm{time}}^{\max}/(\varepsilon_{\mathrm{time}}/N)$;\;
 \RTadd{R70-B3}{$r_{\mathrm{PA}}\leftarrow
 \eta_{\mathrm{PA}}/\varepsilon_{\mathrm{PA}}$ if the PA layer is
 active, else $r_{\mathrm{PA}}\leftarrow 0$};\;
 \RTadd{R70-B3}{$r_\delta\leftarrow\eta_\delta/\varepsilon_\delta$
 if the smoothing layer is active, else
 $r_\delta\leftarrow 0$}
 \If{$\max(r_{\mathrm{time}},r_{\mathrm{PA}},r_\delta)\le 1$}
 \State $\mathit{converged}\leftarrow\mathbf{true}$ \Else
 \State \textit{(time-balance guard: priority only, never affects stop)}
 \If{$\eta_{\mathrm{time}}^{\max}\le
 \RTadd{R73-M5}{\beta_{\mathrm{tb}}}\cdot\max(\eta_{\mathrm{PA}},\eta_\delta)$}
 \State $r_{\mathrm{time}}^{\mathrm{adm}}\leftarrow 0$
 \Else
 \State $r_{\mathrm{time}}^{\mathrm{adm}}\leftarrow r_{\mathrm{time}}$
 \EndIf
 \State \RTadd{R66-rho}{$r_{\max}$}$\leftarrow\max(r_{\mathrm{time}}^{\mathrm{adm}},
 r_{\mathrm{PA}},r_\delta)$
 \If{$r_{\mathrm{time}}^{\mathrm{adm}}\ge\beta_{\mathrm{bal}}\RTadd{R66-rho}{r_{\max}}$}
 \State Refine mesh (Section~\ref{ssec:adaptive-mesh})
 \EndIf
 \If{$r_{\mathrm{PA}}\ge\beta_{\mathrm{bal}}\RTadd{R66-rho}{r_{\max}}$}
 \State Enrich PA bundle via Algorithm~\ref{alg:pa-enrich}
 \EndIf
 \If{$r_\delta\ge\beta_{\mathrm{bal}}\RTadd{R66-rho}{r_{\max}}$}
 \State $\delta\leftarrow\delta/2$
 \EndIf
 \EndIf
\EndWhile
\end{algorithmic}
\end{algorithm}

\RTadd{R53-B1}{\emph{Note.} The first \textbf{if} is an explicit stopping
test on the indicator-to-tolerance ratios; the loop returns \emph{converged} only when
\emph{all three} indicators are at or below their tolerances. The
time-balance guard rescales the time priority $r_{\mathrm{time}}$ to its
admissible value $r_{\mathrm{time}}^{\mathrm{adm}}$ for the purpose of
ranking layers, but it never affects the stopping test. The
\textbf{else} branch refines every layer whose indicator-to-tolerance ratio is within
a factor $\beta_{\mathrm{bal}}$ of the dominant overload, so several
layers may be refined in a single outer iteration when their ratios are
comparable; the next \textbf{while}-iteration automatically re-solves
the TPBVP with the updated mesh, bundle, or~$\delta$.}

\RTadd{R53-B1}{%
\paragraph{Time-balance guard.}
\RTadd{R73-M5}{Algorithm~\ref{alg:outer} carries two conceptually
distinct parameters in $(0,1)$:\ the balance band
$\beta_{\mathrm{bal}}$ of the refinement step (how close a layer's
indicator-to-tolerance ratio must be to the dominant one to be refined in the
current outer iteration), and the time-balance demotion threshold
$\beta_{\mathrm{tb}}$ (how far below the non-time indicators the time
indicator must be before its priority is demoted).\ We use a common
default $\beta_{\mathrm{bal}}=\beta_{\mathrm{tb}}=0.1$ across the
benchmarks of~\S\ref{sec:numerics};\ the two parameters can be tuned
independently.}
The guard $\eta_{\mathrm{time}}^{\max}\le
\RTadd{R73-M5}{\beta_{\mathrm{tb}}}\cdot\max(\eta_{\mathrm{PA}},\eta_\delta)$
uses $\beta_{\mathrm{tb}}\in(0,1)$ (default
\RTadd{R73-M5}{$\beta_{\mathrm{tb}}=0.1$}\Ksadd{12}{; the value $0.1$ encodes
``one order of magnitude'' --- the smallest gap at which mesh
refinement is unambiguously dominated by the non-time error --- and is
used unchanged across the benchmarks of~\S\ref{sec:numerics}. Problems
where the time indicator decays much faster than the bundle/smoothing
indicators may benefit from a smaller $\beta_{\mathrm{tb}}$
(e.g.\ $10^{-2}$); problems where the three indicators decay at
comparable rates from a larger $\beta_{\mathrm{tb}}$
(e.g.\ $0.3$)}) to demote the time priority whenever the time indicator
is already at least one order of magnitude below the dominant non-time
indicator. Without this guard, the outer loop can over-refine the mesh
while the modelling or smoothing error is the actual bottleneck,
because $\eta_{\mathrm{time}}^{\max}$ is max-type and
$\eta_{\mathrm{PA}},\eta_\delta$ are sum-type: the same absolute
tolerance on each can leave the max-type quantity nominally above
threshold even when its contribution to the total error is negligible.
The guard is skipped when $\max(\eta_{\mathrm{PA}},\eta_\delta)=0$,
which is the explicit-gradient mode of
Section~\ref{ssec:explicit-grad}. \Rzzedit{B4-pol}{The guard rescales only}
the priority ratio used to decide \emph{which} layer to refine; the
stopping test remains the unconditional inequality
$\max(r_{\mathrm{time}},r_{\mathrm{PA}},r_\delta)\le 1$ on the
un-rescaled time ratio.%
}

\RTadd{R66-B2}{%
Figure~\ref{fig:outer-flow} renders Algorithm~\ref{alg:outer} as a
flowchart. Each outer iteration solves the TPBVP on the current
$(N,M,\delta)$, computes the three indicator-to-tolerance ratios
$(r_{\mathrm{time}},r_{\mathrm{PA}},r_\delta)$, and either stops
(if $\max(r_{\mathrm{time}},r_{\mathrm{PA}},r_\delta)\le 1$) or
refines, in parallel, every layer whose admissible indicator-to-tolerance ratio is
within a factor $\beta_{\mathrm{bal}}$ of the dominant ratio.}

\RTadd{R66-B2}{%
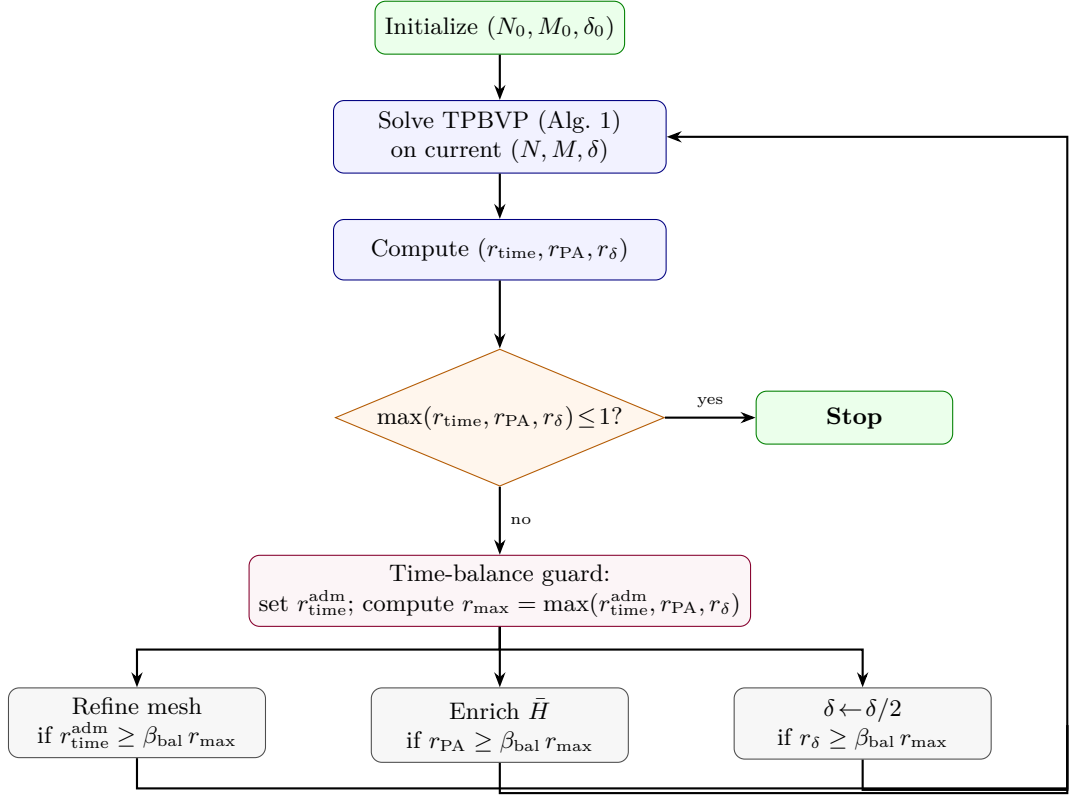
\begin{figure}[htbp]
\centering
\begin{tikzpicture}[>=Stealth, font=\footnotesize,
 node distance=6mm and 12mm,
 proc/.style={rectangle, rounded corners, draw=blue!50!black,
 fill=blue!5, minimum width=44mm, minimum height=8mm,
 align=center},
 dec/.style={diamond, draw=orange!70!black, fill=orange!7,
 aspect=2.4, inner sep=0.5mm, align=center},
 term/.style={rectangle, rounded corners, draw=green!50!black,
 fill=green!8, minimum width=26mm, minimum height=7mm,
 align=center},
 act/.style={rectangle, rounded corners, draw=black!70,
 fill=black!3, minimum width=34mm, minimum height=7mm,
 align=center},
 parallel/.style={rectangle, rounded corners, draw=purple!70!black,
 fill=purple!4, minimum width=44mm, minimum height=8mm,
 align=center},
 arr/.style={-{Stealth[length=2mm]}, thick}]
 \node[term] (start) {Initialize
 $(N_0,M_0,\delta_0)$};
 \node[proc, below=of start] (solve) {Solve TPBVP
 (Alg.~\ref{alg:newton})\\ on current $(N,M,\delta)$};
 \node[proc, below=of solve] (eval) {Compute
 $(r_{\mathrm{time}},r_{\mathrm{PA}},r_\delta)$};
 \node[dec, below=9mm of eval] (dstop)
 {$\max(r_{\mathrm{time}},r_{\mathrm{PA}},r_\delta)\!\le\!1$?};
 \node[term, right=12mm of dstop] (stop) {\textbf{Stop}};
 \node[parallel, below=9mm of dstop] (guard)
 {Time-balance guard:\\
 set $r_{\mathrm{time}}^{\mathrm{adm}}$;\
 compute $r_{\max}=\max(r_{\mathrm{time}}^{\mathrm{adm}},r_{\mathrm{PA}},r_\delta)$};
 \node[act, below=8mm of guard, xshift=-48mm] (atime)
 {Refine mesh\\ if $r_{\mathrm{time}}^{\mathrm{adm}}\ge\beta_{\mathrm{bal}}\,r_{\max}$};
 \node[act, below=8mm of guard] (aPA)
 {Enrich $\bar H$\\ if $r_{\mathrm{PA}}\ge\beta_{\mathrm{bal}}\,r_{\max}$};
 \node[act, below=8mm of guard, xshift=48mm] (ad)
 {$\delta\!\leftarrow\!\delta/2$\\ if $r_\delta\ge\beta_{\mathrm{bal}}\,r_{\max}$};
 \draw[arr] (start) -- (solve);
 \draw[arr] (solve) -- (eval);
 \draw[arr] (eval) -- (dstop);
 \draw[arr] (dstop) -- node[above, font=\tiny]{yes} (stop);
 \draw[arr] (dstop) -- node[right, font=\tiny]{no} (guard);
 \draw[arr] (guard.south) -- ++(0,-3mm) -| (atime.north);
 \draw[arr] (guard.south) -- (aPA.north);
 \draw[arr] (guard.south) -- ++(0,-3mm) -| (ad.north);
 % parallel return path
 \coordinate (ret) at ([xshift=58mm]aPA.east);
 \draw[thick] (atime.south) -- ++(0,-4mm) -| (ret);
 \draw[thick] (aPA.south) -- ++(0,-4mm) -| (ret);
 \draw[thick] (ad.south) -- ++(0,-4mm) -| (ret);
 \draw[arr] (ret) |- (solve.east);
\end{tikzpicture}
\caption{Balanced adaptive outer loop (Algorithm~\ref{alg:outer}).\
The stopping test is the unconditional overload-ratio inequality
$\max(r_{\mathrm{time}},r_{\mathrm{PA}},r_\delta)\le 1$ on the raw
ratios.\ When it fails, the time-balance guard sets the admissible
time ratio $r_{\mathrm{time}}^{\mathrm{adm}}$ for prioritization only
(it never affects the stop), and every layer whose (admissible)
indicator-to-tolerance ratio is within a factor $\beta_{\mathrm{bal}}$ of the
dominant ratio $r_{\max}$ is refined in the same outer iteration
before the next TPBVP solve.\ \RTadd{R66-rho}{The symbol $r_{\max}$ is
used here for the dominant indicator-to-tolerance ratio to avoid clash with the
time-error density $\rho_n$ of~\S\ref{ssec:indicators}.}}
\label{fig:outer-flow}
\end{figure}%
}

\RTadd{R62-B4}{The controller is not a strict time-first priority
chain.\ The unconditional stop test uses the raw indicator-to-tolerance ratios:\ the
loop returns \emph{converged} only when
$\max(r_{\mathrm{time}},r_{\mathrm{PA}},r_\delta)\le 1$.\ When the stop
test fails, the time-balance guard rescales the time priority to its
admissible value $r_{\mathrm{time}}^{\mathrm{adm}}$ for ranking purposes
only;\ all layers whose (admissible) indicator-to-tolerance ratio lies within a
factor $\beta_{\mathrm{bal}}$ of the dominant ratio are refined in the
same outer iteration. Each refinement stage warm-starts the Newton
solver from the interpolated solution on the previous mesh, which
typically provides a good initial guess and reduces the iteration
count.}

\paragraph{Computational cost.}
One Newton iteration comprises: (i)~a forward sweep
($\mathcal{O}(Nd)$ operations), (ii)~residual and Jacobian assembly
($\mathcal{O}(NMd^2)$ for the soft-min derivatives with $M$~planes),
and (iii)~a block-tridiagonal solve
($\mathcal{O}(Sd^3)$ for $S$~shooting intervals). For fixed $d$ and
moderate~$M$, the dominant cost per outer iteration is
$\mathcal{O}(N)$, and the total work is controlled by the product
of outer iterations, Newton steps per solve, and the final mesh
size~$N$.

\MSadd{\paragraph{Post-convergence recovery of the admissible trajectory.}
The outer adaptive loop drives the three value-error indicators below tolerance and returns a
converged costate \RTadd{R62-M1}{$p_h$} together with the discrete state $x_h$ produced by the symplectic
Euler TPBVP solve. Under the value-first scope of~\S\ref{sec:introduction}, the state $x_h$
from the joint Newton solve is treated as a means rather than as a separately-validated output;
\RTadd{R62-B10}{the value-first scope calls for a single post-convergence plug-and-integrate pass at
the converged $p_h$ to obtain an admissible \emph{admissible-trajectory} pair
$(x_h^{\mathrm{adm}},\alpha_h)$:}
\RTadd{R70-B4}{(R1) define the control selector
$\widehat\alpha(t,x,p_h(t)) \in \arg\min_{\alpha\in\mathcal{A}}\bigl\{p_h(t)^\top f(x,\alpha,t)+\ell(x,\alpha,t)\bigr\}$
(or its smoothed surrogate when the exact $\arg\min$ is set-valued at the Borel tie-set);
(R2) integrate $\dot x(t) = f(x(t),\widehat\alpha(t,x(t),p_h(t)),t)$ with $x(0)=x_0$ on the
converged time mesh to obtain the admissible state $x_h^{\mathrm{adm}}$;
(R3) set $\alpha_h(t) := \widehat\alpha(t,x_h^{\mathrm{adm}}(t),p_h(t))$.}
The pair $(x_h^{\mathrm{adm}},\alpha_h)$ is admissible by construction and provides the
\RTadd{R53-M1}{admissible trajectory for the upper bound \RTadd{R70-B4}{$J^{\mathrm{pr}}(\alpha_h)\ge U(0,x_0)$},
treated as a primal admissible value. No pathwise error bound on
$(x_h^{\mathrm{adm}},\alpha_h)$ is claimed; the value-error budget on $|J_h-U(0,x_0)|$
is supported by the three diagnostics $\eta_{\mathrm{time}}$, $\eta_{\mathrm{PA}}$,
$\eta_\delta$ above. This work does not construct a globally validated
lower verification function and therefore does not certify a primal--dual interval.}\
\RTadd{R62-B10}{The post-convergence plug-and-integrate pass (R1)--(R3) is the
value-first scope's prescription for recovering a admissible trajectory
from the converged costate;\ the present implementation returns the
joint Newton TPBVP triple $(x_h,p_h,\alpha_h^{\mathrm{Newton}})$, and the
\Rzzedit{A1d}{post-convergence witness $(x_h^{\mathrm{adm}},\alpha_h)$ is computed as a
post-processing step on the converged costate whenever an admissible witness is required.}}}

% ------------------------------------------------------------------
\MSadd{%
\subsection{PA bundle enrichment}\label{ssec:pa-enrich}

At an outer iteration in which the PA branch is taken, the bundle
$\bar H$ is grown by a problem-data-driven number of new planes rather
than by one plane per iteration. The enrichment
\MSdel{algorithm}\MSadd{Algorithm~\ref{alg:pa-enrich}} ranks
candidate nodes by a \emph{time-weighted} score and filters them by a
local time-separation radius and a ranking-score floor.\MSadd{\
The idea of Algorithm~\ref{alg:pa-enrich} is detailed in the
following steps~(P1)--(P5).}

\begin{enumerate}
\item[(P1)] \emph{Per-node PA gap and score.}\;
 At each node $n\in\{0,\dots,N-1\}$ compute the local PA gap
 \[
 g_n \;:=\; \bar H(p_{n+1},x_n,t_n) - H(p_{n+1},x_n,t_n)\;\ge\;0,
 \]
 and the interval-weighted PA score
 \begin{equation}\label{eq:pa-score}
 s_n \;:=\; g_n\,\Delta t_n.
 \end{equation}
\item[(P2)] \emph{Batch size.}\;
 Let $M$ be the current bundle size and
 $\theta_{\mathrm{PA}}\in(0,1]$ a batch fraction (default
 $\theta_{\mathrm{PA}}=0.1$). The target number of new planes this
 iteration is
 \begin{equation}\label{eq:pa-kadd}
 k_{\mathrm{add}}\;:=\;
 \max\bigl(1,\lceil\theta_{\mathrm{PA}}\,M\rceil\bigr).
 \end{equation}
 The batch fraction $\theta_{\mathrm{PA}}$ is held constant in this
 work; an example-adaptive choice tied to the current bundle size
 $M$, the smoothing scale $\delta$, and rate-of-convergence
 diagnostics is a natural future-work extension.
\item[(P3)] \emph{Score floor.}\;
 Only candidates with $s_n\ge\tau\,\max_m s_m$ for a floor fraction
 $\tau\in[0,1]$ (default $\tau=0.2$) are eligible for enrichment;
 this prevents the batch from being padded with nodes whose
 contribution to $\eta_{\mathrm{PA}}$ is negligible.
\item[(P4)] \emph{Greedy selection with time-separation filter.}\;
 Sort eligible nodes by decreasing $s_n$. Initialize the set of
 selected indices $\mathcal{S}:=\varnothing$. For each candidate
 $n$ in order, accept~$n$ if it is not too close in time to a
 previously selected node, i.e.\ if
 \begin{equation}\label{eq:pa-sep}
 |t_n-t_m|\;>\;c_{\mathrm{sep}}\,\max(\Delta t_{n-1},\Delta t_n)
 \qquad\forall\,m\in\mathcal{S},
 \end{equation}
 where $c_{\mathrm{sep}}\ge 1$ is the time-separation factor
 (default $c_{\mathrm{sep}}=5$). Stop when
 $|\mathcal{S}|=k_{\mathrm{add}}$ or no further candidate satisfies
 (P3)--(P4).
\item[(P5)] \emph{Bundle update.}\;
 For each $n\in\mathcal{S}$ call the oracle at $(p_{n+1},x_n,t_n)$,
 receive \Rzzedit{S22}{$\bigl(H,g,a^\star\bigr)$ in the notation of
 \S\ref{ssec:oracle}, where the supergradient $g\in\partial_pH(\cdot,x,t)|_p$
 equals $\nabla_p H$ at the queried point and equivalently coincides
 with the controlled drift $f(x,a^\star,t)$}, and insert the
 corresponding affine majorant into the bundle via the coefficient
 extension of~\S\ref{ssec:extension}. Duplicate controls
 (same~$a^\star$ as an existing bundle entry) are
 rejected to keep the bundle full-rank.
\end{enumerate}

\begin{algorithm}[H]
\caption{PA bundle enrichment}\label{alg:pa-enrich}
\begin{algorithmic}
\Require Current bundle $\bar H$ with $M$~planes;\;
 trajectory $(x_n,p_n)_n$ and mesh $\{\Delta t_n\}$;\;
 parameters $\theta_{\mathrm{PA}}$, $\tau$, $c_{\mathrm{sep}}$
\Ensure Enriched bundle $\bar H'$ with $M'\ge M$~planes
\State Compute $\{s_n\}$ by~\eqref{eq:pa-score}
\State $k_{\mathrm{add}}\leftarrow
 \max(1,\lceil\theta_{\mathrm{PA}}M\rceil)$
\State $s_{\max}\leftarrow\max_n s_n$;\;
 $\mathcal{C}\leftarrow\{n:s_n\ge\tau\,s_{\max}\}$
\State Sort $\mathcal{C}$ by decreasing $s_n$;\;
 $\mathcal{S}\leftarrow\varnothing$
\For{$n\in\mathcal{C}$ in sorted order}
 \If{$|\mathcal{S}|\ge k_{\mathrm{add}}$}
 \State \textbf{break}
 \EndIf
 \If{$|t_n-t_m|>c_{\mathrm{sep}}\,\max(\Delta t_{n-1},\Delta t_n)$
 for all $m\in\mathcal{S}$}
 \State $\mathcal{S}\leftarrow\mathcal{S}\cup\{n\}$
 \EndIf
\EndFor
\For{$n\in\mathcal{S}$}
 \State Call oracle at $(p_{n+1},x_n,t_n)$;\;
 insert affine majorant into $\bar H$ (skip duplicate controls)
\EndFor
\State \Return updated bundle
\end{algorithmic}
\end{algorithm}

\paragraph{Why time-weighting matters.}
A naive enrichment rule that ranks candidates by the raw nodal gap
$g_n$ alone can miss regions that contribute strongly to the sum
$\eta_{\mathrm{PA}}=\sum_n g_n\Delta t_n$ only through long intervals.
Ranking by~\eqref{eq:pa-score} aligns the enrichment rule with the
quantity that actually drives the stopping criterion, and was the
single largest source of the ``plots say this but the code does that''
mismatch in the earlier implementation (cf.\ the diagnostics listed in
Section~\ref{sec:numerics}).

\paragraph{Why batching matters.}
Single-plane-per-iteration enrichment produces long sequences of
nearly identical outer iterations when the bundle is small relative
to the problem difficulty. The rule $k_{\mathrm{add}}\propto M$ lets
the bundle grow proportionally, while (P3)--(P4) prevent the batch
from collapsing onto a single tight cluster of neighbouring nodes.%
}%
\MSnote{2026-04-24 update, \S3.4--3.7: replace the implicit
one-plane-per-iteration rule by a fully specified PA enrichment
subroutine with interval-weighted score, batch size, score floor,
and time-separation filter.}

% ------------------------------------------------------------------
\MSadd{%
\subsection{Solver modes}\label{ssec:explicit-grad}

The Newton solve and the three-indicator outer loop are written
generically for the smoothed Hamiltonian $H_\delta$ built from a PA
bundle, but in practice the proposed solver is operated in one of
three concrete modes that differ in what plays the role of $H_\delta$
and which indicators are then nonzero. The mode is a property of how
the user supplies the Hamiltonian, not of the solver. We tag the
mode of every numerical example in Section~\ref{sec:numerics}.

\begin{enumerate}
\item \textbf{Exact-smooth mode.}\;
 $H$ is given analytically and is already smooth enough for the
 Newton solve. The solver sets $H_\delta\equiv H$, bypasses both
 the PA bundle and the soft-min smoothing, and drives adaptation
 with the time indicator only. By construction
 $\eta_{\mathrm{PA}}\equiv 0$ and $\eta_\delta\equiv 0$, so the
 corresponding branches of Algorithm~\ref{alg:outer} are skipped and
 the time-balance guard is skipped (because
 $\max(\eta_{\mathrm{PA}},\eta_\delta)=0$). The outer loop reduces
 to the single-indicator adaptive scheme
 of~\cite{KarlssonLarssonSandbergTempone2015}.
\item \textbf{Analytic-smoothing mode.}\;
 $H$ is given analytically but is non-smooth in the costate (or in
 another argument), and a problem-specific smoothing
 $H_\delta^{\mathrm{an}}$ is supplied in closed form. The Newton
 system uses $H_\delta^{\mathrm{an}}$, and the PA bundle is again
 bypassed: $\eta_{\mathrm{PA}}\equiv 0$ structurally, while
 $\eta_\delta>0$ in general because $H_\delta^{\mathrm{an}}\ne H$.
 An initial PA bundle $\bar H$ is still constructed (typically
 collapsing to the exact $H$ in this mode, as in
 \S\ref{ssec:ex-nonsmooth}) so that
 $\eta_\delta = \mathrm{Quad}[\bar H - H_\delta^{\mathrm{an}}]$
 remains well-defined and drives the $\delta$-adaptation; the
 bundle is held fixed across outer iterations because the PA layer
 is structurally inactive.
 Equivalently, this is exact-smooth mode applied to
 $H_\delta^{\mathrm{an}}$ rather than to $H$ itself, with an extra
 smoothing indicator $\eta_\delta$ that controls the analytic
 smoothing bias. \RTadd{R66-S1}{Example~2}
 (Section~\ref{ssec:ex-nonsmooth}) is the prototype:\
 $H(p,x)=x^{10}-|p|$ is replaced by
 $H_\delta^{\mathrm{an}}(p,x)=x^{10}-\sqrt{p^2+\delta^2}$.
\item \textbf{PA-softmin mode.}\;
 The default mode of this work. The Hamiltonian is approached
 through a PA bundle $\bar H$ built by the oracle, and
 $H_\delta$ is the soft-min smoothing~\eqref{eq:H_delta} of that
 bundle. The initial set of PA planes is example-specific and aims
 to cover the parts of the costate range relevant for that
 benchmark (e.g.\ the bang-off scalar in~\S\ref{ssec:ex-nonsmooth}
 collapses to the three-plane set $\{-1, 0, +1\}$, while the
 QP-oracle row of~\S\ref{ssec:ex-B} starts from $M_0 = 8$ planes
 seeded by oracle calls at the bootstrap mesh). Both $\eta_{\mathrm{PA}}\ge 0$ and $\eta_\delta\ge 0$ are
 generically nonzero and the full three-indicator outer loop is
 active. The LQR calibration run (Appendix~\ref{apx:lqr}) and the
 \RTadd{R62-M4}{hypersensitive Example~1
 (Section~\ref{ssec:ex-hypersensitive})}
 are run in this mode. \RTadd{R70-m4}{A degenerate finite-branch
 subcase occurs when the oracle returns only finitely many distinct
 minimizers along the trajectory:\ once the bundle contains all of
 them, the measured PA gap becomes structurally zero
 ($\eta_{\mathrm{PA}}\equiv 0$) even though the PA-softmin
 machinery remains active.}
\end{enumerate}

\noindent
The verification-driven taxonomy clarifies a recurring confusion in earlier
drafts:\ ``$\eta_{\mathrm{PA}}=0$'' is a property of the bundle's
content, not of the solver mode. In exact-smooth and
analytic-smoothing modes it is $\equiv 0$ structurally because no
bundle is constructed;\ in PA-softmin mode it can also be zero
\emph{measured}, namely whenever the bundle has caught up with the
finite set of optimal controls actually visited along the
trajectory.}%
%% === end of Algorithm_v2.tex ===

\section{Numerical Experiments}\label{sec:numerics}

%% === inlined from Numexperiments_v2.tex (R98z-cc arxiv-flat) ===
% ------------------------------------------------------------------
% §4 Numerical Experiments
% ------------------------------------------------------------------
%
% Numerical Experiments section: four validated unconstrained-state
% benchmark runs (LQR, hypersensitive, nonsmooth scalar, singular
% tracking). All numbers in this section are
% read off the deep-report artefacts in
% \texttt{reports/<example>\_deep\_report/tables/final\_summary.json}
% (and the Riccati / Karlsson--et-al.\ references where applicable).
%
% ------------------------------------------------------------------

We test the proposed solver on \RTadd{R53-B3}{four unconstrained-state}
fixed-horizon Bolza problems of increasing structural difficulty:
\begin{itemize}
\item \RTadd{R53-B3}{Example~1 (\S\ref{ssec:ex-hypersensitive}) -- the
 hypersensitive scalar problem of~\cite{KarlssonLarssonSandbergTempone2015},
 which exhibits a long-horizon boundary-layer structure and stresses
 adaptive time refinement.}
\item \RTadd{R53-B3}{Example~2 (\S\ref{ssec:ex-nonsmooth}) -- the scalar
 nonsmooth Hamiltonian benchmark
 of~\cite{KarlssonLarssonSandbergTempone2015}, for which the optimal
 control is bang--off and the analytic smoothing of
 Section~\ref{ssec:explicit-grad} applies.}
\item \RTadd{R53-B3}{Example~3 (\S\ref{ssec:ex-singular}) -- the singular
 tracking problem of~\cite{KarlssonLarssonSandbergTempone2015} with a
 hypergeometric reference, where strongly localized state dynamics
 drive the adaptive mesh.}
\item \RTadd{R53-B3}{Example~4 (\S\ref{ssec:ex-B}) -- a
 \emph{QP-oracle attitude control allocation} benchmark: a planar
 single-axis rigid-body tracking problem with four canted
 reaction-wheel torque commands and a constrained admissible set
 (per-wheel bounds, total-power constraint, null-space-agility
 constraint), in which the Hamiltonian minimisation is evaluated
 by a parametric-QP oracle rather than by a closed-form formula.
 This example exercises all three adaptive branches (mesh
 refinement, PA-plane insertion, and $\delta$-halving) and is the
 most direct test of the oracle-based PA-softmin architecture that
 distinguishes the proposed solver.}
\end{itemize}
\noindent
\RTadd{R53-B3}{A smooth linear--quadratic regulator (LQR) with a known
matrix-Riccati reference is retained as a calibration consistency check
and is reported in Appendix~\ref{apx:lqr}; its bundle saturates after
a few outer iterations and it does not exercise the oracle-based
PA-softmin / soft-min smoothing machinery.}
For each example we report the reference objective (where available),
the computed objective, the final mesh size~$N$, bundle size~$M$,
smoothing level~$\delta$, the Newton iteration count and residual on
the final mesh, and the values of the three a~posteriori
indicators~$\eta_{\mathrm{time}}$, $\eta_{\mathrm{PA}}$,
$\eta_\delta$ from~\S\ref{ssec:indicators}, together with the
companion sum~$\eta_{\mathrm{time,sum}}=\sum_n\bar r_n$ when
relevant.%

\paragraph{Common setup.}
Unless an example specifies otherwise, the starting values for the
first outer iteration are an initial uniform mesh of $N_0$ intervals,
an initial bundle of $M_0$ planes, and an initial smoothing
$\delta_0$;\ values are listed per example in
Table~\ref{tab:example-setup}. The Newton solver uses
$\mathrm{tol}_{\mathrm{Newton}}=10^{-10}$ on the residual norm with
$k_{\max}=50$ and a least-squares fallback after Newton stagnation.
Outer-loop tolerances are example specific and collected in
Table~\ref{tab:example-tolerances}. A single global tolerance triple
is not appropriate:\ the asymmetric nature of the three indicators
(\S\ref{ssec:indicators}) means that their absolute levels for a
given example depend on both the Bolza data and the final mesh size.
To keep the three indicators \emph{comparable in magnitude}---so that
none is effectively ignored---the solver enforces the validation rule
\begin{equation}\label{eq:tol-validate}
\varepsilon_{\mathrm{PA}}\;\le\;2\,\varepsilon_{\mathrm{time}},
\qquad
\varepsilon_\delta\;\le\;2\,\varepsilon_{\mathrm{time}},
\end{equation}
and aborts with an informative error if either condition fails. The
adaptive mesh uses a binary subdivision rule and the per-interval
threshold $\varepsilon_{\mathrm{time}}/N$ throughout
(Sections~\ref{ssec:adaptive-mesh}--\ref{ssec:outer-loop}); we adopt
the shorthand
$\mathrm{tol}_{\mathrm{time}}^\star := \varepsilon_{\mathrm{time}}/N$
in the principal tables of this section so that the comparison
$\eta_{\mathrm{time}} \lessgtr \mathrm{tol}_{\mathrm{time}}^\star$ is
visible at a glance, with the understanding that
$\eta_{\mathrm{time}}$ here always denotes the max-type stopping
quantity $\eta_{\mathrm{time}}^{\max}$ of~\eqref{eq:eta_time}.

\begin{table}[h]\centering
\caption{Initial values used by the proposed solver in
 Section~\ref{sec:numerics}. All four examples respect the
 validation rule~\eqref{eq:tol-validate}.}
\label{tab:example-setup}
\begin{tabular}{l c c c}
\toprule
Example & $N_0$ & $M_0$ & $\delta_0$ \\
\midrule
1 (Hypersensitive) & $29$ & $7$ & $2\!\times\!10^{-2}$ \\
2 (Nonsmooth scalar) & $19$ & $3$ & $2\!\times\!10^{-2}$ \\
3 (Singular tracking) & $29$ & $0$ & $1\!\times\!10^{-2}$ \\
4 (QP-oracle attitude alloc.) & $32$ & $8$ & $5\!\times\!10^{-1}$ \\
\bottomrule
\end{tabular}
\end{table}

\begin{table}[h]\centering
\caption{Outer-loop tolerances used in Section~\ref{sec:numerics}.
 \RTadd{R72-M2}{The validation rule~\eqref{eq:tol-validate}
 ($\varepsilon_{\mathrm{PA}}\le 2\varepsilon_{\mathrm{time}}$,
 $\varepsilon_\delta\le 2\varepsilon_{\mathrm{time}}$) is enforced
 only for active PA and smoothing layers;\ a dash indicates a layer
 bypassed by the solver mode (\S\ref{ssec:explicit-grad}).}
 \RTadd{R62-M3}{In the Nonsmooth row,
 $\eta_{\mathrm{PA}}\equiv 0$ structurally because no PA-bundle
 enrichment is required in the analytic-smoothing mode of
 \S\ref{ssec:explicit-grad};\ the smoothing indicator $\eta_\delta$
 remains active and its tolerance is set by the run configuration.}}
\label{tab:example-tolerances}
\begin{tabular}{l c c c}
\toprule
Example & $\varepsilon_{\mathrm{time}}$
 & $\varepsilon_{\mathrm{PA}}$
 & $\varepsilon_\delta$ \\
\midrule
1 (Hypersensitive) & $5\!\times\!10^{-2}$ & $10^{-2}$ & $10^{-2}$ \\
2 (Nonsmooth scalar) & $10^{-6}$ & \RTadd{R72-M2}{--} & set per run config \\
3 (Singular tracking) & $10^{-4}$ & \RTadd{R72-M2}{--} & \RTadd{R72-M2}{--} \\
4 (QP-oracle attitude alloc.) & $5\!\times\!10^{-3}$ & $10^{-3}$ & $10^{-3}$ \\
\bottomrule
\end{tabular}
\end{table}

\RTadd{R53-M4}{\paragraph{Run configuration.}
The per-example run configurations (Newton tolerance,
line-search safety factor~$\beta$, mesh-refinement factor, overload
parameter $\beta_{\mathrm{bal}}$, \RTadd{R73-M5}{time-balance demotion
threshold $\beta_{\mathrm{tb}}$,} PA bundle initialization, and
$\delta_0$) used to produce Tables~\ref{tab:ex2-summary}--\ref{tab:summary}
are reported alongside each example's text below. The Example~3
analytic-smoothing indicator $\eta_\delta$ is independently verified
against a recomputed integrated-bias indicator
$\eta_\delta^{\mathrm{an}} = \sum_n \Delta t_n |H(p_{n+1}, x_n, t_n)
- H_\delta^{\mathrm{an}}(p_{n+1}, x_n, t_n)|$, with agreement to
ratio $1.024$.}

% ==================================================================
% ==================================================================
\subsection{Example 1: Hypersensitive scalar problem}
\label{ssec:ex-hypersensitive}

\paragraph{Problem statement.}
We consider the scalar control problem
\SZadd{AR4}{of~\cite[\S3.1]{KarlssonLarssonSandbergTempone2015}},
\begin{equation}\label{eq:hyper-prob}
\min_{a(\cdot)} J(a)
\;=\;\gamma\bigl(x(T)-1\bigr)^2
+\int_0^T\bigl(x(t)^2+a(t)^2\bigr)\,dt,
\end{equation}
subject to $\dot x=-x^3+a$, $x(0)=1$, $a\in[a_{\min},a_{\max}]=[-1,3]$,
with $T=25$ and $\gamma=10^6$. The combination of a long horizon and
a stiff terminal penalty produces a boundary-layer structure:\ the
trajectory contracts rapidly to a quasi-steady plateau, sits there
for most of the horizon, and is then steered to satisfy the terminal
penalty in a second boundary layer. Adaptive mesh refinement is
therefore expected to be essential.

\paragraph{Pontryagin system.}
The Hamiltonian is
$H(p,x,a)=p(-x^3+a)+x^2+a^2$. Stationarity gives the
unconstrained minimizer
$a^\star_{\mathrm{unc}}(t)=-p(t)/2$, and the box-projected
optimal control is
$a^\star(t)=\Pi_{[-1,3]}(-p(t)/2)$. The adjoint equation is
$\dot p=3p x^2-2x$ with the gradient terminal condition
$p(T)=2\gamma(x(T)-1)$. The problem is smooth in both state
and costate, but its sensitivity structure is characteristic of
hypersensitive optimal control.

\paragraph{Run summary.}
With outer tolerances $\varepsilon_{\mathrm{time}}=
5\!\times\!10^{-2}$, $\varepsilon_{\mathrm{PA}}=\varepsilon_\delta=
10^{-2}$ and $\mathrm{max\_iters}=40$, the solver reaches a converged
\texttt{STOP} on a final mesh of $N=307$ intervals with $M=30$
planes and smoothing $\delta=3.13\!\times\!10^{-4}$;\ all three
indicators are below their tolerances at termination. The
terminal-state penalty is approximately saturated:\ the computed
terminal state is $x_h(T)=0.99999761$ against the target~$1$, so
the quadratic terminal penalty contributes
$\gamma(x_h(T)-1)^2\approx 5.7\!\times\!10^{-6}$ to the objective
and the integral term carries the rest. The principal results are
collected in Table~\ref{tab:ex2-summary}.

\begin{table}[htbp]
\centering
\caption{Final principal results reported by the proposed solver on
 \RTadd{R62-B2}{Example~1} (hypersensitive). ``Within tol'' refers to the
 per-interval threshold
 $\mathrm{tol}_{\mathrm{time}}^\star=\varepsilon_{\mathrm{time}}/N$.}
\label{tab:ex2-summary}
\begin{tabular}{l l}
\toprule
Quantity & Value \\
\midrule
Mode & PA-softmin \\
Terminal state $x_h(T)$ & $0.99999761$ \\
Objective $J_h$ & $2.27475268$ \\
Mesh intervals $N$ & $307$ \\
Planes $M$ & $30$ \\
Final smoothing $\delta$ & $3.13\!\times\!10^{-4}$ \\
$\eta_{\mathrm{time}}$ & $1.31\!\times\!10^{-4}$ \\
$\mathrm{tol}_{\mathrm{time}}^\star$ & $1.63\!\times\!10^{-4}$ \\
$\eta_{\mathrm{PA}}$ & $5.09\!\times\!10^{-3}$ \\
$\eta_\delta$ & $3.85\!\times\!10^{-3}$ \\
Final action & \texttt{STOP} \\
Status & accepted \\
\bottomrule
\end{tabular}
\end{table}

\begin{table}[htbp]\centering
\caption{\RTadd{R74-K10}{Branch-firing diagnostic for Example~1
 (hypersensitive). The canonical exported log records sequential
 single-action outer iterations (each non-terminal outer iteration
 fires exactly one of the three adaptive branches), in contrast with
 the packed parallel multi-branch structure of
 Table~\ref{tab:ex-B-branches} for Example~4. The row counts give
 the histogram of action types over the $27$-row log; the terminal
 row collapses the \texttt{STOP} action.}}
\label{tab:ex1-branches}
\RTadd{R74-K10}{\begin{tabular}{l c}
\toprule
Outer-iteration outcome & Count \\
\midrule
Time-mesh refinement alone & $8$ \\
PA enrichment alone        & $12$ \\
$\delta$-halving alone     & $6$ \\
Combined operations        & $0$ \\
Terminal \texttt{STOP}     & $1$ \\
\bottomrule
\end{tabular}}
\end{table}

\begin{figure}[htbp]
\centering
\begin{minipage}{0.48\textwidth}\centering
 \includegraphics[width=\textwidth]{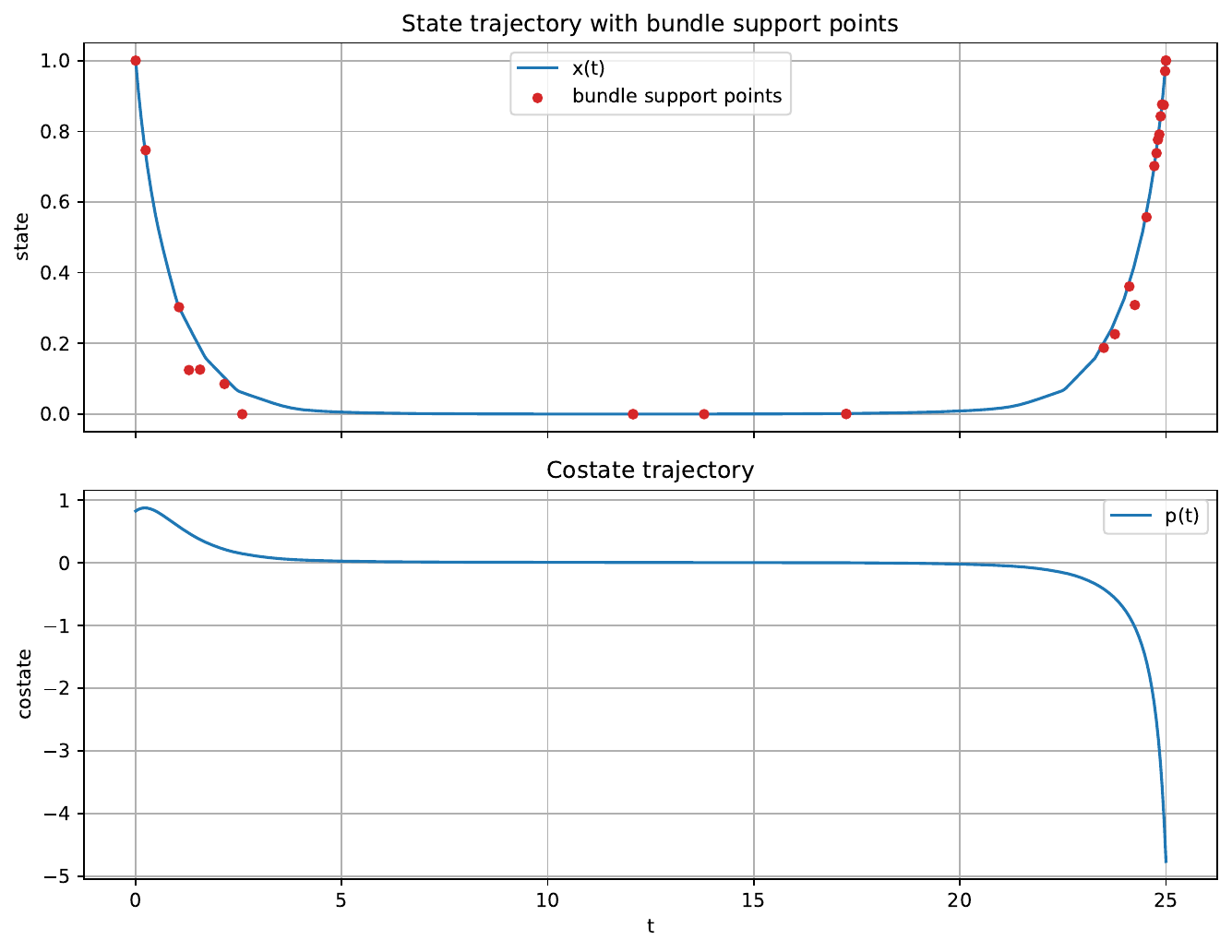}
\end{minipage}\hfill
\begin{minipage}{0.48\textwidth}\centering
 \includegraphics[width=\textwidth]{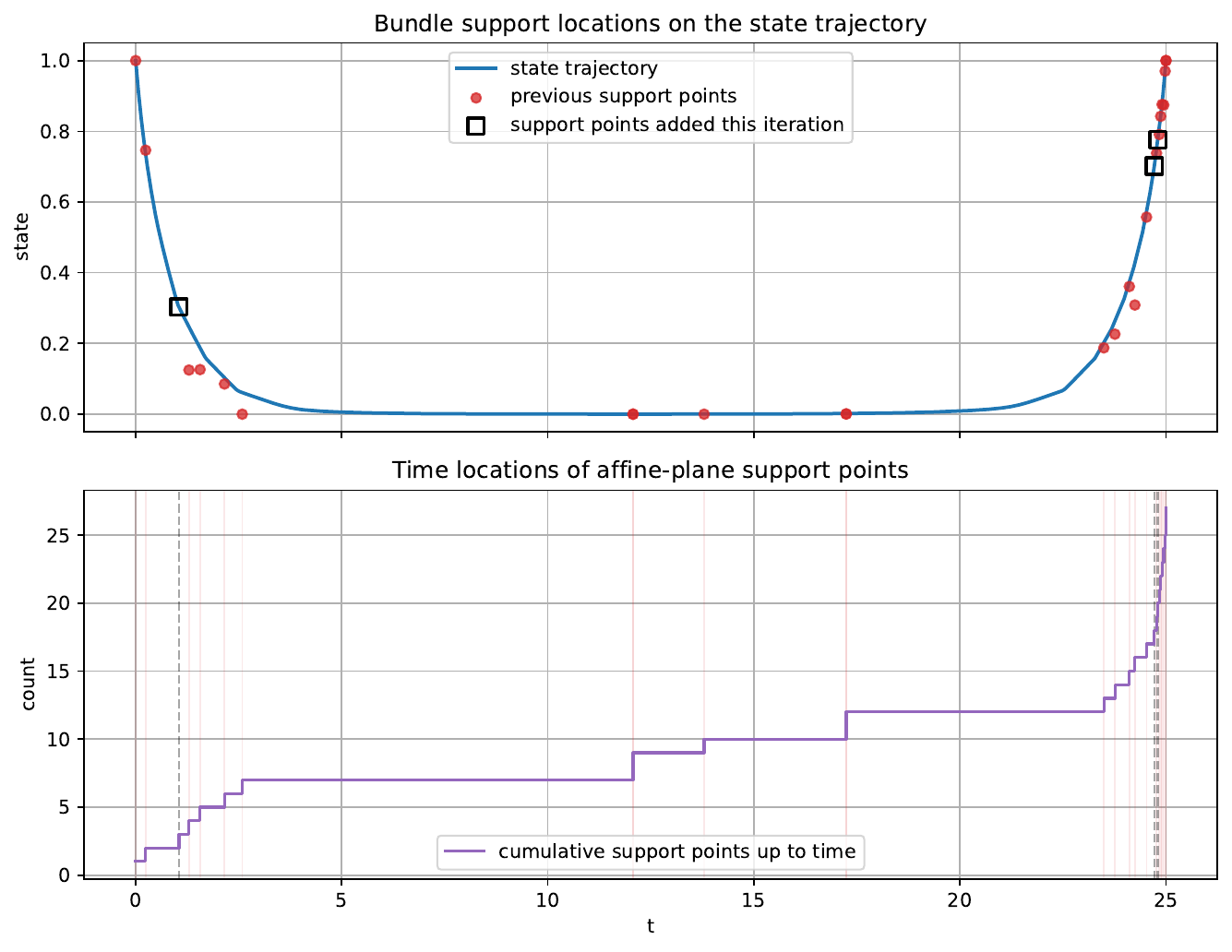}
\end{minipage}
\caption{\RTadd{R62-B2}{Example~1} (hypersensitive). Left:\ state and costate at
 the final iterate ($N=307$, $M=30$, $\delta=3.13\!\times\!10^{-4}$);\
 the three regions of the hypersensitive structure are visible.
 Right:\ bundle support points along the time axis; new affine
 planes concentrate near the two boundary layers, while the long
 interior plateau is covered by a small number of planes.}
\label{fig:ex2-bundle}
\end{figure}

\paragraph{Discussion.}
The hypersensitive example is informative for two reasons. First,
the run now reaches an accepted \texttt{STOP}:\ all three
indicators sit below their respective tolerances at termination,
and in particular the time indicator
$\eta_{\mathrm{time}}=1.31\!\times\!10^{-4}$ falls below the
tightened per-interval threshold
$\mathrm{tol}_{\mathrm{time}}^\star=1.63\!\times\!10^{-4}$ at the
final mesh size $N=307$. Second, the bundle branch remains active
throughout the outer loop, exactly as expected for this geometry:\
the boundary layers near $t=0$ and $t=T$ require additional affine
planes to track the rapid changes of the costate, while the long
quasi-steady interior is well captured by a small number of planes.
The bundle support diagnostic (Figure~\ref{fig:ex2-bundle}, right)
confirms this picture by concentrating new plane insertions near
the two boundary layers. The state and costate at the final iterate
(Figure~\ref{fig:ex2-bundle}, left) display the textbook
three-region structure:\ the initial transient at $t\approx 0$, the
long near-equilibrium plateau, and the terminal transient driven by
the penalty.

% ==================================================================
\subsection{Example 2: Nonsmooth scalar benchmark}
\label{ssec:ex-nonsmooth}
% R98z-bb Anders (Szepessy): bracketed notation-consistency note removed per his comment on the §4.2 opening sentence; the notation $a(t)$ is established directly in the equation $\dot x = a$ below.

\paragraph{Problem statement.}
We consider the scalar Bolza problem
\SZadd{AR4}{of~\cite[\S3.2]{KarlssonLarssonSandbergTempone2015}},
\begin{equation}\label{eq:nonsmooth-prob}
\min_{a(\cdot)}J(a)\;=\;\int_0^T x(t)^{10}\,dt
\qquad\text{subject to}\qquad
\dot x=a,\quad x(0)=\tfrac12,\quad a(t)\in[-1,1],
\end{equation}
on $T=1$. The problem has no terminal penalty. The exact optimal
solution is the bang--off profile
\begin{equation}\label{eq:nonsmooth-exact}
x^\star(t)=\max\!\Bigl(\tfrac12-t,\;0\Bigr),
\qquad
a^\star(t)=
\begin{cases}-1,&0\le t<\tfrac12,\\0,&\tfrac12\le t\le 1,\end{cases}
\end{equation}
with exact value
$J^\star=\frac{(1/2)^{11}}{11}\approx 4.43892\!\times\!10^{-5}$.

\paragraph{Pontryagin structure and analytic smoothing.}
The Hamiltonian (minimization convention) is
$H(p,x,a)=x^{10}+pa$, and minimization over $a\in[-1,1]$ yields the
nonsmooth reduced Hamiltonian
\begin{equation}\label{eq:nonsmooth-H}
H(p,x)\;=\;x^{10}-|p|,
\end{equation}
which is $C^\infty$ in the state and nonsmooth in the costate at
$p=0$. The PA bundle generated by oracle calls with
$a^\star(p)=-\operatorname{sgn}(p)$ collapses to the exact three-plane
set $\{-1,0,+1\}$ after a single bootstrap step:\ enrichment beyond
that set adds no information, since~\eqref{eq:nonsmooth-H} itself is
available. We retain the bundle $\bar H \equiv H$ throughout so that
$\eta_\delta = \mathrm{Quad}[\bar H - H_\delta]$ remains well-defined
and drives the $\delta$-adaptation, even though the PA channel is
structurally inactive. We exercise the explicit-gradient mode of
\S\ref{ssec:explicit-grad} with the analytic smoothing
\begin{equation}\label{eq:nonsmooth-Hdelta}
H_\delta(p,x)\;=\;x^{10}-\sqrt{p^2+\delta^2},
\qquad
\partial_pH_\delta=-\frac{p}{\sqrt{p^2+\delta^2}},
\qquad
\partial_xH_\delta=10\,x^9.
\end{equation}
Under~\eqref{eq:nonsmooth-Hdelta} the indicator
$\eta_{\mathrm{PA}}\equiv 0$ \emph{structurally} because no PA bundle
is constructed in analytic-smoothing mode. The smoothing indicator
$\eta_\delta$, by contrast, is generally nonzero unless the analytic
smoothing is exact, i.e.\ unless $H_\delta^{\mathrm{an}}=H$. In the
present \RTadd{R66-S2}{Example~2} the analytic smoothing replaces $|p|$ by
$\sqrt{p^2+\delta^2}$, so $H_\delta^{\mathrm{an}}\neq H$ and
$\eta_\delta$ remains an active indicator throughout the outer loop;
its reported value
($\eta_\delta=4.60\!\times\!10^{-8}$ at the final $\delta$,
Tables~\ref{tab:ex3-summary} and~\ref{tab:summary})
measures precisely that residual analytic-smoothing bias. With
$\eta_{\mathrm{PA}}\equiv 0$ the outer loop reduces to a two-indicator
scheme (time and $\delta$), recovering the adaptive scheme of~%
\cite{KarlssonLarssonSandbergTempone2015} as a single-indicator
specialization when $\eta_\delta$ is itself below tolerance.

\paragraph{Run summary.}
With outer tolerances
$\varepsilon_{\mathrm{time}}=10^{-6}$
($\varepsilon_{\mathrm{PA}}$ vacuously zero, no bundle; $\varepsilon_\delta$
set from the run configuration to match the integrated-bias scale
of the analytic smoothing), $\mathrm{max\_iters}=25$, and
$\delta_0=2\!\times\!10^{-2}$, the solver reaches a converged
\texttt{STOP} action at iteration $k=25$ on a mesh of
$N=275$ intervals with the smoothing $\delta=7.63\!\times\!10^{-8}$,
exactly $M=3$ planes, and the two active indicators (time and~$\delta$)
each within their respective tolerances. Principal results are listed in
Table~\ref{tab:ex3-summary}.

\begin{table}[htbp]
\centering
\caption{Final principal results reported by the proposed solver on
 \RTadd{R62-B2}{Example~2} (nonsmooth scalar).}
\label{tab:ex3-summary}
\begin{tabular}{l l}
\toprule
Quantity & Value \\
\midrule
Mode & analytic-smoothing \\
Terminal state $x_h(T)$ & $0.10117261$ \\
Objective $J_h$ & $4.47130\!\times\!10^{-5}$ \\
Exact objective $J^\star$ & $4.43892\!\times\!10^{-5}$ \\
Relative objective error & $7.29\!\times\!10^{-3}$ \\
Mesh intervals $N$ & $275$ \\
Planes $M$ & $3$ \\
Final smoothing $\delta$ & $7.63\!\times\!10^{-8}$ \\
$\eta_{\mathrm{time}}$ & $2.56\!\times\!10^{-9}$ \\
$\mathrm{tol}_{\mathrm{time}}^\star$ & $3.64\!\times\!10^{-9}$ \\
$\eta_{\mathrm{time,sum}}=\sum_n\bar r_n$ & $3.34\!\times\!10^{-7}$ \\
$\eta_{\mathrm{PA}}$ & $0$ (bundle exhausted) \\
$\eta_\delta$ & $4.60\!\times\!10^{-8}$ \\
Final action & \texttt{STOP} \\
Status & accepted \\
\bottomrule
\end{tabular}
\end{table}

\begin{table}[htbp]\centering
\caption{\RTadd{R74-K10}{Branch-firing diagnostic for Example~2
 (nonsmooth scalar). Analytic-smoothing mode keeps the bundle at
 $M=3$ throughout, so PA enrichment is structurally absent;\ the
 adaptive vocabulary collapses to $\delta$-halving and time-mesh
 refinement, alternating sequentially. The histogram is over the
 $26$-row log.}}
\label{tab:ex2-branches}
\RTadd{R74-K10}{\begin{tabular}{l c}
\toprule
Outer-iteration outcome & Count \\
\midrule
$\delta$-halving alone     & $18$ \\
Time-mesh refinement alone & $7$ \\
PA enrichment alone        & $0$ (structurally inactive) \\
Combined operations        & $0$ \\
Terminal \texttt{STOP}     & $1$ \\
\bottomrule
\end{tabular}}
\end{table}

\begin{figure}[htbp]
\centering
\begin{minipage}{0.48\textwidth}\centering
 \includegraphics[width=\textwidth]{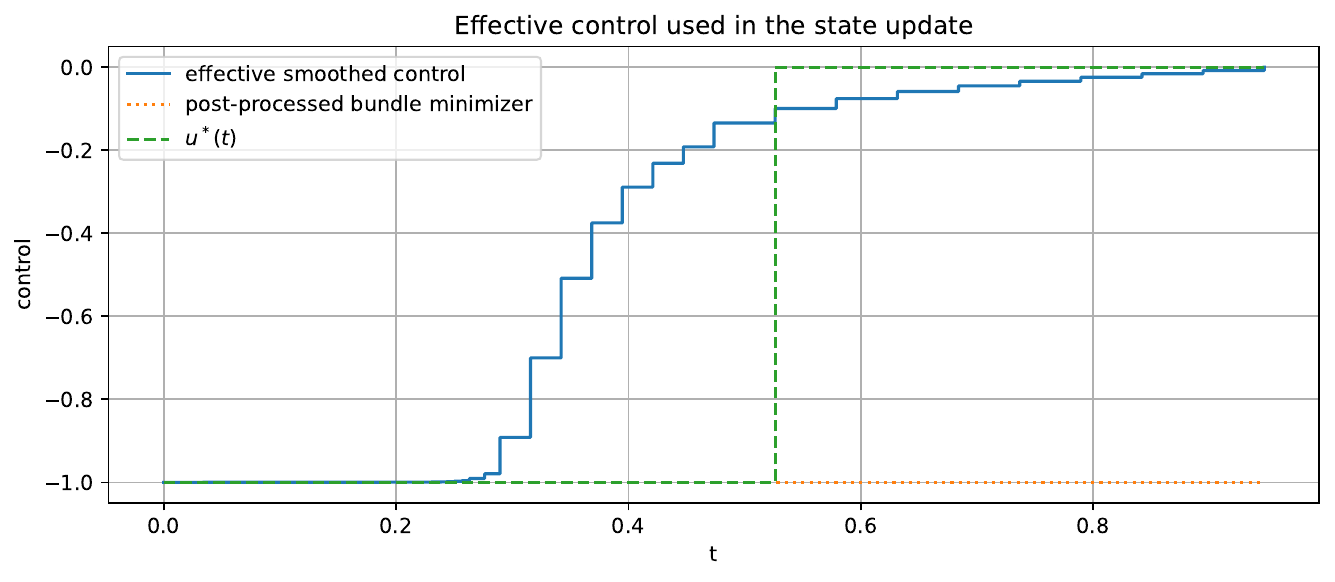}
\end{minipage}\hfill
\begin{minipage}{0.48\textwidth}\centering
 \includegraphics[width=\textwidth]{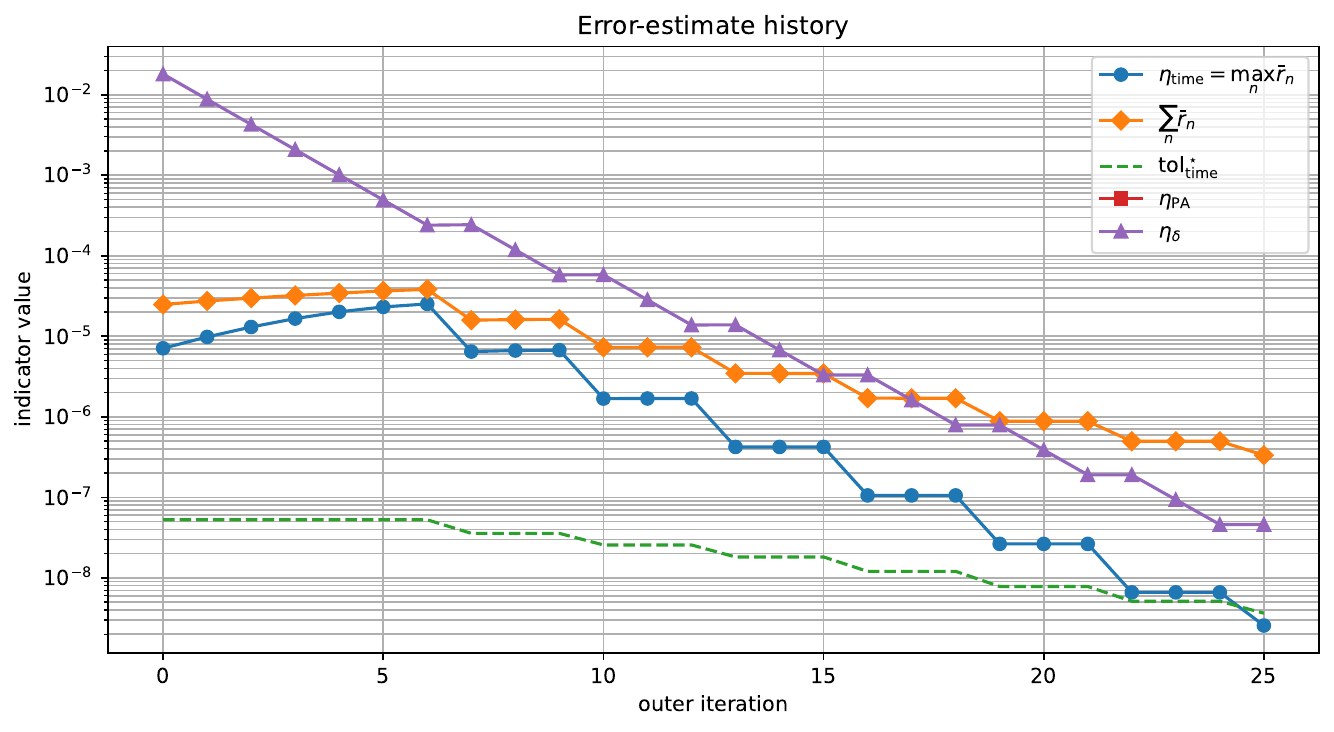}
\end{minipage}
\caption{\RTadd{R62-B2}{Example~2} (nonsmooth scalar). Left:\ effective smoothed
 control $a_\delta(t)=\partial_pH_\delta(p,x)$ used by the discrete
 state update, plotted against the post-processed bundle minimizer
 $a^\star(t)\in\{-1,0,+1\}$;\ the two coincide except in a narrow
 neighborhood of the switching time $t^\star=\tfrac12$. Right:\
 indicator history along the outer loop;\ $\eta_{\mathrm{PA}}$
 remains identically at the floor by the explicit-gradient
 reasoning, $\eta_\delta$ contracts geometrically through the
 $\delta$-halving steps, and $\eta_{\mathrm{time}}$ tracks the mesh
 refinements.}
\label{fig:ex3-control}
\end{figure}

\paragraph{Discussion.}
This benchmark separates two phenomena cleanly. First, the PA
machinery is correctly inactive: three planes already exhaust the
oracle, and the indicator $\eta_{\mathrm{PA}}$ stays identically zero
along the entire outer loop, which is the desired behavior whenever
an analytic smoothing of the exact Hamiltonian is available.
Second, even with the PA part exact, the smoothed effective control
$a_\delta(t)=\partial_pH_\delta(p(t),x(t))=
-p(t)/\sqrt{p(t)^2+\delta^2}$
that drives the discrete state equation does not coincide with the
post-processed bundle minimizer: the latter remains at $-1$ on the
active arc, while the former rounds the kink at $p=0$. This is the
reason the relative objective error remains at $7.29\!\times\!10^{-3}$
in spite of all three indicators being below tolerance:\ the residual
gap is a smoothing error in the state propagation near the switching
time $t^\star=\tfrac12$, not a stopping-rule artifact.
Figure~\ref{fig:ex3-control} (left) shows the effective smoothed
control $a_\delta$ together with the post-processed bundle minimizer,
which differ in a narrow neighborhood of $t^\star$;
Figure~\ref{fig:ex3-control} (right) shows that the indicator history
is dominated by $\eta_\delta$, with $\eta_{\mathrm{PA}}$ identically
at the floor. \RTadd{R66-S2}{Example~2} thus validates the explicit-gradient and
analytic-smoothing path of~\S\ref{ssec:explicit-grad} and shows that
in this regime the outer loop reduces, as designed, to the
single-indicator adaptive scheme
of~\cite{KarlssonLarssonSandbergTempone2015}.

% ==================================================================
\subsection{Example 3: Singular tracking}
\label{ssec:ex-singular}

\paragraph{Problem statement.}
We consider the explicitly time-dependent tracking problem
\SZadd{AR4}{of~\cite[\S3.3]{KarlssonLarssonSandbergTempone2015}},
written in autonomous form by lifting~$t$ to an auxiliary state
$s(t):=t$:
\begin{equation}\label{eq:singular-prob}
\min_{\alpha(\cdot)} J(\alpha)
\;=\;\bigl(x(T)-x_{\mathrm{ref}}(T)\bigr)^2
+\int_0^T\!\Bigl(\frac{\alpha(t)-x(t)}{\sigma(s(t))}\Bigr)^{\!2}\,dt
\end{equation}
subject to
$\dot x=\alpha/\sigma(s)$, $\dot s=1$, with the singular weight
\RTadd{R69-C5}{$\sigma(s)$ (renamed from $d(s)$ to avoid collision with
the state dimension $d$ of~\S\ref{sec:formulation} and the bundle
intercept $d_i$ of~\S\ref{ssec:PA-surrogate})}
$\sigma(s)=\bigl((s-t_0)^2+\varepsilon^2\bigr)^{\beta/2}$,
initial condition $(x(0),s(0))=(x_{\mathrm{ref}}(0),0)$, and the
parameter values
$\varepsilon=10^{-10}$, $\beta=0.75$, $t_0=\tfrac53$, $T=4$. The
exact reference trajectory is generated analytically through a
hypergeometric representation and is used both as the initial guess
and as the basis for the error metrics reported below.

\paragraph{Pontryagin system.}
With costate $p=(p_1,p_2)$, the Hamiltonian is
$H(p,x,s,\alpha)=p_1\,\alpha/\sigma(s)+p_2+\bigl((\alpha-x)/\sigma(s)\bigr)^2$.
Stationarity
$\partial_\alpha H=p_1/\sigma(s)+2(\alpha-x)/\sigma(s)^2=0$ gives the
unconstrained minimizer
$\alpha^\star(t)=x(t)-p_1(t)\,\sigma(s(t))/2$, and after substitution
the reduced Hamiltonian becomes
\begin{equation}\label{eq:singular-H}
H(p,x,s)\;=\;\frac{p_1\,x}{\sigma(s)}-\frac{p_1^2}{4}+p_2,
\end{equation}
which is smooth. The example is singular in the sense that the
optimal control inherits a delicate localization through the weight
$\sigma(s)$:\ as $s\to t_0$ from either side, $\sigma(s)$ contracts towards
$\varepsilon^\beta$, and the dynamics, the running cost, and the
optimal control all develop strongly localized features near
$s=t_0$. Adaptive mesh refinement is therefore essential.

\paragraph{Run summary.}
With outer tolerances
$\varepsilon_{\mathrm{time}}=\varepsilon_{\mathrm{PA}}=
\varepsilon_\delta=10^{-4}$, $\mathrm{max\_iters}=18$,
explicit-gradient mode active throughout (since the analytic form
of~\eqref{eq:singular-H} is available), and the initial state and
costate guesses taken from the exact reference trajectory, the solver
reaches the \texttt{final\_resolve} action at the iteration budget
with $N=2122$ intervals (and no plane bundle, $M=0$). Main
numbers are collected in Table~\ref{tab:ex4-summary}.

\begin{table}[htbp]
\centering
\caption{Final principal results reported by the proposed solver on
 \RTadd{R62-B2}{Example~3} (singular tracking).}
\label{tab:ex4-summary}
\begin{tabular}{l l}
\toprule
Quantity & Value \\
\midrule
Mode & exact-smooth \\
Objective $J_h$ & $4.42\!\times\!10^{-14}$ \\
Exact objective & $0$ \\
Terminal state error & $-2.10\!\times\!10^{-7}$ \\
Max state error & $7.77\!\times\!10^{-2}$ \\
Max control error & $7.77\!\times\!10^{-2}$ \\
Max costate-$1$ error & $4.85\!\times\!10^{-3}$ \\
Mesh intervals $N$ & $2122$ \\
Planes $M$ & $0$ (exact-smooth mode) \\
$\eta_{\mathrm{time}}$ & $3.90\!\times\!10^{-6}$ \\
$\mathrm{tol}_{\mathrm{time}}^\star$ & $4.71\!\times\!10^{-8}$ \\
$\eta_{\mathrm{time,sum}}=\sum_n\bar r_n$ & $8.76\!\times\!10^{-6}$ \\
Final action & \texttt{final\_resolve} \\
Status & budget-limited \\
\bottomrule
\end{tabular}
\end{table}

\begin{table}[htbp]\centering
\caption{\RTadd{R74-K10}{Branch-firing diagnostic for Example~3
 (singular tracking). Exact-smooth mode with $M=0$ (no PA bundle)
 and $\delta = 10^{-2}$ held fixed throughout, so the adaptive
 vocabulary collapses to time-mesh refinement plus the terminal
 bookkeeping action;\ the run exhausts the iteration budget before
 the per-cell time tolerance is reached. The histogram is over the
 $19$-row log.}}
\label{tab:ex3-branches}
\RTadd{R74-K10}{\begin{tabular}{l c}
\toprule
Outer-iteration outcome & Count \\
\midrule
Time-mesh refinement alone     & $18$ \\
PA enrichment alone            & $0$ (structurally inactive) \\
$\delta$-halving alone         & $0$ (structurally inactive) \\
Combined operations            & $0$ \\
Terminal \texttt{final\_resolve} & $1$ \\
\bottomrule
\end{tabular}}
\end{table}

\begin{figure}[htbp]
\centering
\begin{minipage}{0.48\textwidth}\centering
 \includegraphics[width=\textwidth]{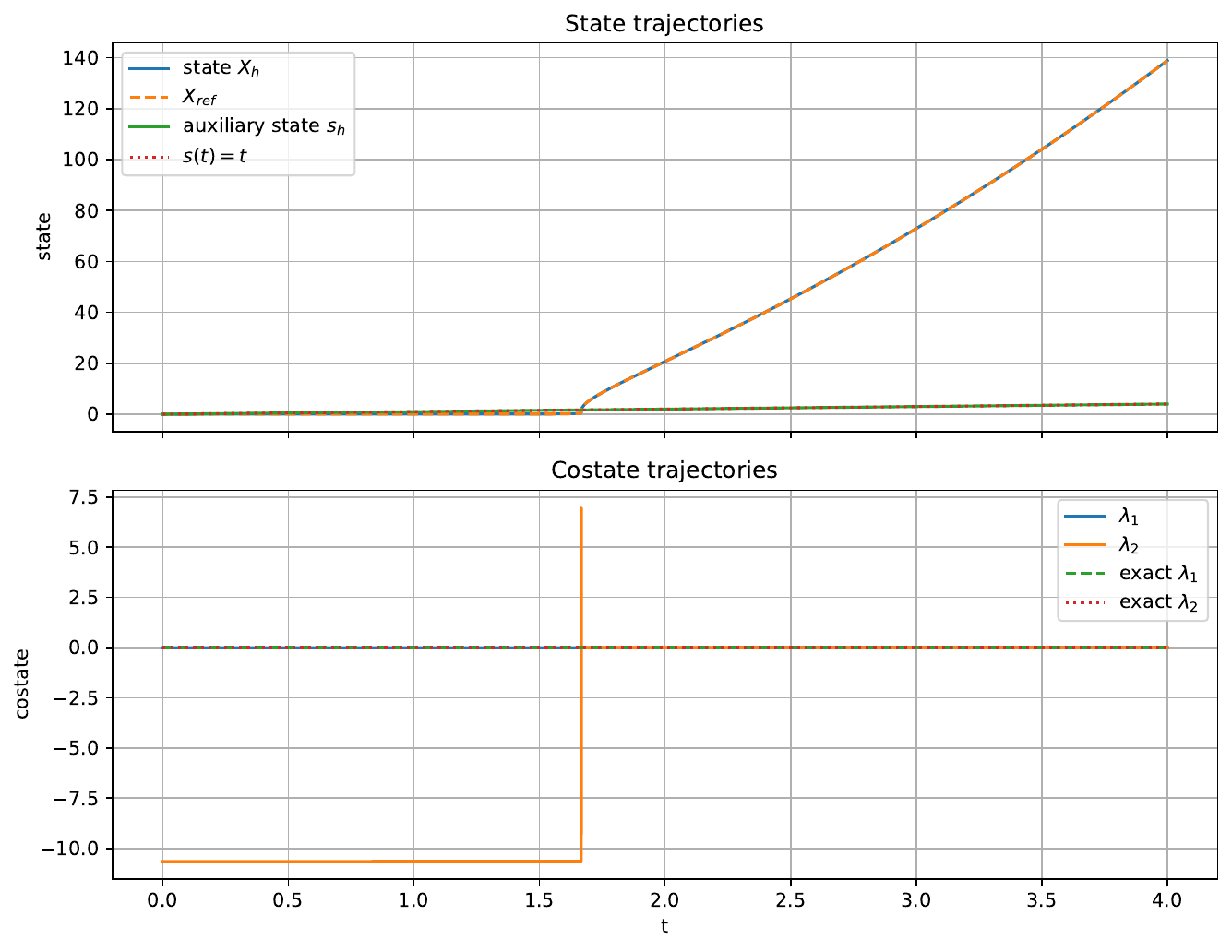}
\end{minipage}\hfill
\begin{minipage}{0.48\textwidth}\centering
 \includegraphics[width=\textwidth]{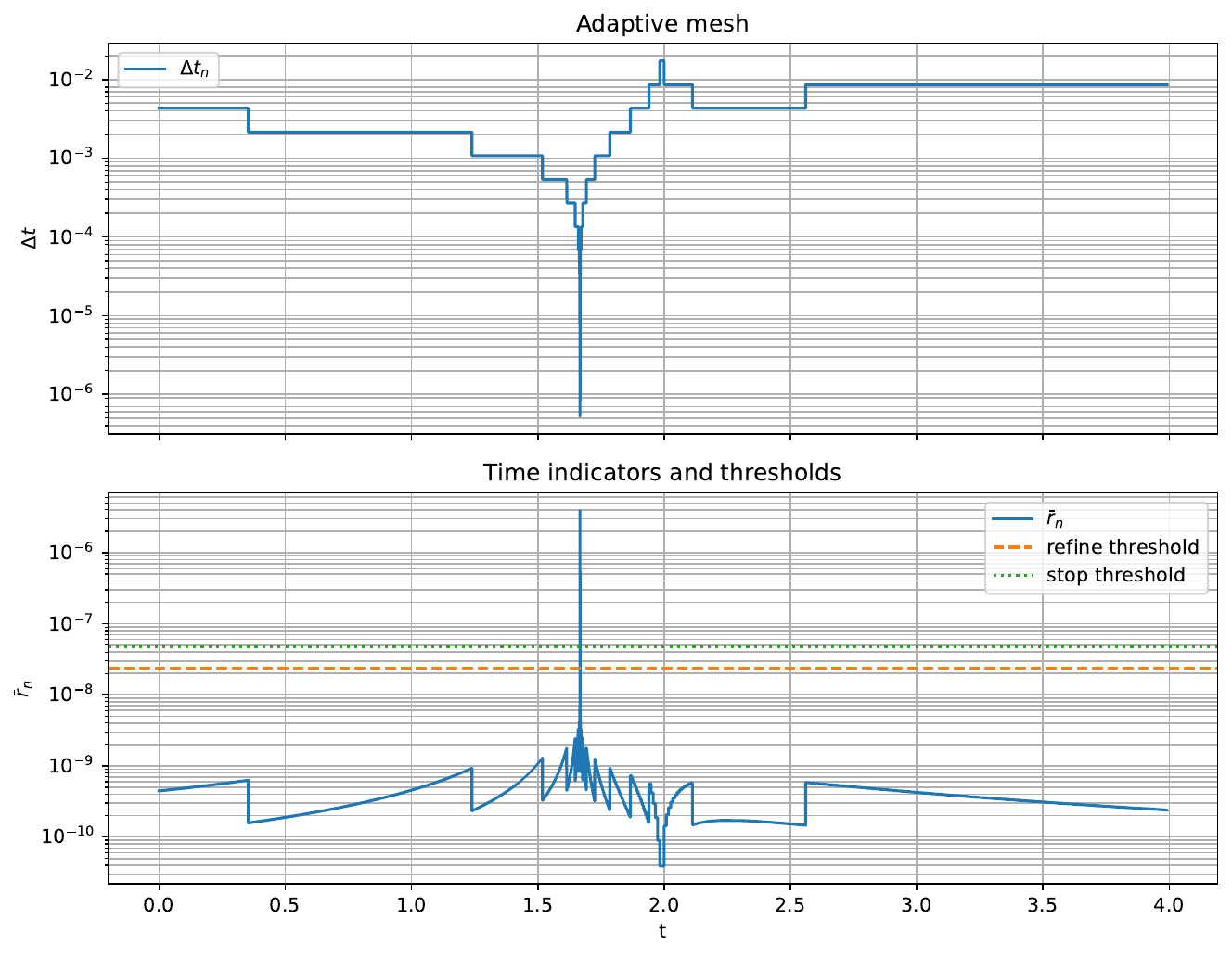}
\end{minipage}
\caption{\RTadd{R62-B2}{Example~3} (singular tracking). Left:\ state and first
 costate component at the final iterate, with the singular point
 $s=t_0=\tfrac53$ marked. Right:\ adaptive mesh and per-interval
 time indicator at the final iterate;\ the step distribution clusters
 strongly around the singular point and is comparatively coarse far
 from it.}
\label{fig:ex4-mesh}
\end{figure}

\paragraph{Discussion.}
The primal approximation is excellent:\ the objective is approximately
machine zero, the terminal state error is at the level of
$2\!\times\!10^{-7}$, and the first costate component is correct to
five digits. The maximum state error and the maximum control error
remain at the level of $8\!\times\!10^{-2}$ in a thin neighborhood of
the singularity at $s=t_0$, where the weight $\sigma(s)$ contracts to its
floor of $\varepsilon^\beta$ and any small mismatch in $\alpha$ is
amplified by~$1/\sigma(s)$. The adaptive mesh diagnostic in
Figure~\ref{fig:ex4-mesh} (right) shows that the time-step
distribution clusters strongly around $s=t_0$, exactly where the
local error density~$\rho_n$ peaks; mesh intervals more than a unit
of~$s$ from the singularity are kept comparatively coarse. Finally,
the second costate component~$p_2$ has limited diagnostic value in
this autonomous reformulation:\ it is the conjugate variable of the
auxiliary state $s$, and its discrete representation is sensitive to
the lifting and to round-off in the reference trajectory. We
therefore report it for completeness in the deep diagnostics but do
not include it in the main error metrics above. \RTadd{R66-S3}{Example~3}
demonstrates that, with explicit-gradient access, the proposed solver
reaches a primal-accurate solution on a strongly localized singular
benchmark with controlled mesh adaptation and without bundle
construction.

% ==================================================================
\RTadd{R53-B3}{%
\subsection{Example 4: QP-oracle attitude control allocation}
\label{ssec:ex-B}

\begin{figure}[htbp]
\centering
\IfFileExists{qp_oracle_setup.jpg}{%
  \includegraphics[width=0.75\textwidth]{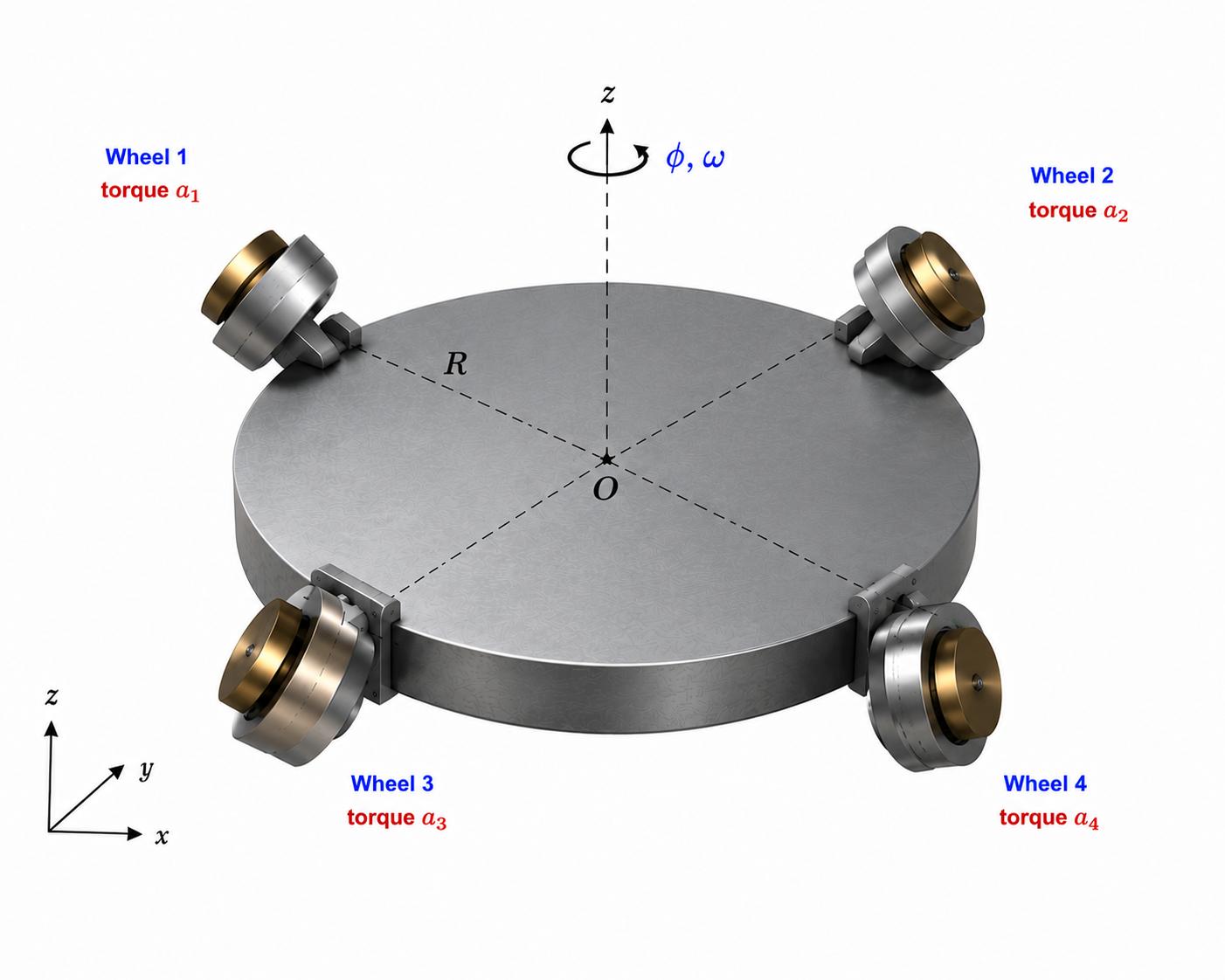}%
}{%
  \fbox{\parbox[c][6cm][c]{0.75\textwidth}{\centering
    \textit{Figure pending: QP-oracle attitude allocation setup
    (planar rigid body with four canted reaction wheels at radius
    $R$; drop the rendering at
    \texttt{qp\_oracle\_setup.jpg} in this directory to enable
    inclusion).}}}%
}
\caption{Schematic of the QP-oracle attitude control allocation example: a planar single-axis rigid body of inertia $J$ rotating about the $z$-axis with attitude $\phi$ and rate $\omega$, driven by four reaction-wheel torque commands $a_1, a_2, a_3, a_4$. Dynamics: $\dot\phi = \omega$, $\dot\omega = J^{-1}\,b^\top a$, with torque-to-axis coupling vector $b\in\mathbb{R}^4$ encoding the canted wheel geometry. Admissible set: per-wheel saturation $\|a\|_\infty\le a_{\max}$, total Euclidean-norm bound $\|a\|_2\le p_{\max}$, and a null-space-agility constraint $\|N a\|_\infty \le \nu$ forbidding spinning the wheels without on-axis torque; canonical constants $a_{\max} = 1.0$, $p_{\max} = 1.5$, $\nu = 0.3$. Schematic; wheel spin axes are canted with non-trivial $z$-component $b_i$ to produce attitude torque about the $z$-axis.}
\label{fig:qp-oracle-setup}
\end{figure}

\paragraph{Problem statement.}
We consider a planar single-axis rigid-body tracking problem with four
canted reaction-wheel torque commands $a(t)\in\mathbb{R}^4$
\Rzzedit{B2-S35}{(introduced and analysed in detail in the
costate-optimization foundation
companion~\cite{KsousLalvayCompanion2026}, which draws on the classical
reaction-wheel attitude-control formulation of~\cite{Wie2008}, the
control-allocation literature~\cite{HarkegardGlad2005,Bodson2002}, and
the parametric-QP solver theory
of~\cite{BemporadMorariDuaPistikopoulos2002,TondelJohansenBemporad2003})}
\Rzzedit{Sandberg-F8}{(see Figure~\ref{fig:exB-state-costate} for the
converged state and costate at the terminal iterate)}. The state
is \RTadd{R69-C6}{$x = (\phi,\omega)\in\mathbb{R}^2$} (attitude and rate); the dynamics
are $\dot\phi = \omega$, $\dot\omega = J^{-1}\,b^\top a$, where the
torque-to-axis coupling $b\in\mathbb{R}^4$ encodes the canted wheel
geometry. The cost penalizes attitude tracking error and control
effort:
\[
J(a) \;=\; \tfrac12\!\int_0^T\!\bigl(q\,(\phi-\phi_{\mathrm{ref}}(t))^2
+ r\,\|a\|_2^2\bigr)\,dt + \text{terminal penalty},
\]
with horizon $T=1$, $J=1$, $q=10$, $r=10^{-2}$. The admissible set is
\[
\mathcal{A} \;=\; \bigl\{ a\in\mathbb{R}^4 :\;
  \|a\|_\infty \le a_{\max},\;
  \|a\|_2 \le p_{\max},\;
  \|N a\|_\infty \le \nu \bigr\},
\]
with canonical constants $a_{\max} = 1.0$, $p_{\max} = 1.5$, and
$\nu = 0.3$, and $N \in \mathbb{R}^{3 \times 4}$ an orthonormal basis
for the null space of $b^\top$; the three constraint families encode,
in order, per-wheel torque saturation, a total Euclidean-norm bound
on the control vector, and a null-space-agility constraint forbidding
spinning the wheels without producing on-axis torque. The
intersection of the $\ell_\infty$ box, the $\ell_2$ ball of radius
$p_{\max}$, and the $\ell_\infty$ null-space slab is a non-trivial
QP. The Hamiltonian minimization
\[
H(p,\RTadd{R69-C6}{x},t) \;=\; \min_{a\in\mathcal{A}}\Bigl\{
p_\phi\,\omega + p_\omega\,J^{-1}b^\top a
+ \tfrac12 q\,(\phi-\phi_{\mathrm{ref}}(t))^2
+ \tfrac12 r\,\|a\|_2^2\Bigr\}
\]
is evaluated by a parametric-QP oracle rather than by a closed-form
formula. This is the only benchmark in the suite where the oracle
assumptions are exercised: no analytic $\arg\min$ is available.

\paragraph{Reference value (self-reference convention).}
\RTadd{R70-m2}{Because no closed-form reference exists, we use a
numerical self-reference convention:\ for the chosen tolerance level
$\theta_1$ we report the converged value $J_h$ and compare it to a
tighter self-reference $J_{\mathrm{ref}}$ produced by the same solver
at a stricter tolerance, in the spirit of the adaptive-benchmark
discussions of~\cite{KarlssonLarssonSandbergTempone2015}.} We use
$J_{\mathrm{ref}}=0.241166551748$ as the self-reference target for
the $\theta_1$ row.\
\RTadd{R62-M5}{This is not an analytic ground truth;\ the per-run
parameters (final mesh, bundle size, smoothing level, tolerance
triple, branch-firing counts, and terminal indicator values) are
reported in the principal results below and in
Table~\ref{tab:ex-B-branches} so that the comparison is fully
characterised from the manuscript itself.}

\paragraph{Run summary ($\theta_1$ tolerance).}
Under the overload-ratio outer-loop controller of
Algorithm~\ref{alg:outer}, the $\theta_1$ tolerance level
terminates with status \texttt{STOP} after $17$ outer iterations on
a final mesh of $N=328$ intervals with $M=15$ PA planes and final
smoothing $\delta=1.526\!\times\!10^{-5}$. The computed objective
is $J_h=0.24136739775749513$ against the tighter self-reference
$J_{\mathrm{ref}}=0.241166551748$, giving an absolute gap
$|J_h-J_{\mathrm{ref}}|=2.008\!\times\!10^{-4}$. The branch-firing
diagnostic for this run records $5$ time-mesh refinements,
$5$ PA-plane insertions, and $15$ $\delta$-halvings, with $9$
single-branch iterations, $5$ two-branch parallel iterations, and
$2$ three-branch parallel iterations (distributed as listed in
Table~\ref{tab:ex-B-branches}). At termination all three
indicators are below tolerance:
$\eta_{\mathrm{time}}^{\max}=3.91\!\times\!10^{-7}$,
$\eta_{\mathrm{PA}}=4.06\!\times\!10^{-6}$,
$\eta_\delta=5.22\!\times\!10^{-6}$. The empirical lower-comparability
constant settles at $c\approx 0.0996$, giving raw step-count speedup
$\mathcal{R}_{\mathrm{raw}}\approx 1.48$ and conservative speedup
$\mathcal{R}_c\approx 0.148$ -- consistent with the manuscript-side
convention of~\S\ref{ssec:quasi-norm-advantage}.

\begin{table}[h]\centering
\caption{Branch-firing diagnostic for the $\theta_1$ run.
Each row counts the number of outer iterations in which the named
combination of adaptive operations fired; the columns separate
single-branch, two-branch, and three-branch (parallel) outer
iterations. The terminal row collapses the STOP action.}
\label{tab:ex-B-branches}
\begin{tabular}{l c}
\toprule
Outer-iteration outcome & Count \\
\midrule
$\delta$-halving alone & $9$ \\
PA enrichment $+$ $\delta$-halving & $2$ \\
time refine $+$ $\delta$-halving & $2$ \\
time refine $+$ PA enrichment & $1$ \\
time refine $+$ PA enrichment $+$ $\delta$-halving & $2$ \\
STOP & $1$ \\
\bottomrule
\end{tabular}
\end{table}

\paragraph{Discussion.}
Example~4 is the most direct test of the oracle-based PA-softmin
architecture, because no closed-form Hamiltonian minimiser exists and
every outer iteration pays the cost of QP oracle calls. The
parallel-refinement structure of the controller is particularly
visible here:\ the $\theta_1$ run records $5+2=7$ outer iterations
in which two or three adaptive branches fired together, including
two three-branch (time$+$PA$+\delta$) iterations, so the bulk of
the work is concentrated in parallel updates rather than in a long
sequence of single-branch steps.
}

\RTadd{R57}{Figures~\ref{fig:exB-state-costate}--\ref{fig:exB-dt-evolution}
visualise four diagnostic layers of the converged $\theta_1$ run:\
the state and costate trajectories at the terminal iterate, the PA
bundle anchor times $\hat t_K$ in $(t,\phi)$ and $(t,\omega)$, the
time-error-density~$\bar\rho(t)$ evolution across outer iterations,
and the mesh~$\Delta t_n(t)$ evolution. The constraint-saturation
transition times $t_{c,1:3}$ are overlaid as dotted vertical lines
wherever they appear, marking the dynamically-distinguished instants
of the benchmark.

\begin{figure}[t]
 \centering
 \includegraphics[width=0.85\linewidth]{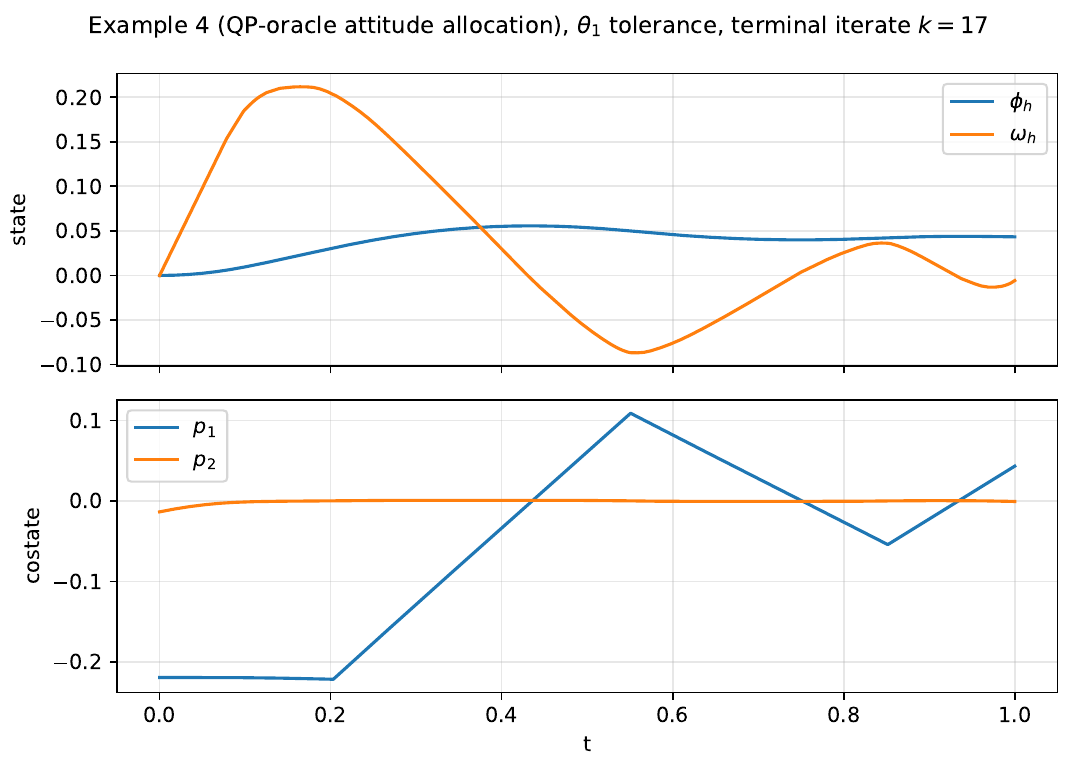}
 \caption{Example~4, $\theta_1$ tolerance, terminal iterate $k=17$:\
 state $(\phi_h,\omega_h)$ (top panel) and costate $(p_1,p_2)$
 (bottom panel) along the converged trajectory. $\phi$ tracks
 smoothly to its reference; the costate is piecewise linear on
 the final mesh, with switching structure reflecting the
 constraint-saturation transitions in the QP oracle's optimal
 control.}
 \label{fig:exB-state-costate}
\end{figure}

\begin{figure}[t]
 \centering
 \includegraphics[width=0.95\linewidth]{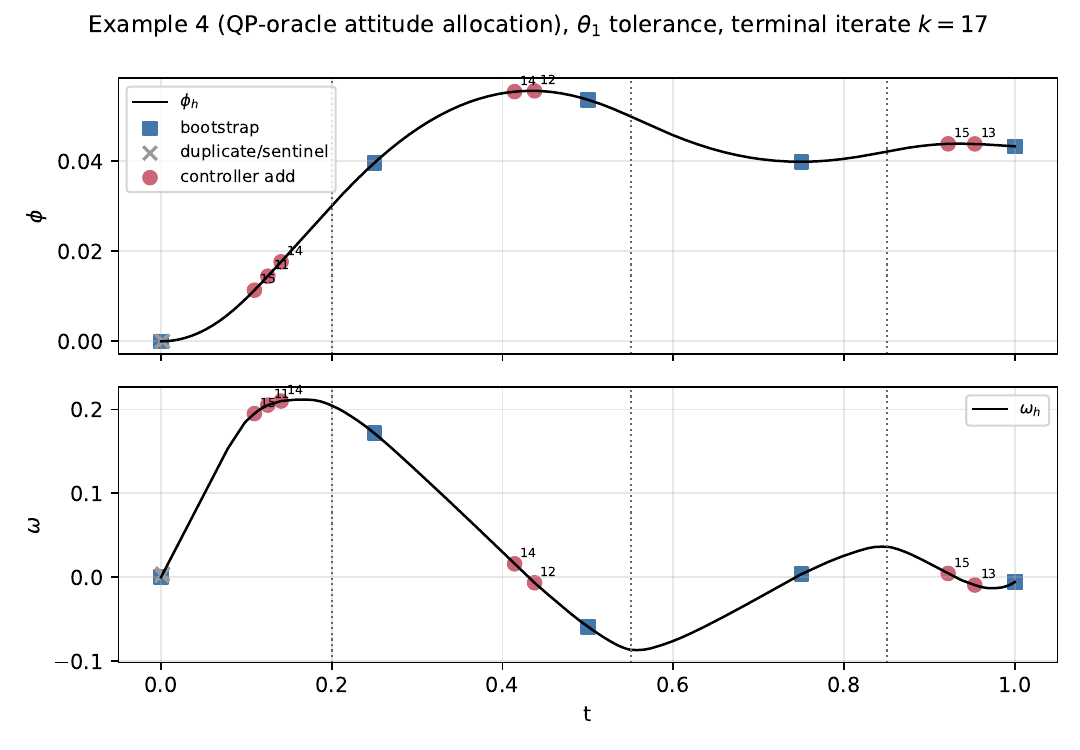}
 \caption{Example~4, $\theta_1$ tolerance, terminal iterate $k=17$:\
 PA bundle anchor knots $\hat t_K$ overlaid on the converged
 state trajectory in $(t,\phi)$ (top) and $(t,\omega)$ (bottom).
 Blue squares mark bootstrap-seed anchors;\ red circles mark
 anchors added by the controller via the
 \texttt{add\_plane} branch, with the corresponding outer-iteration
 index $k$ shown next to each. Crosses (if any) mark
 duplicate/sentinel controls retained in the bundle for
 provenance. The dotted vertical lines mark the
 constraint-saturation transition times $t_{c,1:3}$;\ the
 controller concentrates new anchors precisely at and just
 before these transitions, which is where the QP oracle's
 optimal control changes its active set and the PA bundle
 needs the most planes to retain a tight surrogate.}
 \label{fig:exB-bundle-anchors}
\end{figure}

\begin{figure}[t]
 \centering
 \includegraphics[width=0.85\linewidth]{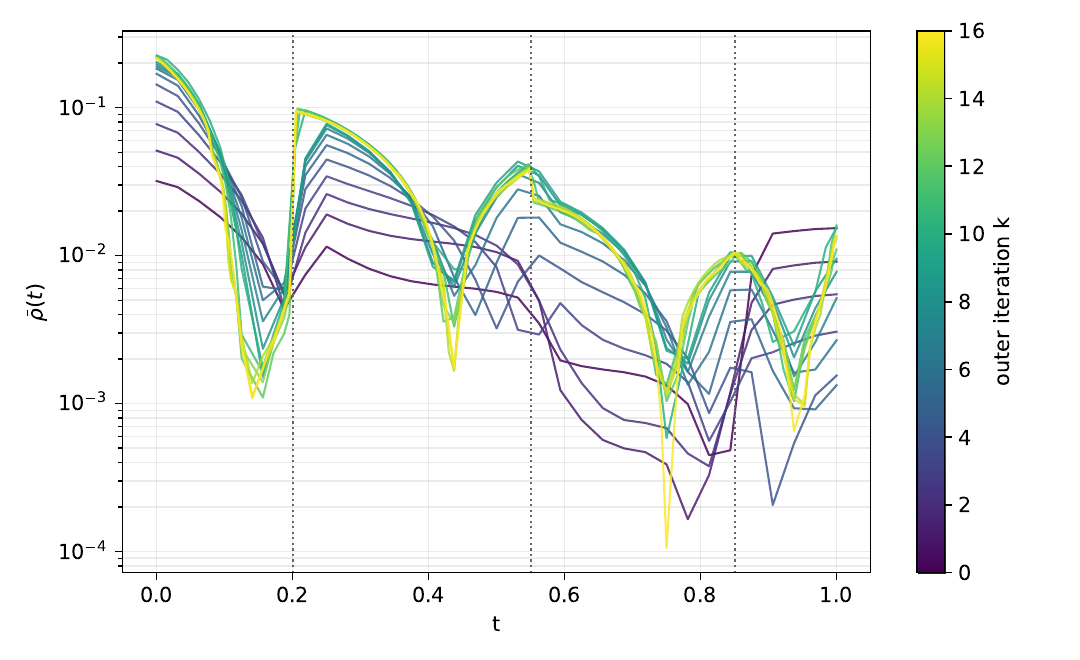}
 \caption{Example~4, $\theta_1$ tolerance:\ time-error-density
 $\bar\rho(t)$ at the mesh nodes, traced through outer iterations
 \RTadd{R62-E4}{$k=0,\ldots,16$} (colour gradient dark$\to$light). At early
 iterations $\bar\rho(t)$ is broad and smooth;\ at late iterations
 it sharpens into three peaks located precisely at the
 constraint-saturation transitions $t_{c,1:3}$ (dotted vertical
 lines). This is the direct evidence that the adaptive mesh
 refinement targets the dynamically distinguished times of the
 benchmark, and explains why mesh refinement is the dominant
 work source in the late phase of the outer loop.}
 \label{fig:exB-rho-evolution}
\end{figure}

\begin{figure}[t]
 \centering
 \includegraphics[width=0.85\linewidth]{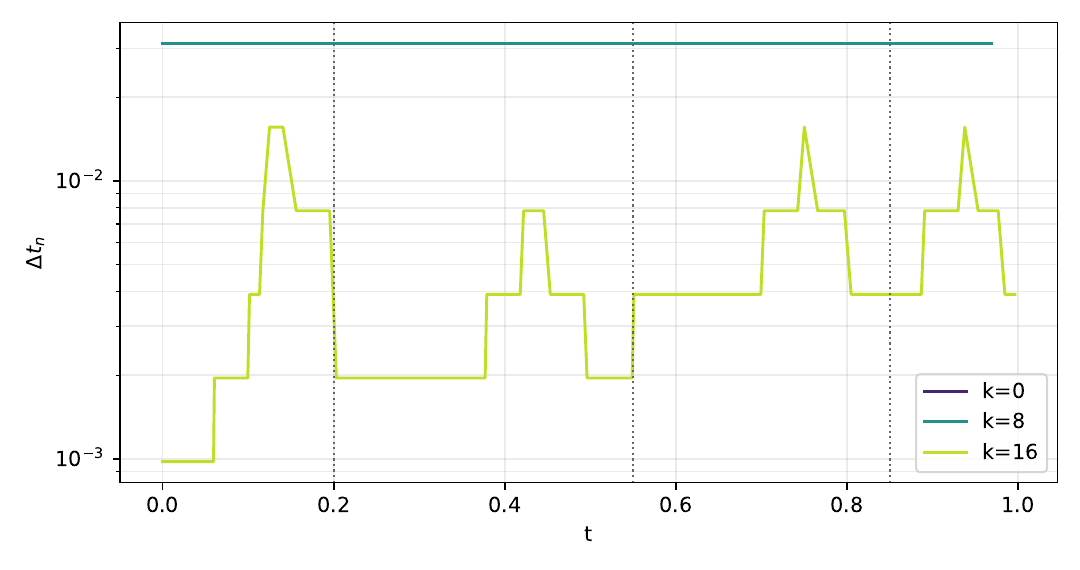}
 \caption{Example~4, $\theta_1$ tolerance:\ time-mesh spacing
 $\Delta t_n$ versus $t$ at three representative outer iterations
 ($k=0$ initial uniform mesh, $k=8$ mid-phase, $k=16$ near the
 terminal STOP). The initial and mid-phase meshes are
 indistinguishable on this log scale (both uniform); the late mesh
 has subdivided strongly around each of the three
 constraint-saturation transitions $t_{c,1:3}$ (dotted vertical
 lines), reaching $\Delta t_n\approx 10^{-3}$ locally while
 remaining $\sim 4\times 10^{-3}$ elsewhere. Together with
 Figure~\ref{fig:exB-rho-evolution} this gives a complete picture
 of the adaptive feedback loop:\ $\bar\rho(t)$ identifies the
 error hot-spots and $\Delta t_n(t)$ adapts to resolve them.}
 \label{fig:exB-dt-evolution}
\end{figure}
}

% ==================================================================
\subsection{Adaptive vs.\ uniform efficiency:\ the $L^{1/2}$ quasi-norm
 proxy}\label{ssec:quasi-norm-advantage}

\begin{figure}[t]\centering
\includegraphics[width=0.95\textwidth]{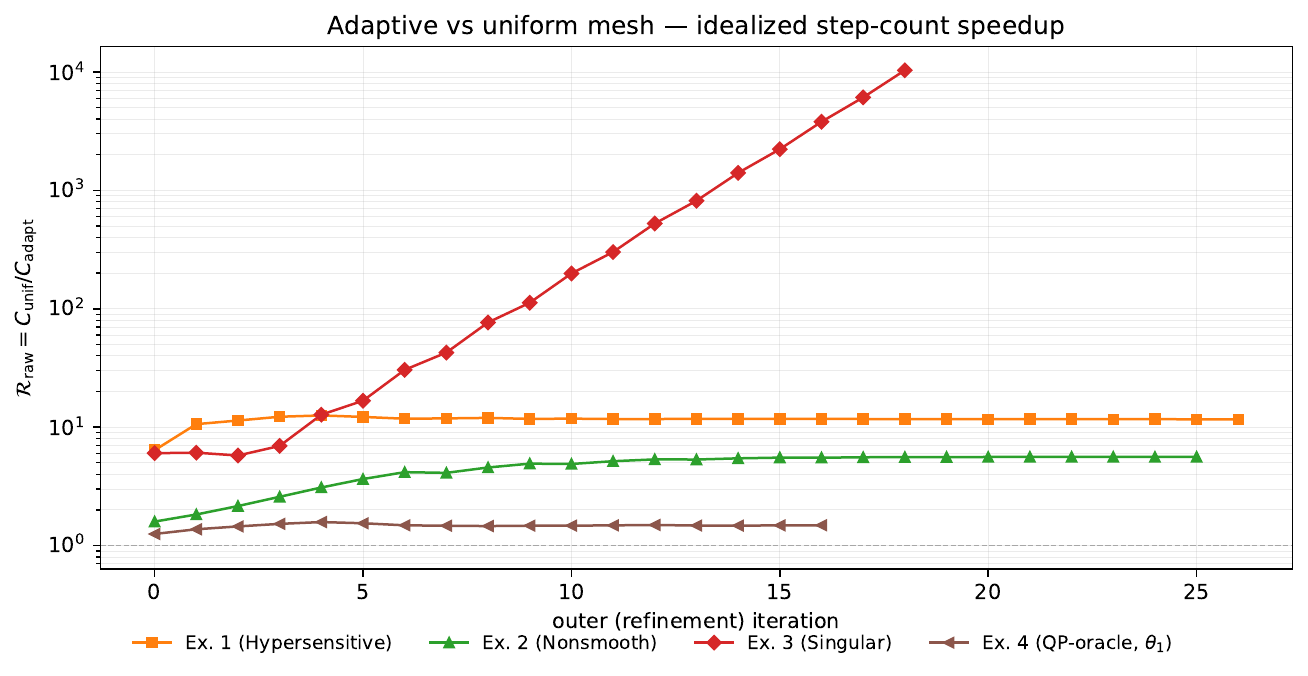}
\caption{\RTadd{R53-M5}{Asymptotic step-count speedup of the
 equidistributing adaptive mesh over a uniform mesh across outer
 (refinement) iterations for the four benchmarks of~\S\ref{sec:numerics}.
 The dashed reference line marks the break-even threshold $y=1$:\
 values above the line indicate that the adaptive mesh asks for fewer
 cells than the uniform mesh at a given tolerance.\ The figure shows
 the idealized speedup
 $\mathcal{R}_{\mathrm{raw}}=C_{\mathrm{unif}}/C_{\mathrm{adapt}}$
 of~\eqref{eq:Rraw};\ the conservative $c$-corrected
 $\mathcal{R}_c=c\,\mathcal{R}_{\mathrm{raw}}$ and the per-example
 empirical~$c$ are reported in the prose of each example's
 subsection.}}
\label{fig:quasi-norm-advantage}
\end{figure}

A complementary, equidistribution-based read on the four main
benchmarks (Examples~1--4 of~\S\ref{ssec:ex-hypersensitive}--\ref{ssec:ex-B})
is offered by the equidistribution coefficients of
\cite{KarlssonLarssonSandbergTempone2015}.\ The natural,
$\bar\rho$-derived, dimensionally clean quantity to compare across
runs is the \emph{asymptotic step-count coefficient}:\ the
proportionality constant $C(\bar\rho)$ such that, for symplectic
Euler with leading order \RTadd{R69-C3}{$q_{\mathrm{int}}=1$ (we use
$q_{\mathrm{int}}$ for the integrator order to avoid clash with the
PMP costate $p$ and with the tracking weight $q$ of the QP-oracle
cost in~\S\ref{ssec:ex-B})}, the number of cells required to
bring the total time-discretization error below~$\mathrm{TOL}$
satisfies
\begin{equation}
 N \;\approx\; \frac{C(\bar\rho)}{\mathrm{TOL}}.
 \label{eq:Nstep}
\end{equation}
For the two extreme strategies the coefficient is explicit:
\begin{equation}
 C_{\mathrm{unif}}(\bar\rho)
 \;:=\; T\,\|\bar\rho\|_{L^1}
 ,\qquad
 C_{\mathrm{adapt}}(\bar\rho)
 \;:=\; \|\bar\rho\|_{L^{1/2}}
 \;=\;\Bigl(\!\int_0^T |\bar\rho(t)|^{1/2}\,dt\Bigr)^{2},
 \label{eq:Ccoefs}
\end{equation}
with $C_{\mathrm{unif}}$ obtained by summing the per-cell error
$\bar\rho_n h^2$ on a uniform mesh of width $h=T/N$ and
$C_{\mathrm{adapt}}$ obtained by optimally equidistributing the
per-cell error.\ Both coefficients have the units of~$\mathrm{TOL}$,
so they are directly comparable row-by-row, and Jensen on the
concave map $x\mapsto x^{1/2}$ enforces
\begin{equation}
 C_{\mathrm{adapt}}(\bar\rho) \;\le\; C_{\mathrm{unif}}(\bar\rho),
 \label{eq:Cjensen}
\end{equation}
i.e.\ the adaptive mesh asymptotically needs at most as many cells
as the uniform mesh at every tolerance.\ The dimensionless
\emph{adaptive-vs-uniform step-count ratio}
\begin{equation}
 \mathcal{R}_{\mathrm{raw}}(k)
 \;:=\;
 \frac{C_{\mathrm{unif}}(\bar\rho^{(k)})}{C_{\mathrm{adapt}}(\bar\rho^{(k)})}
 \;=\;
 \frac{T\,\|\bar\rho^{(k)}\|_{L^1}}{\|\bar\rho^{(k)}\|_{L^{1/2}}}
 \;\ge\; 1
 \label{eq:Rraw}
\end{equation}
then reports the asymptotic step-count speedup of the adaptive over
the uniform mesh, as published in
\cite[Thm.~2.9]{KarlssonLarssonSandbergTempone2015}.\
For a theorem-faithful (non-asymptotic) read, the non-asymptotic upper bound of
\cite[Thm.~2.9, eq.~(2.22)]{KarlssonLarssonSandbergTempone2015}
replaces $C_{\mathrm{adapt}}$ by its conservative version
\begin{equation}
 C^{c}_{\mathrm{adapt}}(\bar\rho,c)
 \;:=\; \|\bar\rho/c\|_{L^{1/2}}
 \;=\; c^{-1}\,C_{\mathrm{adapt}}(\bar\rho),
 \label{eq:Ccadapt}
\end{equation}
where the empirical lower-comparability constant
$c^{(k)}\in(0,1]$ measured on the same trace
(\S\ref{ssec:adaptive-mesh}) controls the worst-case spread of the
per-interval error density across one subdivision;\ smaller~$c$
means a coarser local match between parent-and-children densities,
and therefore a more conservative step-count estimate.\ The
corresponding ratio
\begin{equation}
 \mathcal{R}_c(k)
 \;:=\;
 \frac{C_{\mathrm{unif}}(\bar\rho^{(k)})}{C^{c}_{\mathrm{adapt}}(\bar\rho^{(k)},c^{(k)})}
 \;=\; c^{(k)}\,\mathcal{R}_{\mathrm{raw}}(k)
 \label{eq:Rc}
\end{equation}
is the conservative speedup that survives without invoking the
asymptotic limit.\
Figure~\ref{fig:quasi-norm-advantage} traces $\mathcal{R}_{\mathrm{raw}}$
across outer iterations for the four benchmarks
of~\S\ref{sec:numerics} on a logarithmic scale, so the singular
benchmark (whose gain spans four orders of magnitude) does not crowd
the smoother problems out of view.\
Table~\ref{tab:quasi-norm-advantage} reports the terminal-iterate
step-count coefficients themselves together with the two speedup ratios.
\RTadd{R62-M6}{The ratios reported in this subsection concern only the
asymptotic number of time cells predicted by the time-error density;\
they do not include PA oracle calls, bundle enrichment, Newton
iterations, sparse-LU factorisation cost, or QP-solve cost.}

\begin{table}[t]\centering
\caption{\RTadd{R53-M5}{Terminal-iterate asymptotic step-count coefficients
 $C(\bar\rho)$ of~\eqref{eq:Nstep}--\eqref{eq:Ccoefs} and the idealized
 step-count speedup $\mathcal{R}_{\mathrm{raw}}$ (symplectic Euler,
 \RTadd{R69-C3}{$q_{\mathrm{int}}=1$}, so $1/(\RTadd{R69-C3}{q_{\mathrm{int}}}+1)=1/2$).\ $T$ is the horizon of each benchmark.\
 The two coefficients $C_{\mathrm{unif}}$ and $C_{\mathrm{adapt}}$
 share the units of $\mathrm{TOL}$ and are Jensen-comparable row by
 row;\ their ratio gives the idealized step-count speedup of the
 equidistributing adaptive mesh over the uniform mesh.\ The empirical
 comparability constant~$c$ (MSTZ~2.22) and the corresponding
 conservative speedup $\mathcal{R}_c=c\,\mathcal{R}_{\mathrm{raw}}$
 are reported per-example in the prose
 of~\S\ref{ssec:ex-hypersensitive}--\ref{ssec:ex-B} rather than in
 this principal table, since the table's role is to capture the
 idealized speedup that Jensen on the concave map $x\mapsto x^{1/2}$
 enforces.}}
\label{tab:quasi-norm-advantage}
\RTadd{R53-M5}{\setlength{\tabcolsep}{6pt}
\begin{tabular}{l c c c c}
\toprule
Run & $T$ &
$C_{\mathrm{unif}}$ &
$C_{\mathrm{adapt}}$ &
$\mathcal{R}_{\mathrm{raw}}$ \\
\midrule
1 (Hypersensitive) & $25$ & $5.10\!\times\!10^{1}$ & $4.39$ & $11.6$ \\
2 (Nonsmooth) & $1$ & $4.92\!\times\!10^{-4}$ & $8.78\!\times\!10^{-5}$ & $5.60$ \\
3 (Singular) & $4$ & $4.97\!\times\!10^{1}$ & $4.81\!\times\!10^{-3}$ & $1.03\!\times\!10^{4}$ \\
4 (QP-oracle attitude alloc.) & $1$ & $3.25\!\times\!10^{-2}$ & $2.20\!\times\!10^{-2}$ & $1.48$ \\
\bottomrule
\end{tabular}}
\end{table}

\RTadd{R53-M5}{The terminal-iterate numbers in
Table~\ref{tab:quasi-norm-advantage} let one read, for each
benchmark, how many cells the asymptotic theory predicts each mesh
strategy needs to bring the time-discretization error below a
prescribed tolerance:\ via~\eqref{eq:Nstep}, the uniform mesh asks
for $N\approx C_{\mathrm{unif}}/\mathrm{TOL}$ cells and the
equidistributing adaptive mesh asks for
$N\approx C_{\mathrm{adapt}}/\mathrm{TOL}$, with the ratio of the
two being the idealized speedup $\mathcal{R}_{\mathrm{raw}}$.\ The
singular benchmark (Example~3) is the dominant entry:\ its uniform
coefficient $C_{\mathrm{unif}}\approx 4.97\!\times\!10^{1}$ and its
adaptive coefficient $C_{\mathrm{adapt}}\approx 4.81\!\times\!10^{-3}$
give a four-orders-of-magnitude reduction in the asymptotic step
count and a $\mathcal{R}_{\mathrm{raw}}\approx 1.03\!\times\!10^{4}$
speedup, consistent with the visual clustering of mesh points around
the singular point reported in Figure~\ref{fig:ex4-mesh}---adaptive
refinement is essential at any tolerance one cares about.\ The
hypersensitive (Example~1) and nonsmooth (Example~2) benchmarks
deliver a strong, robust asymptotic advantage of roughly one order of
magnitude ($\mathcal{R}_{\mathrm{raw}}\approx 11.6$ and $5.60$,
respectively).\
The QP-oracle attitude control allocation (Example~4) gives a more
modest idealized gain ($\mathcal{R}_{\mathrm{raw}}\approx 1.48$) and
offers a different diagnostic signal:\ its empirical $c$ stabilises
near~$0.10$ (the $c$-corrected speedup $\mathcal{R}_c\approx 0.148$
is reported in the prose of~\S\ref{ssec:ex-B}). \RTadd{R66-M7}{For
comparison, the empirical $c$ values are:\ LQR appendix consistency check
$c\approx 0.78$, hypersensitive $0.76$, nonsmooth $0.79$, singular
$0.11$, and QP-oracle attitude allocation $0.10$. These values are
read off the converged iterates and discussed in the corresponding
example subsections~\S\ref{ssec:ex-hypersensitive}--\ref{ssec:ex-B} and
in Appendix~\ref{apx:lqr}.}}

% ==================================================================
\subsection{Summary}\label{ssec:numerics-summary}

\RTadd{R53-B3}{Table~\ref{tab:summary} collects the principal results
across the four unconstrained-state examples. All four runs use the
proposed solver in its adaptive form. The explicit-gradient and
analytic-smoothing path of~\S\ref{ssec:explicit-grad} is active in
Example~2 (nonsmooth); Example~3 (singular tracking) runs in
exact-smooth mode (mode~(i)) with $M=0$ because the analytic reduced
Hamiltonian is available directly. The analytic-smoothing path is
applicable wherever an analytic Hamiltonian gradient is available.
Example~4 (QP-oracle attitude control allocation) is the only
benchmark that exercises the oracle assumptions.}

\begin{table}[t]\centering
\caption{Summary of the numerical results. ``$J^{\mathrm{ref}}$'' is
 the reference objective when available; the time, PA, and smoothing
 indicators $\eta_{\mathrm{time}}, \eta_{\mathrm{PA}}, \eta_\delta$
 are evaluated at the converged iterate, with \texttt{---} entries
 marking structurally inactive layers. The Status column reports the
 outer-loop exit branch and its acceptance verdict (see the legend
 below the table).}
\label{tab:summary}
\resizebox{\textwidth}{!}{%
\begin{tabular}{l c c c c c c c c c c}
\toprule
Ex.\ &
$J^{\mathrm{ref}}$ & $J_h$ & $|J_h-J^{\mathrm{ref}}|$ &
$N$ & $M$ & $\delta$ &
$\eta_{\mathrm{time}}$ &
$\eta_{\mathrm{PA}}$ &
$\eta_\delta$ &
\RTadd{R66-M8}{Status} \\
\midrule
1 (Hypersensitive) & (n/a) & $2.27475$ & --
 & $307$ & $30$ & $3.13\!\times\!10^{-4}$
 & $1.31\!\times\!10^{-4}$ & $5.09\!\times\!10^{-3}$ & $3.85\!\times\!10^{-3}$
 & \RTadd{R66-M8}{\texttt{STOP}/accepted} \\
2 (Nonsmooth) & $4.439\!\times\!10^{-5}$ & $4.471\!\times\!10^{-5}$ & $3.24\!\times\!10^{-7}$
 & $275$ & $3$ & $7.63\!\times\!10^{-8}$
 & $2.56\!\times\!10^{-9}$ & $0$ & $4.60\!\times\!10^{-8}$
 & \RTadd{R66-M8}{\texttt{STOP}/accepted} \\
3 (Singular) & $0$ & $4.42\!\times\!10^{-14}$ & $4.42\!\times\!10^{-14}$
 & $2122$ & $0$ & --
 & $3.90\!\times\!10^{-6}$ & -- & --
 & \RTadd{R66-M8}{budget-limited \texttt{final\_resolve}} \\
4 (QP-oracle, $\theta_1$) & $2.4117\!\times\!10^{-1}$ & $2.4137\!\times\!10^{-1}$ & $2.01\!\times\!10^{-4}$
 & $328$ & $15$ & $1.53\!\times\!10^{-5}$
 & $3.91\!\times\!10^{-7}$ & $4.06\!\times\!10^{-6}$ & $5.22\!\times\!10^{-6}$
 & \RTadd{R66-M8}{\texttt{STOP}/accepted} \\
\bottomrule
\end{tabular}}
\end{table}

\RTadd{R70-M5}{\noindent\emph{Status legend.}\;
\texttt{STOP}/accepted denotes the converged exit of
Algorithm~\ref{alg:outer} with all three indicator-to-tolerance
ratios at or below one.\ \texttt{final\_resolve}/post-budget denotes
a final post-budget refinement or resolve action that the
diagnostic gates accept, distinct from the converged \texttt{STOP}
exit;\ the ``post-budget'' suffix records that this exit branch is
taken after the iteration budget is exhausted, not by the converged
stop test.\ \RTadd{R72-M7}{In Example~3 (singular tracking) the
status is budget-limited:\ the run exits through the
\texttt{final\_resolve} branch because the iteration budget is
reached before the time indicator falls below the per-cell
threshold.}}

\RTadd{R55}{Example~4 (QP-oracle attitude control allocation) numbers
above are the converged $\theta_1$ run under the overload-ratio
controller of Algorithm~\ref{alg:outer};\ see~\S\ref{ssec:ex-B} for
the branch-firing diagnostic and the per-iteration
parallel-refinement structure.}

\RTadd{R74-K10}{\paragraph{Branch-firing breakdown across the four
benchmarks.}
The per-example branch-firing histograms are reported in
Tables~\ref{tab:ex1-branches},~\ref{tab:ex2-branches},
\ref{tab:ex3-branches}, and~\ref{tab:ex-B-branches}, one per
benchmark. Examples~1, 2, and~3 record sequential single-action
outer iterations:\ each non-terminal outer iteration fires exactly
one of the three adaptive branches (time-mesh refinement, PA
enrichment, or $\delta$-halving). Example~4 (QP-oracle) is the only
benchmark in the suite where the outer loop fires two or three
branches in parallel within a single outer iteration, reflecting the
overload-ratio controller of Algorithm~\ref{alg:outer} on the only
fully mode-(iii) row.}

\RTadd{R53-B3}{The four examples illustrate the asymmetric roles of
the three indicators predicted in~\S\ref{ssec:indicators} and
exercised by the outer loop of Algorithm~\ref{alg:outer}, and they
cover the three solver modes of \S\ref{ssec:explicit-grad}. In
Example~1 (hypersensitive), the solver runs in PA-softmin mode:\ the
bundle and smoothing indicators are active in the early outer
iterations and the time indicator becomes the binding constraint
once the mesh has been refined. The run reaches an accepted
\texttt{STOP} with all three indicators below their tolerances.
In Example~2 (nonsmooth) the analytic-smoothing mode reduces the
outer loop to a two-indicator scheme (time and~$\delta$):\ the PA
bundle is exact once the three planes $\{-1,0,+1\}$ are present, so
$\eta_{\mathrm{PA}}\equiv 0$ \RTadd{R53-M9}{\emph{structurally}, while
$\eta_\delta$ remains active and is set by the run configuration};
$\eta_\delta$ measures the residual analytic-smoothing bias. In
Example~3 (singular tracking) exact-smooth mode is active---the
analytic reduced Hamiltonian~\eqref{eq:singular-H} is plugged in
directly---and the time indicator drives the local clustering of mesh
points around the singular weight. \RTadd{R66-S4}{In Example~4
(QP-oracle attitude alloc.) all three branches engage, with several
outer iterations refining two or three layers simultaneously;\ PA
enrichment is strongest early, $\delta$-halving remains active over
many iterations, and time subdivision resolves the final transition
structure.}}

\paragraph{Deep diagnostics.}
The solver implementation has been instrumented to produce, for
every outer iteration of every example,
(a)~per-interval PA and smoothing contribution plots,
(b)~bundle-support-point overlays on the optimal-control plot,
(c)~history curves for $\eta_{\mathrm{time}}$,
$\eta_{\mathrm{time,sum}}$, $\eta_{\mathrm{PA}}$, and
$\eta_\delta$, and (d)~PA-enrichment tables listing selected node
indices, times, gaps, local~$\Delta t$, ranking scores,
accept/reject status, bundle size after insertion, and the number
of candidates rejected by the time-separation filter or the
duplicate-control filter. These per-iteration diagnostics inform
the convergence narratives of \S\S\ref{ssec:ex-hypersensitive}--%
\ref{ssec:ex-B}.
%% === end of Numexperiments_v2.tex ===

\section{Discussion and Outlook}\label{sec:discussion}

%% === inlined from Discussion_v2.tex (R98z-cc arxiv-flat) ===
% ------------------------------------------------------------------
% §5 Discussion and Outlook
% ------------------------------------------------------------------

\subsection{Relation to prior work}

\MSadd{\paragraph{Recombination, not replacement.}
The present solver should be read as a \emph{recombination} of classical optimal-control
algorithmic ingredients into a Hamiltonian-centered architecture, not as a replacement for any of
its three benchmark classes. Indirect shooting~\cite{Pontryagin1962,Bryson1975}, direct
transcription~\cite{Betts1998,Elnagar1995,Rao2014}, and Hamilton--Jacobi--Bellman PDE
solvers~\cite{FlemingSoner2006} are mature method families. The novelty here is the
\emph{combination}: a concave piecewise-affine bundle surrogate of the Hamiltonian, log-sum-exp
smoothing of that surrogate, symplectic-Euler discretization, full-space all-at-once damped
Newton on the canonical $(x,p)$ pair, and a unified adaptive controller balancing time-mesh,
PA-bundle, and smoothing layers. Each ingredient is classical; the assembly is not. Whether this
assembly demonstrably outperforms strong direct-collocation baselines on canonical benchmarks is
an empirical question that head-to-head comparison work --- pending future numerical rounds --- is
designed to settle. The mathematical content of this work is algorithmic rather than a
new convergence doctrine.}

\emph{Symplectic Pontryagin discretization.}\;
\SZadd{AR7}{The work}~\cite{SandbergSzepessy2006} studied a symplectic
Pontryagin scheme for deterministic optimal control and gave a
convergence analysis under suitable regularity assumptions on the
Hamiltonian.\SZadd{AM10}{\ Broader background on symplectic integration
of Hamiltonian systems, including the constrained case, is given by
\cite{LeimkuhlerReich2004}.}%
\MSnote{2026 audit (p.~8, minor~3.5): the earlier phrasing ``provided
convergence rates for regularized problems as $\delta,\Delta t\to 0$''
paraphrased more than the cited work literally contains; we now refer
generically to their convergence analysis and keep the use-made-of
this result at a conceptual level.}\;
Our solver builds on that discretization but replaces the assumed
smooth Hamiltonian by a soft-min-smoothed PA bundle model and adds
explicit adaptive refinement of all three error sources (time, model,
smoothing). The error ``knobs'' in our outer loop correspond
conceptually to the terms in their analysis, but we make them
computable and controllable.

\Rzzedit{bundle}{\emph{Bundle methods in non-smooth convex
optimization.}\;
The piecewise-affine bundle surrogate $\bar H$ of \S\ref{ssec:PA-surrogate}
shares its data structure with the \emph{bundle methods} developed for
non-smooth convex optimization since~\cite{Lemarechal1975};
\cite{Kiwiel1985} and \cite{HiriartUrrutyLemarechal1993} are standard
references. There, a collection of affine cuts accumulated from past
subgradient calls is used to build a cutting-plane / bundle model of a
non-smooth convex objective, and the smoothing approach
of~\cite{Nesterov2005} (which we cite in~\S\ref{ssec:smoothing} for the
log-sum-exp inequality) is a standard way to make such bundle models
differentiable. We note this parallel as a structural analogy rather
than a methodological inheritance:\ in bundle methods the bundle is the
optimization target and is paired with a proximal or trust-region step,
whereas in the present work the bundle is a Hamiltonian surrogate
parameterized in $(x,t)$ via the control-based extension
of~\S\ref{ssec:extension} and is paired with a Newton solve on the
discrete TPBVP residual of~\S\ref{ssec:newton}. We arrived at the same
piecewise-affine data structure from the Hamiltonian-approximation side,
not by specializing a bundle-method algorithm.}

\emph{Adaptive time stepping for PMP.}\;
\SZadd{AR4}{The work}~\cite{KarlssonLarssonSandbergTempone2015} introduced
a~posteriori error estimates for the symplectic Euler discretization
of PMP and demonstrated an adaptive time-stepping algorithm. We adopt
their error-density-driven step control (Section~\ref{ssec:indicators},
Section~\ref{ssec:adaptive-mesh}) and embed it into a larger adaptive
scheme that simultaneously refines the Hamiltonian approximation and
the smoothing level.

\emph{Multiple shooting and classical TPBVP solvers.}\;
On the numerical side, the full-space / multiple-shooting viewpoint
adopted in Section~\ref{ssec:newton} is related to the classical
multiple-shooting methods for optimal control pioneered by Bock and
Plitt~\cite{BockPlitt1984}, with whom the block-tridiagonal segment
structure and its linear-in-$S$ cost \Rzzedit{B2-S36}{(where $S$
denotes the number of multiple-shooting segments)} originate.\MSdel{ This work does not treat state-constrained optimal control; the corresponding state-constrained indirect-methods background and the smooth state-constrained-indirect approaches related to the present solver are discussed in the forthcoming costate-optimization companion work that extends the framework to explicit state constraints.}

\emph{Classical shooting vs.\ this method.}\;
Traditional shooting methods for optimal control solve the TPBVP by
Newton iteration, which requires differentiability of the shooting
function. When the Hamiltonian is nonsmooth (e.g., bang--bang
controls), the shooting Jacobian is undefined and Newton may fail. By
introducing the smoothed Hamiltonian $H_\delta$
(Section~\ref{ssec:smoothing}), we ensure differentiability while
preserving the concavity structure. The PA surrogate avoids expensive
exact Hamiltonian minimizations at each Newton step.

\emph{Direct methods.}\;
\MSdel{Compared to direct transcription~\cite{Betts1998,Elnagar1995,Rao2014},
which reformulates the problem as a finite-dimensional nonlinear
program (NLP) solved via sparse KKT systems, our method preserves the
Hamiltonian structure of the PMP and avoids the explosion in problem
size that direct transcription entails for fine meshes. On the other
hand, direct methods handle constraints more systematically via
general-purpose NLP solvers. The proposed method avoids the explicit
state-space grid of classical HJB discretizations but its cost still
depends on the Newton linear algebra, the surrogate size, and the
oracle model.}%
\MSadd{Direct transcription methods~\cite{Betts1998,Elnagar1995,Rao2014}
reformulate the problem as a finite-dimensional nonlinear program
(NLP) in the joint state--control vector and call a general-purpose
NLP solver with sparse KKT linear algebra. Relative to that family,
the present solver trades breadth for structure. On the \emph{plus}
side, the Newton system of Section~\ref{ssec:newton} works on the
$w=(x,p)$ canonical pair only (no control enters as an NLP unknown),
inherits the block-tridiagonal sparsity of the symplectic Euler
discretization, and is driven by a concave smoothed Hamiltonian that
keeps the search direction well defined even near bang--bang regimes.
On the \emph{minus} side, direct methods are more flexible: inequality
constraints and mixed state--control constraints are handled
uniformly by the NLP solver's active-set or interior-point machinery,
whereas the present Pontryagin-based solver is restricted to the
unconstrained-state regime\MSdel{; state-constrained extensions are the subject of a forthcoming costate-optimization companion work}. Mesh
size is not a fair differentiator either way: modern sparse NLP
solvers scale well in the number of collocation nodes, and the
operation count of our solver likewise grows linearly in~$N$
(Section~\ref{ssec:newton}). The two approaches should be
read as complementary, with the choice dictated by how tightly one
wants to preserve the Hamiltonian structure versus how heterogeneous
the constraint set is. Separately, we avoid the explicit state-space
grid of HJB discretizations, but our cost still scales with the
Newton linear algebra, the bundle size, and the oracle cost.}%
\MSnote{2026 audit (p.~7, major~6): soften the ``our method preserves
the Hamiltonian structure and avoids the NLP explosion'' phrasing
to a balanced comparison that names trade-offs in both directions.}%
\MSadd{\ The present implementation also uses sparse LU
factorization (\texttt{splu}) of the block-tridiagonal Newton system
and a nonlinear-least-squares fallback (\texttt{least\_squares}) when
damped Newton stalls (Section~\ref{ssec:newton});\ both are standard
ingredients in modern NLP solvers and close the
``implementation-maturity'' gap between indirect and direct
approaches.}%
\MSnote{2026-04-24 update: name the sparse-LU assembly and the
least-squares fallback as part of the Direct methods comparison so
the maturity gap with modern NLP solvers is explicit.}

\emph{HJB and max-plus methods.}\;
Our trajectory-optimization solver finds one optimal trajectory rather
than the full value function, which is why it scales polynomially in
time discretization and state dimension (for fixed oracle cost). The
downside is local optimality: PMP finds locally optimal solutions and
requires an initial guess, whereas HJB~\cite{FlemingSoner2006}
computes the global optimum at the cost of the curse of
dimensionality. Max-plus~\cite{Akian2008} and occupation-measure
relaxations~\cite{Lasserre2008} attempt to mitigate this, but face
their own scalability issues.

\subsection{Limitations}\label{ssec:limitations}

\emph{Local optimality.}\;
Like any shooting/Newton method, the solver finds local optima. If the
problem is nonconvex, convergence to a suboptimal trajectory depends
on the initial guess. The smoothing helps by eliminating trivial
nonsmoothness but does not convexify the problem. Running from
multiple initial guesses or employing homotopy is recommended.

\MSdel{Unconstrained-state scope limitation paragraph removed at R10 per user directive (the paragraph was a forward pointer to a state-constrained companion; replaced implicitly by the unconstrained-state scope already stated in the abstract and standing assumptions).}

\emph{Jacobian conditioning.}\;
The full-space all-at-once Newton Jacobian of~\S\ref{ssec:newton} can
become ill-conditioned at small $\delta$ (Remark~\ref{rem:H_delta_hessian})
and for long horizons or nearly singular dynamics. Careful scaling of
state and costate variables, continuation strategies, and---in the
multiple-shooting realisation of~\S\ref{ssec:newton}---shorter
shooting intervals can mitigate this.

\emph{State and path constraints.}\;
The present solver assumes unconstrained-state problems; extension to
state-constrained and path-constrained regimes is treated in a
state-constrained companion work and is outside the scope here.

\emph{No second-order verification.}\;
We check first-order PMP conditions only. Verifying second-order
sufficiency (conjugate point condition or definiteness of the shooting
Jacobian) is a possible extension.

\RTadd{R66-S5}{\emph{PA-majorant one-sided interpretation.}\;
The one-sided interpretation of $\eta_{\mathrm{PA}}$ relies on the
control-based majorant extension of~\S\ref{ssec:extension}.\ If the
method is transferred to black-box or learned Hamiltonians without
such an extension, $\eta_{\mathrm{PA}}$ becomes diagnostic rather than
a certified one-sided model-gap indicator;\ the design of adaptive
PMP shooting in this regime is left to future work (\S5.3).}

\subsection{Outlook}\label{ssec:outlook}

The encouraging results open several directions: using the computed
PMP solution to seed value-function approximations (e.g., training a
neural network near the optimal trajectory), extending the method to
controlled diffusions using the companion stochastic framework, and
embedding the solver within a hierarchical or model-predictive control
scheme. An ablation study quantifying the marginal contribution of
each adaptive component (time, PA, $\delta$) would further strengthen
the empirical case. Two further directions are particularly natural:
scaling the solver to higher state dimensions ($d \gg 5$), where the
Newton linear algebra and the PA-bundle sampling cost both grow, and
extending the framework to explicit state and path constraints,
treated separately in a state-constrained companion work. \RTadd{R51}{A further direction is to drop the
analytic-data assumption underlying the coefficient extension of
\S\ref{ssec:extension}: in settings where $H$ and its gradients are
available but a clean $(f,\ell)$ decomposition is not---black-box
oracles, learned Hamiltonians (neural value functions, Gaussian-process
means), mean-field reductions---a local-model surrogate replaces the
global PA majorant, the inequality $\bar H\ge H$ holds only at sample
points, and $\eta_{\mathrm{PA}}$ becomes a diagnostic rather than a
certified one-sided indicator; the design of adaptive PMP shooting in
this regime is left to future work.}

\MSdel{Outlook ``Next benchmarks within the unconstrained-state scope'' paragraph removed at R10 per user directive (the paragraph was a forward pointer to state-constrained variants and to a costate-optimization companion; the standalone version of this work carries no forward pointer beyond the unconstrained scope explicitly developed here).}
%% === end of Discussion_v2.tex ===

\section{Conclusion}\label{sec:conclusion}

%% === inlined from Conclusion_v2.tex (R98z-cc arxiv-flat) ===
% ------------------------------------------------------------------
% §6 Conclusion
% ------------------------------------------------------------------

We have presented a Pontryagin-based numerical solver for
deterministic optimal control problems in Bolza form. The method
constructs a concave piecewise-affine bundle surrogate of the
Hamiltonian, smooths it via a log-sum-exp formula to obtain a
Hamiltonian~$H_\delta$ that is $C^\infty$ and concave in the
costate and sufficiently differentiable in the state under the
standing assumptions, and solves the resulting
symplectic Euler TPBVP by damped Newton iteration. A unified
adaptive outer loop controls three distinct error sources---time
discretization, PA modeling, and smoothing bias---via computable
a~posteriori indicators.

\RTadd{R62-B1}{Four unconstrained-state benchmark problems---hypersensitive
dynamics, nonsmooth Hamiltonian switching, singular time localization,
and QP-oracle over-actuated attitude control allocation---illustrate
the solver's workflow:\ the adaptive mesh concentrates intervals near
control switches and singular features, the bundle grows only where
the PA gap is significant, and the smoothing parameter is annealed as
the solution stabilizes.\ A smooth LQR/Riccati calibration check is
reported separately in Appendix~\ref{apx:lqr}.}

This work is an algorithmic note for the unconstrained-state
regime. The solver is designed for moderate-dimensional problems where
the Hamiltonian is accessed via an oracle, the control set is compact,
and a reasonable initial guess is available.
\MSadd{Consistent with the value-first scope of~\S\ref{sec:introduction}, \RTadd{R53-M1}{the
reported output of the adaptive cycle is a numerical approximation $J_h$ to the value
function $U(0,x_0)$, supported by the three value-error indicators; the discrete state, costate,
and control are computational means. A post-convergence plug-and-integrate pass
(\S\ref{ssec:outer-loop}) recovers an admissible trajectory that serves as a primal
upper-bound witness on $U(0,x_0)$. Certified lower bounds and certified primal--dual
value-gaps require an independent lower verification function and are outside the scope
of this work.}}
\MSdel{State-constrained extensions of the present framework---including the viability-based tangential Hamiltonian, the relaxed discrete-step variant, the log-barrier interior-point smoothing, the constrained value-equivalence theorem, and the boundary-layer a~posteriori indicators---are developed in a forthcoming costate-optimization companion work that extends the present framework to explicit state constraints.}
Extensions to stochastic dynamics, hybrid systems, and rigorous convergence analysis are deferred to future work.

\RTadd{R53-M4}{}
%% === end of Conclusion_v2.tex ===

\appendix

%% === inlined from Appendix_LQR.tex (R98z-cc arxiv-flat) ===
% ------------------------------------------------------------------
% Appendix A — LQR/Riccati calibration consistency check
% Relocated from §4.1 to a calibration appendix at R53 Phase B (B3).
% R98z-pp-apply (2026-06-17): extended to two-mode comparison.
%   Mode-(i) closed-form smooth   — data from R98z-pp-extension (4-row trajectory).
%   Mode-(iii) PA-softmin         — data from R98z-pp        (11-row trajectory).
%   Three-curve control overlay (fig:lqr-control-comparison) added.
%   eq:lqr-Jref reference precision boosted to twelve decimals.
% ------------------------------------------------------------------

\section{LQR/Riccati calibration consistency check}\label{apx:lqr}

\RTadd{R53-B3}{This appendix records a smooth linear--quadratic
regulator (LQR) calibration run that exercises the solver on the
simplest, fully analytic Bolza data with a matrix-Riccati reference.}
\RTadd{R98z-pp-A}{Since LQR's Hamiltonian is smooth and explicit on the optimal trajectory (the box constraint is structurally inactive at the LQR optimum, as recorded below), the natural mode is the closed-form smooth Hamiltonian of \S\ref{ssec:pmp} (denoted \emph{mode-(i)} here for brevity), with the PA bundle and the $\delta$-smoothing layers structurally absent. We also exercise the universal PA-softmin scheme (\emph{mode-(iii)} of \S\ref{ssec:explicit-grad}) on the same LQR row as a compatibility test of the bundle-and-smoothing layers on a smooth benchmark, and report both modes side-by-side. The four main benchmarks of~\S\ref{sec:numerics} each engage at least one nontrivial adaptive branch; the LQR row is retained here as a consistency check rather than as a main benchmark.}

\paragraph{Problem statement.}
We solve the fixed-final-time Bolza problem
\[
\min_{u(\cdot)} J(u)\;=\;x(T)^\top Q_T\,x(T)
+\int_0^T\bigl(x(t)^\top Q\,x(t)+u(t)^\top R\,u(t)\bigr)\,dt
\]
subject to $\dot x=Ax+Bu$, $x(0)=x_0$, $t\in[0,T]$, with
\[
A=\begin{bmatrix}0&1\\0&0\end{bmatrix},\quad
B=\begin{bmatrix}0\\1\end{bmatrix},\quad
Q=Q_T=I_2,\quad
R=10^{-2},\quad
x_0=(1,0)^\top,\quad
T=1,
\]
and box control bounds $u\in[u_{\min},u_{\max}]=[-11,5]$.

\paragraph{Pontryagin system and reference.}
The Hamiltonian (minimization convention,
\S\ref{ssec:pmp}) is
$H(p,x,u)=p^\top(Ax+Bu)+x^\top Qx+u^\top Ru$. The unconstrained
stationarity condition $\partial_uH=B^\top p+2Ru=0$ gives the interior
minimizer $u^\star_{\mathrm{unc}}(t)=-\tfrac12 R^{-1}B^\top p(t)
=-50\,p_2(t)$, and the box-projected control is
$u^\star(t)=\Pi_{[-11,5]}(-50\,p_2(t))$. The adjoint system is
$\dot p=-A^\top p-2Qx$, $p(T)=2Q_Tx(T)$. The matrix Riccati equation
$-\dot P=A^\top P+PA-PBR^{-1}B^\top P+Q$, $P(T)=Q_T$, integrated
exactly via the Hamiltonian-matrix exponential,
provides the reference objective
\begin{equation}\label{eq:lqr-Jref}
J^\star \;=\; x_0^\top P(0)\,x_0 \;=\; 1.113742492015,
\end{equation}
\RTadd{R98z-pp-B}{accurate to twelve decimals (the legacy value $1.113744628163$ obtained from a $1000$-point trapezoidal trajectory-cost integration was off by approximately $2.14\!\times\!10^{-6}$; this appendix reports gaps against the precise value).}
The closed-loop reference trajectory satisfies
$\dot x=(A-BR^{-1}B^\top P(t))x$, $u^\star=-R^{-1}B^\top P(t)x$,
$p^\star=2P(t)x$. The optimal trajectory remains strictly interior to
the imposed box, so the bounds are inactive at the optimum and the
algorithmic restriction $u\in[-11,5]$ does not change the problem
being solved.

\RTadd{R98z-pp-C}{\paragraph{Two-mode calibration runs.}
The solver was run on the LQR row in two complementary modes, with the same problem data and starting from the same initial mesh $N_0=20$.
\emph{Mode-(i), closed-form smooth.} The reduced Hamiltonian is evaluated analytically on the unconstrained quadratic branch (the box is inactive at the optimum); $\delta=0$ structurally and the PA bundle is structurally absent. The time-mesh layer is the only active layer. Tolerance: $\varepsilon_{\mathrm{time}}=2.5\!\times\!10^{-5}$, marking threshold $s_{\mathrm{time}}=0.5$, RK4 on the augmented closed-loop state/cost equation. The solver terminates with a clean \texttt{STOP} after $4$ outer iterations on a final mesh of $N=46$ intervals (Tables~\ref{tab:lqr-mode-i-trajectory} and~\ref{tab:lqr-branches}), with computed objective $J_h=1.113745702206$ and absolute gap $|J_h-J^\star|=3.21\!\times\!10^{-6}$, comfortably below the per-interval threshold and three orders of magnitude below the gap that the next-best controller calibration would deliver under the same $\varepsilon_{\mathrm{time}}$.
\emph{Mode-(iii), PA-softmin compatibility test.} The full bundle-and-smoothing scheme is exercised on the same LQR row to confirm that those layers introduce no spurious artifact at the value level on a smooth benchmark. Initial values \Ksadd{R98z-bb-M0}{$M_0=9$}, $\delta_0=1.5\!\times\!10^{-1}$, tolerances $\varepsilon_{\mathrm{time}}=2\!\times\!10^{-2}$, $\varepsilon_{\mathrm{PA}}=\varepsilon_\delta=10^{-2}$, $\mathrm{max\_iters}=15$. The solver terminates with a clean \texttt{STOP} after $11$ outer iterations on a final mesh of $N=102$ intervals with $M=10$ planes and smoothing $\delta=4.69\!\times\!10^{-3}$ (Tables~\ref{tab:lqr-mode-iii-trajectory} and~\ref{tab:lqr-branches}); computed objective $J_h=1.122936958852$ and absolute gap $|J_h-J^\star|=9.19\!\times\!10^{-3}$. All three indicators are below their tolerances at termination.
Principal final values for both modes are listed side-by-side in Table~\ref{tab:lqr-summary}.}

\begin{table}[htbp]
\centering
\caption{Final principal results reported by the proposed solver on
 the LQR calibration run in both modes. ``Within tol'' refers to the per-interval
 threshold $\mathrm{tol}_{\mathrm{time}}^\star=\varepsilon_{\mathrm{time}}/N$.}
\label{tab:lqr-summary}
\begin{tabular}{lcc}
\toprule
Quantity & Mode-(i) closed-form smooth & Mode-(iii) PA-softmin \\
\midrule
Objective $J_h$ & $1.113745702206$ & $1.122936958852$ \\
Riccati reference $J^\star$ & $1.113742492015$ & $1.113742492015$ \\
Absolute gap $|J_h-J^\star|$ & $3.21\!\times\!10^{-6}$ & $9.19\!\times\!10^{-3}$ \\
Outer iterations $k_{\mathrm{end}}$ & $3$ & $10$ \\
Mesh intervals $N$ & $46$ & $102$ \\
Planes $M$ & --- (absent) & $10$ \\
Final smoothing $\delta$ & $0$ (absent) & $4.69\!\times\!10^{-3}$ \\
Time tolerance $\varepsilon_{\mathrm{time}}$ & $2.5\!\times\!10^{-5}$ & $2\!\times\!10^{-2}$ \\
$\eta_{\mathrm{time}}$ & $2.77\!\times\!10^{-7}$ & $1.20\!\times\!10^{-4}$ \\
$\mathrm{tol}_{\mathrm{time}}^\star$ & $5.43\!\times\!10^{-7}$ & $1.96\!\times\!10^{-4}$ \\
$\eta_{\mathrm{PA}}$ & $0$ (absent) & $2.32\!\times\!10^{-3}$ \\
$\eta_\delta$ & $0$ (absent) & $2.70\!\times\!10^{-3}$ \\
Final action & \texttt{STOP} & \texttt{STOP} \\
Status & accepted & accepted \\
\bottomrule
\end{tabular}
\end{table}

\RTadd{R98z-pp-D}{\begin{table}[htbp]
\centering
\caption{Per-iteration trajectory of the mode-(i) closed-form smooth LQR calibration run. Only the time layer is active: $\delta_k = 0$ and $M_k$ is structurally absent throughout. The controller takes three time-refine actions and then stops; the objective gap drops by about a decade per outer iteration.}
\label{tab:lqr-mode-i-trajectory}
\footnotesize
\begin{tabular}{cccccc}
\toprule
$k$ & $N_k$ & $J_h^{(k)}$ & $\eta_{\mathrm{time,max}}^{(k)}$ & $\mathrm{tol}_{\mathrm{time}}^{\star,(k)}$ & action \\
\midrule
$0$ & $20$ & $1.1140436330$ & $8.14\!\times\!10^{-5}$ & $1.25\!\times\!10^{-6}$ & \texttt{time\_refine} \\
$1$ & $33$ & $1.1137682995$ & $8.58\!\times\!10^{-5}$ & $7.58\!\times\!10^{-7}$ & \texttt{time\_refine} \\
$2$ & $41$ & $1.1137476824$ & $5.57\!\times\!10^{-6}$ & $6.10\!\times\!10^{-7}$ & \texttt{time\_refine} \\
$\mathbf 3$ & $\mathbf{46}$ & $\mathbf{1.1137457022}$ & $\mathbf{2.77\!\times\!10^{-7}}$ & $\mathbf{5.43\!\times\!10^{-7}}$ & $\mathbf{\texttt{STOP}}$ \\
\bottomrule
\end{tabular}
\end{table}

\begin{table}[htbp]
\centering
\caption{Per-iteration trajectory of the mode-(iii) PA-softmin LQR compatibility-test run. All three layers fire across the eleven outer iterations: five $\delta$-halvings, four time-refinements, one PA-enrichment, and the terminal \texttt{STOP}.}
\label{tab:lqr-mode-iii-trajectory}
\footnotesize
\begin{tabular}{ccccccc}
\toprule
$k$ & $\delta_k$ & $N_k$ & $M_k$ & $J_h^{(k)}$ & $\eta_{\mathrm{time,max}}^{(k)}$ & action \\
\midrule
$0$  & $0.15$  & $20$  & $9$  & $1.2476260091$ & $2.44\!\times\!10^{-2}$ & \texttt{time\_refine} \\
$1$  & $0.15$  & $37$  & $9$  & $1.1997785116$ & $7.15\!\times\!10^{-3}$ & \texttt{delta\_halve} \\
$2$  & $7.50\!\times\!10^{-2}$ & $37$  & $9$  & $1.1788683149$ & $7.54\!\times\!10^{-3}$ & \texttt{delta\_halve} \\
$3$  & $3.75\!\times\!10^{-2}$ & $37$  & $9$  & $1.1721298618$ & $7.66\!\times\!10^{-3}$ & \texttt{time\_refine} \\
$4$  & $3.75\!\times\!10^{-2}$ & $62$  & $9$  & $1.1458541654$ & $1.92\!\times\!10^{-3}$ & \texttt{delta\_halve} \\
$5$  & $1.88\!\times\!10^{-2}$ & $62$  & $9$  & $1.1426438289$ & $1.92\!\times\!10^{-3}$ & \texttt{time\_refine} \\
$6$  & $1.88\!\times\!10^{-2}$ & $86$  & $9$  & $1.1304938669$ & $4.80\!\times\!10^{-4}$ & \texttt{delta\_halve} \\
$7$  & $9.38\!\times\!10^{-3}$ & $86$  & $9$  & $1.1277586662$ & $4.79\!\times\!10^{-4}$ & \texttt{delta\_halve} \\
$8$  & $4.69\!\times\!10^{-3}$ & $86$  & $9$  & $1.1259248592$ & $4.79\!\times\!10^{-4}$ & \texttt{PA\_enrich} \\
$9$  & $4.69\!\times\!10^{-3}$ & $86$  & $10$ & $1.1274603533$ & $4.79\!\times\!10^{-4}$ & \texttt{time\_refine} \\
$\mathbf{10}$ & $\mathbf{4.69\!\times\!10^{-3}}$ & $\mathbf{102}$ & $\mathbf{10}$ & $\mathbf{1.1229369589}$ & $\mathbf{1.20\!\times\!10^{-4}}$ & $\mathbf{\texttt{STOP}}$ \\
\bottomrule
\end{tabular}
\end{table}}

\RTadd{R98z-pp-E}{\begin{table}[htbp]
\centering
\caption{Branch-firing diagnostic for the mode-(i) LQR run (left) and the mode-(iii) compatibility test (right). The mode-(i) tally collapses to time-refine and \texttt{STOP} only, reflecting the structurally inactive PA and $\delta$ layers; the mode-(iii) tally exhibits the expected mixed-layer activity on the same LQR row.}
\label{tab:lqr-branches}
\begin{tabular}{lc@{\hspace{1.5cm}}lc}
\toprule
\multicolumn{2}{c}{Mode-(i) closed-form smooth} & \multicolumn{2}{c}{Mode-(iii) PA-softmin} \\
\midrule
Outer-iteration outcome & Count & Outer-iteration outcome & Count \\
\midrule
\texttt{delta\_halve\_alone}     & $0$ & \texttt{delta\_halve\_alone}     & $5$ \\
\texttt{time\_refine\_alone}     & $3$ & \texttt{time\_refine\_alone}     & $4$ \\
\texttt{PA\_enrich\_alone}       & $0$ & \texttt{PA\_enrich\_alone}       & $1$ \\
compound (any)                   & $0$ & compound (any)                   & $0$ \\
\texttt{STOP}                    & $1$ & \texttt{STOP}                    & $1$ \\
\midrule
Total                            & $\mathbf{4}$ & Total                            & $\mathbf{11}$ \\
\bottomrule
\end{tabular}
\end{table}}

\RTadd{R98z-pp-F}{\begin{figure}[htbp]
\centering
\includegraphics[width=0.92\textwidth]{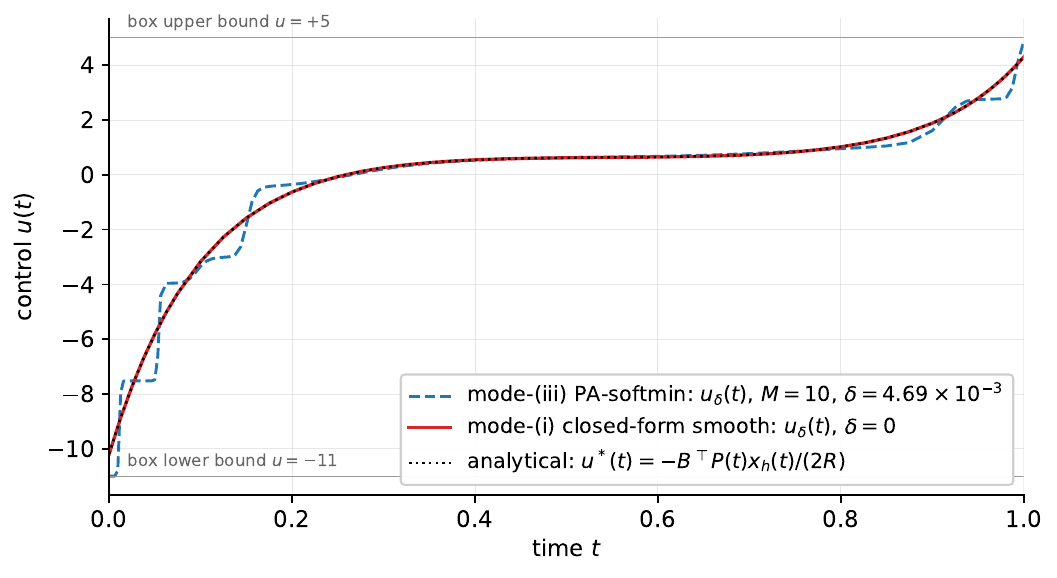}
\caption{Effective control on the final converged trajectory: mode-(i) closed-form smooth $u_\delta(t)$ (red solid, $47$ nodes on the final $N=46$ mesh), mode-(iii) PA-softmin $u_\delta(t)$ (blue dashed, $103$ nodes on the final $N=102$ mesh), and analytical Riccati feedback $u^\star(t)=-R^{-1}B^\top P(t)\,x_h(t)$ (black dotted, evaluated on the mode-(i) mesh). The full-cost convention applies throughout (PMP costate $p=2P(t)x$, effective control $u_\delta=-p_2/(2R)$). The mode-(i) curve and the analytical $u^\star$ overlap to round-off (maximum nodal deviation $1.51\!\times\!10^{-11}$); the mode-(iii) curve shows a visible deviation, peaking at the early-time saturation $u_\delta(0)=-11$ at the box lower bound against the analytical $u^\star(0)=-10.18$ (maximum nodal deviation $1.66$). The mode-(iii) deviation visualizes the residual smoothing budget of the PA-softmin scheme at the converged $\delta=4.69\!\times\!10^{-3}$, exactly the regime in which the small-$\delta$ Jacobian-conditioning caveat of \S\ref{ssec:limitations}, Remark~\ref{rem:H_delta_hessian} is salient.}
\label{fig:lqr-control-comparison}
\end{figure}}

\RTadd{R98z-pp-G}{\paragraph{Discussion: mode-(i) and §\ref{ssec:limitations} Remark~\ref{rem:H_delta_hessian}.} On LQR's smooth Hamiltonian, the controller correctly stays in mode-(i) with $\delta=0$ and the Jacobian-conditioning caveat of \S\ref{ssec:limitations}, Remark~\ref{rem:H_delta_hessian} does not arise: the run takes only three time-refinement actions, the value gap drops by approximately a decade per outer iteration to $|J_h-J^\star|=3.21\!\times\!10^{-6}$ (Table~\ref{tab:lqr-mode-i-trajectory}), and the effective control $u_\delta(t)$ tracks the analytical $u^\star(t)$ to round-off everywhere on $[0,T]$ (red and black curves overlap in Figure~\ref{fig:lqr-control-comparison}). The mode-(i) calibration confirms that the smooth-Hamiltonian arm of the framework reaches arbitrary accuracy on this benchmark, limited only by the RK4 budget on the time mesh.

\paragraph{Discussion: mode-(iii) compatibility test.} The PA-softmin compatibility test on the same LQR row stops cleanly with all three indicators below their respective tolerances after eleven outer iterations (Table~\ref{tab:lqr-mode-iii-trajectory}): five $\delta$-halvings drive the smoothing parameter from $\delta_0=1.5\!\times\!10^{-1}$ to $\delta=4.69\!\times\!10^{-3}$, four time-refinements grow the mesh from $N_0=20$ to $N=102$, and a single PA enrichment grows the bundle from $M=9$ to $M=10$. The converged $\delta=4.69\!\times\!10^{-3}$ sits in the small-$\delta$ regime where the Jacobian-conditioning caveat of Remark~\ref{rem:H_delta_hessian} bites; the corresponding effective control $u_\delta(t)$ deviates visibly from $u^\star(t)$ (blue dashed vs.\ black dotted in Figure~\ref{fig:lqr-control-comparison}), most strikingly at $t=0$ where the smoothed control saturates at the box lower bound $u_\delta(0)=-11$ even though the analytical $u^\star(0)=-10.18$ sits strictly interior to the box. The deviation does not propagate catastrophically to the objective ($|J_h-J^\star|=9.19\!\times\!10^{-3}$, three orders of magnitude looser than mode-(i) but consistent with the chosen $\varepsilon_{\mathrm{time}}=2\!\times\!10^{-2}$), and the controller correctly STOPs at the indicator-level rather than chasing the artifact further. The compatibility test is therefore consistent with the framework: each mode is calibrated to its own appropriate tolerance regime, demonstrating mode-polymorphism on a smooth benchmark while making explicit, via Figure~\ref{fig:lqr-control-comparison}, the residual smoothing-and-conditioning budget that Remark~\ref{rem:H_delta_hessian} flags.}

\begin{figure}[htbp]
\centering
\begin{minipage}{0.48\textwidth}\centering
 \includegraphics[width=\textwidth]{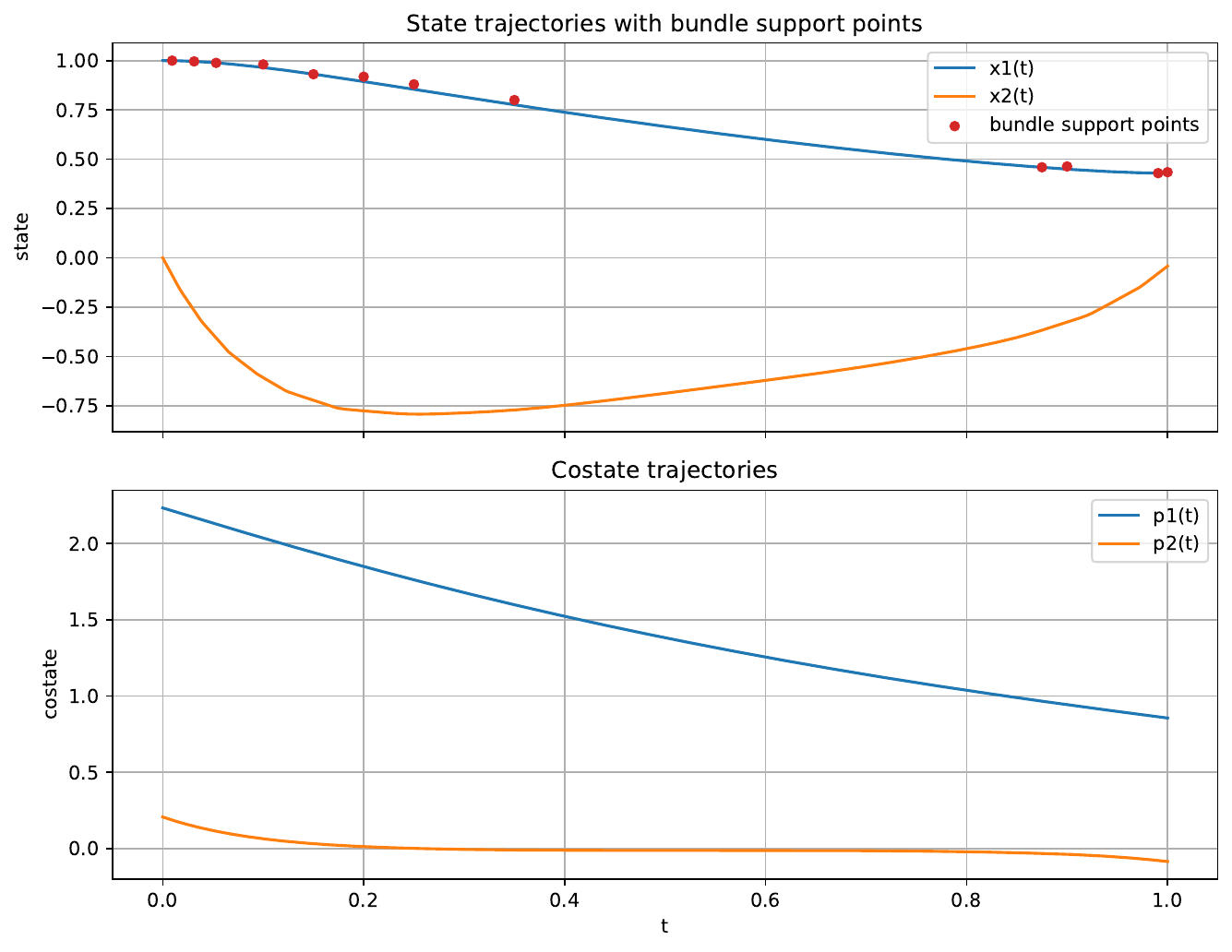}
\end{minipage}\hfill
\begin{minipage}{0.48\textwidth}\centering
 \includegraphics[width=\textwidth]{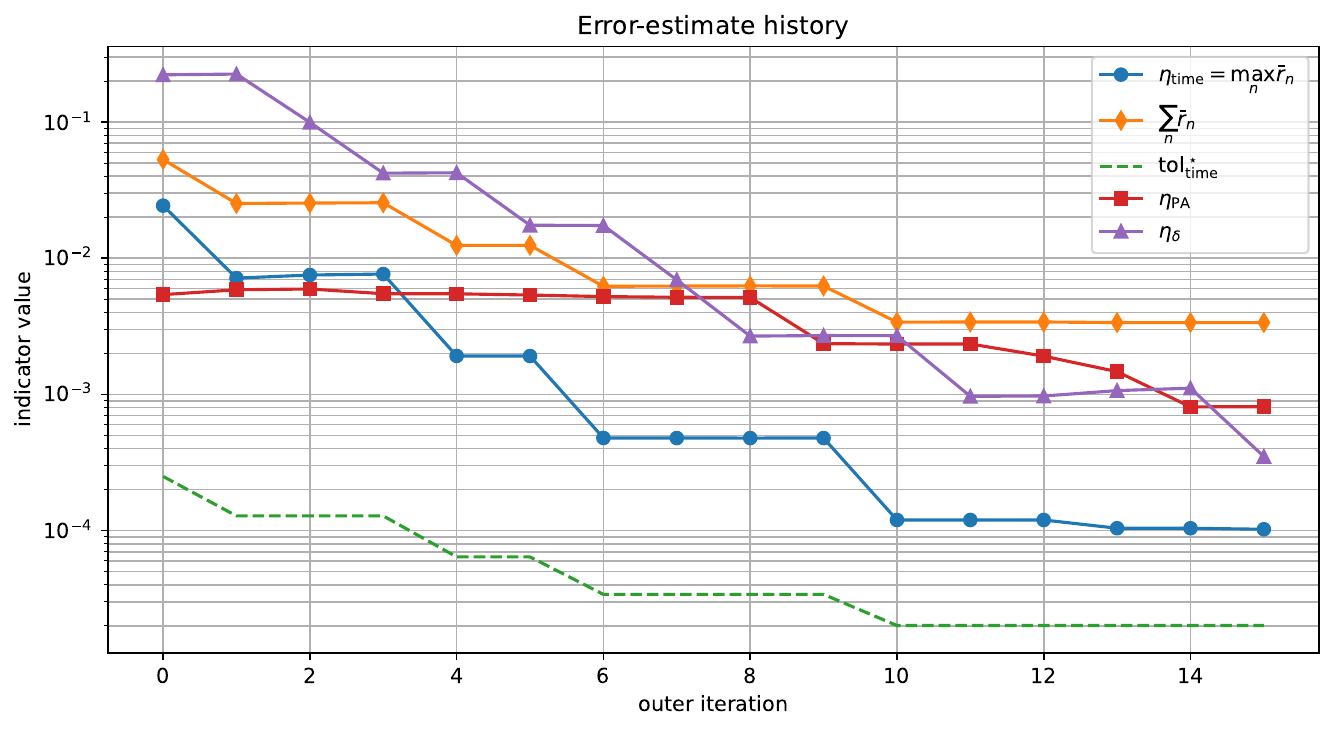}
\end{minipage}
\caption{LQR mode-(iii) PA-softmin calibration run. Left:\ state and costate at the final
 iterate ($N=102$, $M=10$, $\delta=4.69\!\times\!10^{-3}$).
 Right:\ history of the global indicators along the adaptive outer
 loop;\ the figure distinguishes the stopping quantity
 $\eta_{\mathrm{time}}=\max_n\bar r_n$ from the companion sum
 $\sum_n\bar r_n$.}
\label{fig:lqr-state-costate}
\end{figure}

Figure~\ref{fig:lqr-state-costate} shows the state and costate at the
final iterate of the mode-(iii) compatibility run, and the indicator history along the outer loop:\ the
bundle indicator drops quickly to a nontrivial PA-bundle plateau, the
smoothing indicator contracts geometrically by the $\delta$-halving
steps, and the time indicator decreases under successive mesh
refinements until all three pass their thresholds simultaneously.

%% === end of Appendix_LQR.tex ===

% R98z-bb (user comment): duplicate AI-editing-tools acknowledgement removed; the first occurrence as a title-page footnote at the top of the manuscript (\Rzzedit{S1}) is kept as the canonical disclosure.

\bibliographystyle{siamplain}
\bibliography{refs_pmp}

\end{document}